\documentclass[USenglish]{article}	
\usepackage[utf8]{inputenc}
\usepackage[big,online]{dgruyter}
\usepackage{verbatim}
\usepackage{lmodern} 
\usepackage{microtype}
\usepackage[numbers,square,sort&compress]{natbib}
\usepackage{amsmath,amsfonts,amssymb,amsthm,bm}
\usepackage{float,graphicx,caption,subcaption}
\usepackage{newtxtext,newtxmath}
\usepackage{algpseudocode}
\usepackage{filecontents}
\usepackage{amsmath}
\usepackage[normalem]{ulem}
\usepackage{xcolor}
\usepackage{subcaption}
\theoremstyle{dgthm}

\theoremstyle{dgdef}

\newtheorem{example}{Example}

\newcommand{\RE}{\mathbb{R}}

\newcommand{\ceil}[1]{\left\lceil #1 \right\rceil}

\newcommand{\re}[1]{\begin{color}{blue}#1\end{color}}
\newcommand\BGN{\begin{color}{blue}}
\newcommand\ND{\end{color}}

\begin{document}

% %%%--------------------------------------------%%%
% 	\articletype{Research Article}
% 	\received{Month	DD, YYYY}
% 	\revised{Month	DD, YYYY}
%   \accepted{Month	DD, YYYY}
%   \journalname{De~Gruyter~Journal}
%   \journalyear{YYYY}
%   \journalvolume{XX}
%   \journalissue{X}
%   \startpage{1}
%   \aop
%   \DOI{10.1515/sample-YYYY-XXXX}
%%%--------------------------------------------%%%

\title{A data-driven method for parametric PDE eigenvalue problems using Gaussian Process with different covariance functions}
\runningtitle{A data-driven method for PDE eigenvalue problems using GPR}
%\subtitle{Insert subtitle if needed}

\author[1]{Moataz Alghamdi}
%\ use * to mark the author as the corresponding author
\author[2]{Fleurianne Bertrand}
\author[1,3]{Daniele Boffi} 
\author*[1,4]{Abdul Halim} 
\runningauthor{M. Alghamdi et al.}
\affil[1]{\protect\raggedright 
Applied Mathematics and Computational Sciences (AMCS)\\
King Abdullah University of Science and Technology\\
Thuwal,Kingdom of Saudi Arabia \\}
%email: moataz.alghamdi@kaust.edu.sa}
\affil[2]{\protect\raggedright Faculty of Mathematics, TU Chemnitz, Chemnitz, Germany\\
email: fleurianne.bertrand@math.tu-chemnitz.de}
\affil[3]{\protect\raggedright Department of Mathematics ``F. Casorati''\\
University of Pavia\\
via Ferrata 1, 27100 Pavia, Italy}
\affil[4]{\protect\raggedright Department of Mathematics \\
Memari College, West Bengal, India\\
email: abdul.math91@gmail.com}
%\communicated{...}
%\dedication{...}
	
\abstract{We use a Gaussian Process Regression (GPR) strategy that was recently developed \cite{Guo2019,MengwuGJan,dd_gpr} to analyze different types of curves that are commonly encountered in parametric eigenvalue problems. We employ an offline-online decomposition method. In the offline phase, we generate the basis of the reduced space by applying the proper orthogonal decomposition (POD) method on a collection of pre-computed, full-order snapshots at a chosen set of parameters. Then, we generate our GPR model using four different Mat\'{e}rn covariance functions. In the online phase, we use this model to predict both eigenvalues and eigenvectors at new parameters. We then illustrate how the choice of each covariance function influences the performance of GPR. Furthermore, we discuss the connection between Gaussian Process Regression and spline methods and compare the performance of the GPR method against linear and cubic spline methods. We show that GPR outperforms other methods for functions with a certain regularity.}


\maketitle

\section{Introduction} 

PDE eigenvalue problems (EVPs) are an important class of problems frequently arising in science and engineering. An important subclass of such EVPs are those dependent on a set of parameters that may come from the material properties, geometric domain, or initial or boundary conditions. They are commonly known as parametric eigenvalue problems (pEVPs). They are distinguished by the fact that their solutions are multi-output, parameter-dependent eigenvectors and eigenvalues, allowing for the possibility of intersections. Given their dependence on parameters, numerical approaches used to solve such eigenvalue problems require in general solving them multiple times at some chosen parameters. One thus needs to discretize the domain, solve for the eigenvalues and eigenvectors, and assemble the corresponding matrices at each such parameter, making this process computationally very expensive. 

Reduced-order modeling techniques have been traditionally developed to overcome such computational difficulties. Model order reduction (MOR) techniques, in particular, allow for fast but accurate computations. We are interested in MOR techniques which use an offline-online decomposition. In this class of techniques, one generates in the offline phase a collection of snapshots of eigensolutions evaluated at a chosen set of parameters. Such snapshots are generally generated in two different ways (or a combination thereof): using a) adaptive greedy algorithms such as in~\cite{alghamdi2022greedy} or b) the proper orthogonal decomposition (POD) technique. Then, a surrogate is constructed, which is used to approximate the eigensolution at any intermediate parameter. This is done in the online phase. One such popular method is the reduced basis (RB) method which has been developed to approximate one \cite{Machiels-Maday-Oliveira-Patera-Rovas-2000} or more \cite{Vallaghe-Huynh-Knezevic-Nguyen-Patera-2015,Pau07a,Pau07b,Horgeretal17} eigenvalues. Stochastic Galerkin \cite{Hakula-Kaarnioja-Laaksonen-2015} and stochastic collocation methods \cite{Andreev-Schwab-2012} have also been proposed for stochastic EVPs. However, in comparison to progress made for source problems (see for example~\cite{Quarteronietal16, QuarteroniRozza07}), eigenvalue problems (for more, see ~\cite{Buchanetal13, GermanRagusa19}) are, in general, much less studied despite their importance.

Furthermore, a greedy RB method has also been developed to deal with affine and non-affine parametric eigenvalue problems \cite{Fumagalli-Manzoni-Prolini-Verani-2016}, focusing only on a single eigenpair. More recently, a data-driven approach~\cite{dd_gpr} that uses Gaussian Process Regression (GPR) has been proposed to deal with both affine and non-affine problems, an approach inspired by the following previous works focused on source problems \cite{MengwuGJan,Guo2019}. 

% In the context of reduced-order models, GPR is used as a surrogate for the parameter-to-eigensolution map in the online phase. 
In this paper, we explore further applying the Gaussian Process Regression (GPR) in parametric eigenvalue problems. GPR is a statistical approach that is normally used in fields such as supervised machine learning. However, its connection to spline methods, which are commonly used in PDE problems, has been known for close to fifty years \cite{Kimeldorf1970}. Nevertheless, such a connection has not been fully taken advantage of in solving eigenvalue problems and source problems. Our paper thus aims to add to previously mentioned contributions to improve our understanding of this connection.

This paper is structured as follows. We introduce parametric eigenvalue problems in Section~\ref{se:ps}. Then, we introduce GPR and discuss its connection to spline methods in~Sections \ref{se:GPR} and \ref{se:GPRSplines}. We also make comparisons in cases where uniform and non-uniform grids are used. We emphasize how the statistical nature of GPR can be utilized to allow for better performance in both cases. Indeed, GPR is representative of statistical approaches used in fields such as machine learning. Spline methods, on the other hand, are the most popular choice for PDE problems. Bridging this rather large gap then will allow for utilizing applied or theoretical techniques developed for PDE problems in machine learning and vice versa.

In addition, we investigate in Section~\ref{se:Num} issues related to applying GPR to three different pEVPs. As will be shown, the performance of GPR is heavily influenced by the choice of the covariance function. We thus focus on analyzing the behavior of the solutions for multiple different covariance functions. In particular, as intersections between eigenvalues are possible, the solutions to be approximated might not be smooth. Consequently, we thoroughly discuss the connection between the regularity of the covariance function and the eigensolutions.
\section{Parametric eigenvalue problems}
\label{se:ps}
%%%%%%%%%%%%%%%%%%%%%%%%%%%%%%%%%%%%%%%%%%%%%%%%%%%%%%%%%%%%

Let $\mathcal{M} \subset \RE^d$ be closed and bounded. In this paper, we consider only the case $d \leq 2$. For each $\bm{\mu} \in \mathcal{M}$, we define two symmetric and bilinear forms
\[
\aligned
&a(\cdot,\cdot;\bm{\mu}):V\times V\to\RE,\\
&b(\cdot,\cdot;\bm{\mu}):H\times H\to\RE.
\endaligned
\]
where we assume the existence of the following Hilbert triplet
$V\subset H\simeq H'\subset V'$.

We seek to solve the following problem: for all $\bm{\mu} \in \mathcal{M}$, find real eigenvalues $\lambda(\bm{\mu})$ and non-vanishing eigenfunctions $u(\bm{\mu})\in V$ such that:
\begin{equation}  
\label{eq:main}
a(u,v;\bm{\mu})=\lambda(\bm{\mu})b(u,v;\bm{\mu})\quad\forall v\in V.
\end{equation}

To ensure that all eigenvalues %$\{\lambda_j(\mu)\}_{j=1}^\infty$
have finite-dimensional eigenspaces, we make the following assumptions so that we have only a compact solution operator. First, we require that V is compact in H. Second, for all $\bm{\mu} \in \mathcal{M}$, the bilinear form $a(\cdot,\cdot;\bm{\mu})$ is elliptic in $V$, while $b(\cdot,\cdot;\bm{\mu})$ is equivalent to the inner product in $H$. %\quad\forall\mu\in \mathcal{M}$

For each $\bm{\mu} \in \mathcal{M}$, we compute our high-fidelity eigensolutions, namely the eigenvalues and their associated eigenvectors. We do this by discretizing using finite elements. In other words, we consider a finite-dimensional subspace $V_h\subset V$ of dimension $n_h$. 
Then, for all $\bm{\mu} \in \mathcal{M}$, we consider the eigensolutions for the following generalized eigenvalue problem: find the real eigenvalues $\lambda_h(\bm{\mu})$ and non-vanishing eigenfunctions $u_h(\bm{\mu})\in V_h$ such that

\begin{equation}
    \label{eq:mainh}
    a(u_h,v;\bm{\mu})=\lambda_h(\mu)b(u_h,v;\bm{\mu})\quad\forall v\in V_h.
\end{equation}

% \subsection{\sout{Reduced Order Modelling}}

% \sout{Reduced-order modeling is commonly used to solve parametric eigenvalue problems. In reduced-order modeling, one has to solve the discretized problem for a few selected parameters. Then, using such solutions, a reduced space is formed. Solving the problem for any new parameter is then sought in the reduced space. The most common approaches for generating such a reduced space are projection-based approaches.}

% \sout{The success of projection-based reduced order modeling is based on the assumption that the parametric eigenvalue problem is affine-parameter dependent. The variational form for the affine-parameter dependent problems is written as:}
% \begin{align}
%     \sum\limits_{q=1}^{n_a}\theta_a^q
%    (\bm{\mu})a_q(w,v)=\lambda(\bm{\mu}) \sum\limits_{q=1}^{n_b}\theta_b^q
%    (\bm{\mu})b_q(w,v),\quad \forall w,v \in V,
% \end{align}
% \sout{where $a_q(\cdot,\cdot),b_q(\cdot,\cdot):V\times V \to \mathbb{R}$ are bilinear forms. If the problem is not affine, then such an approach is not applicable. For non-affine problems, a data-driven model has been shown to be more suitable for solving parametric problems~\cite{dd_gpr}.}

\subsection{Solutions as data sets}

We generate our data set by numerically solving the parametric eigenvalue problem for the first $n_l$ eigenvalues at $n_s$ pre-selected $d$-dimensional parameters $ \bm{\mu}_1,\bm{\mu}_2,\dots,\bm{\mu}_{n_s} \in \mathcal{M}$. For convenience, we label this parameter set $\mathcal{M}_{tr}$. Such $d$-dimensional parameters can be selected using different methods such as uniform tensorial sampling, Latin hypercube sampling, or sampling on Smolyak sparse grids.

Let $j = 1,\dots,n_l$. Our data is then a collection of $n_l$ matrices $\bm{W}_j$. Each such matrix consists of the eigenvectors corresponding to the $j$-th lowest eigenvalue at each parameter $\bm{\mu} \in \mathcal{M}_{tr}$. In other words, $\bm{W}_j = \{(\lambda_j(\bm{\mu}_i),\bm{u}_{h,j}(\bm{\mu}_i)) | i = {1, \cdots , n_s}\}$ where $\lambda_j(\bm{\mu}_i)$ is a scalar and $\bm{u}_{h,j}(\bm{\mu}_i)$ is an $n_h \times 1$ vector. $\bm{W}_j$ can thus be represented as an $(n_h+1) \times n_s$ matrix where the first row corresponds exclusively to the eigenvalues $\lambda_j(\bm{\mu}_1), \dots, \lambda_j(\bm{\mu}_{n_s})$ while the rest of the $n_h+1$ rows correspond to the eigenfunctions $\bm{u}_{h,j}(\bm{\mu}_1), \dots,\bm{u}_{h,j}(\bm{\mu}_{n_s})$. Indeed, the latter has been traditionally used to define the snapshot matrix $\bm{S}_j = [\bm{u}_{h,j}(\bm{\mu}_1)|\dots|\bm{u}_{h,j}(\bm{\mu}_{n_s})]$ which is of size $n_h \times n_s$. 

Next, we construct a reduced basis space for the column space associated with each snapshot matrix $\bm{S_j}$, namely Col($\bm{S_j}$). To do so, we use the principal orthogonal decomposition technique (POD) which generates a low rank (rank = $L_j$) approximation for each such column space, Col(\bm{$S_j$}). We choose here the POD technique because it allows us to perform model order reduction for both linear and nonlinear eigenvalue problems. This is not the case when using other approaches such as Greedy algorithms which require the use of estimators/indicators that are not always known~\cite{Guo2019}.

To perform the reduction, we need to find the POD basis $\zeta_{1,j},\dots, \zeta_{L_j,j}$ associated with each column space Col($\bm{S}_j$). To do so, we perform a singular value decomposition on each $\bm{S}_j$. This allows us to write $\bm{S}_j$ as $\bm{S}_j = \bm{U}_j \bm{\Sigma}_j \bm{Z}_j^T$. Here, $\bm{U}_j = [\phi_{1,j},\dots,\phi_{n_h,j}] \in \RE^{n_h} \times \RE^{n_h}$  and $\bm{Z}_j = [\psi_{1,j},\dots,\psi_{n_s,j}] \in \RE^{n_s} \times \RE^{n_s}$ are orthogonal matrices. Furthermore, $\bm{\Sigma}_j = diag\{\sigma_{1,j}, \sigma_{2,j}, \dots, \sigma_{n_s,j}\}$ such that $\sigma_1 \geq \sigma_2 \geq \dots \geq \sigma_{n_s} $. 

% \sout{The reduced basis $\zeta_{1,j},\dots, \zeta_{L,j}$ is to be generated as follows. One can write $\bm{S}_j \phi_{i,j}$ = $\sigma_j \psi_{i,j}$ and $\bm{S}_j^{T} \psi_{i,j}$ = $\sigma_j \phi_{i,j}$ where $i = 1,\dots,L$ Thus, we have $\bm{S}_j \bm{S}_j^{T} \psi_{j,i}$ = $\sigma_j^2 \psi_{j,i}$ and $\bm{S}_j^{T}\bm{S}_j \phi_{j,i}$ = $\sigma_j^2 \phi_{j,i}$. Indeed, $\sigma_j^2$ are eigenvalues for both $\bm{S}^T \bm{S}$ and $\bm{S} \bm{S}^T$. The matrix $\bm{C}_j = \bm{S_j}^T \bm{S_j}$ is called correlation matrix, and the matrix $\bm{S} \bm{S}^T$ is a projection matrix}.

Now, let $\mathbb{V}_j  = \{ \bm{V}_j \in \bm{R}^{n_h \times L_j}: \bm{V}_j^{T}\bm{V}_j = I_{L_j} \}$ where $L_j \leq n_s$. The projection error $\sum_{i=1}^{n_s}\|u_{h,j}-\bm{V}_j\bm{V}_j^{T} u_{h,j}\|_{L^2(R^{n_h})} = \sum_{i=L_j+1}^{n_s} \sigma_{i,j}^2$ is then to be minimized. Note that $\bm{V}_j \bm{V}_j^T$ is a projection matrix. By the Schmidt–Eckart–Young theorem~\cite{Eckart1936}, the error is minimized by letting $\bm{V}_j$ contain the first $j$ columns of $\bm{U}_j$. One can then define an associated error tolerance $\epsilon_{POD}$ satisfying $\frac{\sum_{i=1}^{L_j} \sigma_{i,j}^2}{\sum_{i=1}^{n_s} \sigma_{i,j}^2} > 1-\epsilon^2_{POD}$.
Throughout this paper, we use $\epsilon_{POD} = 10^{-8}$. For a detailed discussion on POD, one may also refer to~\cite{Quarteronietal16, WithGopal}.

% { \color{red} This is known to be the optimal approximation space in the $L_2(\mathcal{M})$ sense []. By the Schmidt–Eckart–Young theorem [], such a basis minimizes the projection error of the snapshots among all other orthogonal bases in [rephrase]. Such error then is given by $\Sigma_{N+1}^{n_s} \sigma_i^2$ where $n_s >> N_j$. Furthermore, one can define the rank $N$ by setting a tolerance $\epsilon$ such that 
% $\frac{\Sigma_i^{N} \sigma_i^2}{\Sigma_i^{n_s} \sigma_i^2} \leq 1-\epsilon$. For a detailed discussion on POD, one may refer to~\cite{Quarteronietal16, WithGopal}.}

% \sout{Given $\bm{V}$, we can now construct the reduced snapshot matrix as follows. $\bm{S}_j = [V_j^T u_{h,j}(\bm{\mu}_1)|\dots| V_j^T u_{h,j}(\bm{\mu}_{n_s})]$. One thus can define the reduced version of $W_j$, $\hat{W}_j$, by replacing the snapshot matrix with the reduced snapshot matrix. Thus, $\hat{W}_j$ is of size $ (L+1) \times n_s$. Thus, the output variables here are $\hat{y}_j(\bm{\mu}_i)$ where $j = 1,\dots,L+1$.}

%The reduced (POD) basis are thus $\zeta_{1,j} = V_j^Tu_{h,j}(\bm{\mu}_1), \dots, \zeta_{L,j} = V_j^Tu_{h,j}(\bm{\mu}_{n_s})$

% \sout{Let $\mathbb{V}$ be the transform matrix between the basis of reduced space and the basis of the finite-dimensional space $V_h$.}

Following the proposed strategy in \cite{Guo2019,dd_gpr}, we use, in the offline phase, a collection of Bayesian linear regression models $\hat{\pi}_j: \mathcal{M} \to \mathbb{R}^{L_j+1}$ to predict the $j^{\text{th}}$ eigensolution at any parameter $\bm{\mu}_{i}^*$. To do so, we define two collections of data sets each containing $n_l$ training sets as follows.

$$\mathcal{Q}_{tr,j} = \{ \{\bm{\mu}, \bm{V}_j^{T}\bm{u}_{h,j}(\bm{\mu})\}: \bm{\mu} \in \mathcal{M}_{tr} \}$$
$$\mathcal{T}_{tr,j} = \{ \{\bm{\mu}, \lambda_j(\bm{\mu})\}: \bm{\mu} \in \mathcal{M}_{tr} \}$$
where $j = 1,\dots,n_l$. These will be the data sets used to train the model. For convenience, we call them the training data sets. {Similar data sets, namely the test data sets, will also be generated. These data sets will be used only for evaluating the performance of the GPR model.} They will be introduced in Section~\ref{se:Num} as needed. The machinery for generating such regression models will be discussed in Section~\ref{se:GPR}.
%Here, we use $\hat{W}$ as our database. % using  the reduced solutions in the reduced space of dimension $N_j$. %The construction of the regression model is done in the online stage. 

It is important to note here that there are works such as ~\cite{WangB1,Chen_2023} that discuss the possibility of using Bayesian linear modeling for data sets which are not assumed to be uncorrelated. In our case, we assume that there is no correlation between the $n_l$ data sets. Furthermore, we generate the POD basis as discussed previously which forms an orthonormal set. We, however, plan to investigate any possible advantages of using such approaches in future works.

% \sout{It is important to emphasize that we treat each row of $Q_{tr,j}$ as an independent data set. Thus, each pair $(\bm{\mu}_i,y_j(\bm{\mu}_i))$ is treated as an observation.}

% \sout{This approach is not the only possible one. For example, other works such as consider cases where all data sets are considered simultaneously. We plan to investigate this possibility in future works.}

% \sout{For now, we are left with only $N_j+1$ independent data sets, each containing $n_s$ observations. Furthermore, each data set is split into two data sets: a training set and a test set.} %For convenience, we refer to a training parameter as $\bm{\mu_i}$ and for a test parameter as $\bm{\mu_i}^{\star}$.

\section{Gaussian Process Regression} \label{se:GPR}

For each data set, we seek to train a machine-learning model using the training set. A well-trained model will make predictions that well approximate the test set. Toward this goal, we adopt a Bayesian approach known as Gaussian Process Regression (GPR). 

Bayesian methods \cite{BayesianM,Frazier2018} have been developed to combine prior information with available data to generate an updated or posterior model that considers the data. This posterior distribution can also always in turn be used as a prior distribution when new data becomes available. This framework then fits adaptive methods. For example, such Bayesian models can be used as surrogate models in the online phase of adaptive reduced-order models.

Traditionally, splines \cite{WahbaS,SilvermanGreen,WahbaCarven,5717055} have been used to make such predictions. Splines are known to produce accurate estimates for uniform, dense grids. However, they can lead to erroneous results in non-uniform or sparse grids. Unless new parameters are uniformly added to such grids, this issue cannot be remedied. This is problematic since adaptive schemes can generate non-uniform grids with sparse regions such as in~\cite{alghamdi2022greedy}. Bayesian modeling then allows for framing this deficiency in data as uncertainty. One can then exploit statistical relationships between the data points to enhance the predictive power of the model.

For the description of GPR, we rely on the detailed discussion given in \cite{gp_ml}. However, we present only the parts relevant to our purposes.

A Gaussian probability distribution is a continuous probability distribution for a real-valued random variable. Its probability density function $p(\bm{x})$ is given by:
\begin{equation*}\label{eq:gaussian_pdf}
\begin{split}
p(\bm{x}) = \frac{1}{(2\pi)^{k/2} |\bm{\Sigma}|^{1/2}} \exp (- \frac{1}{2} \left(\bm{x}-\bm{\mu}\right)^{T} \bm{\Sigma}^{-1} \left(\bm{x}-\bm{\mu} \right)),
\end{split}    
\end{equation*}
where $\bm{\mu} \in \RE^{k}$ is the mean vector and $\bm{\Sigma} \in \RE^{k\times k}$ is the variance of the distribution over the $k$-dimensional vector $\bm{x}$.
A Gaussian process (GP) is a generalization of the Gaussian probability distribution. A GP is a distribution over functions. In particular, a GP specifies a prior over functions, allowing us to control what functions we sample to predict the true solution (e.g., functions with a certain regularity). A data set is interpreted as a finite-dimensional (of dim $n_s$) realization of the chosen GP. Thus, such an $n_s$-dimensional realization follows a multivariate Gaussian distribution by definition. Using Bayes' theorem, one can use such a data set to update the prior GP to get a posterior GP, which is informed by the available data.

A GP is completely characterized by a mean function $m(\bm{\mu})$ and a covariance or kernel function $k(\bm{\mu},\bm{\mu}')$. Thus, we denote a GP by $\mathcal{GP}(m(\bm{\mu}),k(\bm{\mu},\bm{\mu}'))$. This makes utilizing this mathematical structure much more convenient. In particular, the choice of the mean function and the covariance function will dictate what type of functions can be sampled. Given this mathematical structure, we define our problem as follows:

Let $\bm{M} = [\bm{\mu}_1| \dots| \bm{\mu}_{n_s}]$ and $\bm{y} = ({y}_1,\dots,{y}_{n_s})^T$. We seek to construct a regression model using the Gaussian process that approximates an unknown function $f: \mathcal{M} \to \RE$ satisfying f($\bm{\mu}_i$) = ${y}_i$ for $i = 1,\dots,n_s$. We then use this model to estimate the data outputs $\bm{y}_{\star}$ corresponding to the parameter $\bm{\mu}_{\star}$ which is an element of a chosen test data set. 
% with the following data $(\bm{\mu}_i,{y}_i:=f(\bm{\mu}_i)), i=1,2,\dots,n_s.$
% To do so, we consider the following Gaussian process $\mathcal{GP}(m(\bm{\mu}),k(\bm{\mu},\bm{\mu}'))$. 
 Then,
\begin{equation}\label{gprmat}
\begin{bmatrix}
\bm{y}\\ {\bm{y}}_{\star}
\end{bmatrix}
\sim \mathcal{N}\Bigg(\begin{bmatrix}
m(\bm{M})\\ m(\bm{\mu}_{\star})
\end{bmatrix},\begin{bmatrix}
\Sigma & \Sigma_{\star}^T\\  \Sigma_{\star} &\Sigma_{\star\star}\\
\end{bmatrix}\Bigg),
\end{equation}
where $m(\bm{M})=(m(\bm{\mu}_1),\dots, m(\bm{\mu}_{n_s}))^T$, $\Sigma(i,j)=k(\bm{\mu}_i,\bm{\mu}_j),~\Sigma_{\star}=[k(\bm{\mu}_{\star},\bm{\mu}_1),\dots, k(\bm{\mu}_{\star},\bm{\mu}_{n_s})]$ and $\Sigma_{\star\star}=k(\bm{\mu}_{\star},\bm{\mu}_{\star}),$
 The conditional probability $p( {\bm{y}}_{\star}| \bm{y})$
also follows the Gaussian distribution~\cite{cbishop}:
$${\bm{y}}_{\star}|\bm{\bm{y}}\sim \mathcal{N}(m(\bm{\mu}_{\star})+\Sigma_{\star}\Sigma^{-1}(\bm{y}-m(X)),\Sigma_{\star\star}-\Sigma_{\star}\Sigma^{-1}\Sigma_{\star}^T).$$
The best estimate for $\bm{y}_{\star}$ is the mean of this conditional Gaussian distribution. In other words, it is $\bar{\bm{y}}_{\star}=m(\bm{\mu}_{\star})+\Sigma_{\star}\Sigma^{-1}(\bm{y}-m(\bm{M})).$ 
% Note that the mean is modified using the information of the given data. Also, the term $\Sigma_{\star}\Sigma^{-1}\Sigma_{\star}^T\geq 0$, so the posterior variance is less than the prior variance. 
{For noisy data, the matrix $\Sigma$ of the equation~\eqref{gprmat} is replaced by $\Sigma_n=\Sigma+\sigma_n^2 \mathbb{I}_{n_s}$.} In other words, we assume additive noise of Gaussian type with mean zero and variance $\sigma_n$. Here $\mathbb{I}_{n_s}$ denotes the identity matrix of size $n_s$. Finally, we use MATLAB command \texttt{fitrgp} to generate our GPR models.

An important feature of GPR is that it acts as a linear smoother \cite{HastieT}. In other words, the predicted value is a linear combination of the output or response variable $y$. Given a data set $\bm{y}$ parameterized by $\bm{\mu}$, the predicted value is then given by the general formula:
%{\color{red} change superscript}
\begin{equation} \label{eq:linearsmoother}
f(\bm{\mu_{\star}}) = \sum_{i=1}^{n} k(\bm{\mu_i},\bm{\mu_{\star}}) y_i.
\end{equation}

It is obvious that finding suitable kernels is the most important step in implementing this approach. In this paper, we will compare the following kernels. 
\begin{enumerate}
    \item Squared exponential (SE) kernel, defined as
    $$k(\bm{\mu},\bm{\mu}'|\sigma_f,\ell)=\sigma_f^2\exp(-\frac{|\bm{\mu}-\bm{\mu}'|^2}{2 \ell^2}).$$
     \item Absolute exponential (Exp) kernel, defined as
    $$k(\bm{\mu},\bm{\mu}'|\sigma_f,\ell)=\sigma_f^2\exp(-\frac{|\bm{\mu}-\bm{\mu}'|}{\ell}).$$
     \item Mat\'{e}rn kernel defined as
    $$k(\bm{\mu},\bm{\mu}'|\sigma_f,\ell)=\frac{\sigma_f^2}{2^{\nu-1}\Gamma(\nu)} \Big( \frac{\sqrt{2\nu}}{\ell} |\bm{\mu}-\bm{\mu}'|\Big)^{\nu}K_{\nu}\Big( \frac{\sqrt{2\nu}}{\ell}|\bm{\mu}-\bm{\mu}'|\Big),$$
    where $K_{\nu}$ is the modified Bessel function of the second kind.
    \begin{itemize}
    
        \item For $\nu=\frac{1}{2}$, the Mat\'{e}rn kernel coincides with the absolute exponential function;
         \item For $\nu=\frac{3}{2}$, $k(\bm{\mu},\bm{\mu}'|\sigma_f,\ell)=\sigma_f^2(1+\frac{\sqrt{3}|\bm{\mu}-\bm{\mu}'|}{\ell})\exp(-\frac{\sqrt{3}|\bm{\mu}-\bm{\mu}'|}{\ell})$;
         \item For $\nu=\frac{5}{2}$, $k(\bm{\mu},\bm{\mu}'|\sigma_f,\ell)=\sigma_f^2(1+\frac{\sqrt{5}|\bm{\mu}-\bm{\mu}'|}{\ell}+\frac{\sqrt{5}|\bm{\mu}-\bm{\mu}'|^2}{3\ell^2})\exp(-\frac{\sqrt{5}|\bm{\mu}-\bm{\mu}'|}{\ell})$,
    \end{itemize}

\end{enumerate}
where $\sigma_f$ is the standard deviation of the output data and $\ell$ is the length scale.  {If the mean function is constant that is $m(\bm{\mu})=\beta_0$ then
the hyperparameters of the GPR
% with this linear mean and any covariance is 
$\bm{\theta} = (\beta_0,\sigma_f,\ell,\sigma_n)$ are obtained by using the marginal likelihood function as an objective function. In particular, we find the optimal hyperparameters $\bm{\theta}$ of each kernel choice $k(\mu_i,\mu_j|\bm{\theta})$ by minimizing the following negative log marginal likelihood function:

$$L(\bm{\theta}) = -\frac{1}{2} \bm{y}^{T}K^{-1}\bm{y} - \frac{1}{2} \log(|K|) - \frac{n_s}{2} \log(2\pi) $$
where 
\begin{equation*}
K = 
\begin{bmatrix}
k(\mu_1,\mu_1) + \sigma_n^2 & \dots & k(\mu_1,\mu_{n_s})\\
\vdots & & \vdots\\
k(\mu_{n_s},\mu_1) & \dots & k(\mu_{n_s},\mu_{n_s}) + \sigma_n^2
\end{bmatrix}
\end{equation*}}

{In future works, we will explore the effects of generating new kernels from old kernels \cite{duvenaud_2014}}. For now, it is important to note that the regularity of a kernel is inherited by the sampled functions. If one wishes to make predictions using a data set that contains discontinuities, kernels with discontinuities will produce better predictions. %{\color{red} To this end, we introduce the following definition:
% \begin{defn}[Mean Square Continuity and Differentiability]~\cite{gp_ml}
% Let $\{\bm{\mu}_k\} _{k=1}^{\infty}$ be a sequence converging to $\bm{\mu}_{\star}$ in $\RE^d$ then a
% process $f(\bm{\mu})$ is continuous in mean square at $\bm{\mu}_{\star}$ if $E[|f(\bm{\mu}_k)-f(\bm{\mu}_{\star})|^2] \to 0$ as $k\to \infty$. A random field is continuous in
% mean square at $\bm{\mu}_{\star}$ iff its covariance function $k(\bm{\mu},\bm{\mu}')$ is continuous
% at the point $\bm{\mu}=\bm{\mu}'=\bm{\mu}_{\star}$. For stationary covariance functions (only depends on $\bm{\mu}-\bm{\mu}'$) this reduces to checking continuity at $k(0)$.
% \end{defn}
% }
The Gaussian process with squared exponential covariance function is infinitely differentiable~\cite{gp_ml}. That is, the GP has mean square derivatives of all orders and is thus very smooth. A Gaussian process with Matérn $\nu$ covariance function is $ \ceil{\nu}-1 $ times differentiable in the mean-square sense. Thus the Gaussian process with the absolute exponential kernel is mean-square continuous, with Mat\'{e}rn 3/2 kernel is one time differentiable, and with Mat\'{e}rn 5/2 kernel is twice differentiable in the mean-square sense. Such properties will prove useful when we explore using these kernels to generate meaningful predictions in Section~\ref{se:Num}.

% Thus, we can hypothesize that using if we want to approximate a non-smooth function using GPR, then the GPR with the absolute exponential kernel or Mat\'{e}rn kernel will be better than the squared exponential kernel which we observed in our numerical experiments.

\section{A comparison between splines and GPR} \label{se:GPRSplines}

An important goal of this work is to show that GPR should be considered a serious alternative to the more popular spline methods in the field of PDEs. Indeed, the existence of correspondence between GPR and spline methods has been known since the 1970s \cite{Kimeldorf1970}. Furthermore, GPR is a supervised learning algorithm that is becoming increasingly important in fields such as machine learning. Given the increasing interest, it is expected that advances in understanding such a connection will have a large impact on both fields. 

Here, we present a comparison between the performance of GPR and spline methods. In particular, we discuss an example where our data is sampled from a familiar setting: a smooth function with a large curvature, namely:

$$ f(\mu) = 1 + \frac{\sin(\mu)}{\mu}.$$
% We compare the performance of GPR, cubic spline interpolation, and linear interpolation. 
The kernel used for GPR in this example is the squared exponential kernel. %We consider two cases: a uniform and a non-uniform grid.
The training set is chosen in two different ways. The first is the uniform grid case where $\mu$ belongs to a uniform discretization of the parameter space $\mathcal{M} = [-\pi,3 \pi]$ with step size $h$. The second is a non-uniform case. In both cases, the test set $\{f(\mu_{\star}),\mu_{\star} \}$ is generated on a uniform grid with a step size of $0.01$.

In both cases, GPR proves superior to splines. This is, however, not always the case for any data set generated from a parametric eigenvalue problem, and there are situations where splines perform better. A throughout comparison of the two approaches as well as theoretical discussions, will be the object of future investigations.
\subsection{Case I}
\begin{figure}
\begin{subfigure}{0.38\textwidth}
\centering
\includegraphics[height=5cm,width=6.5cm]{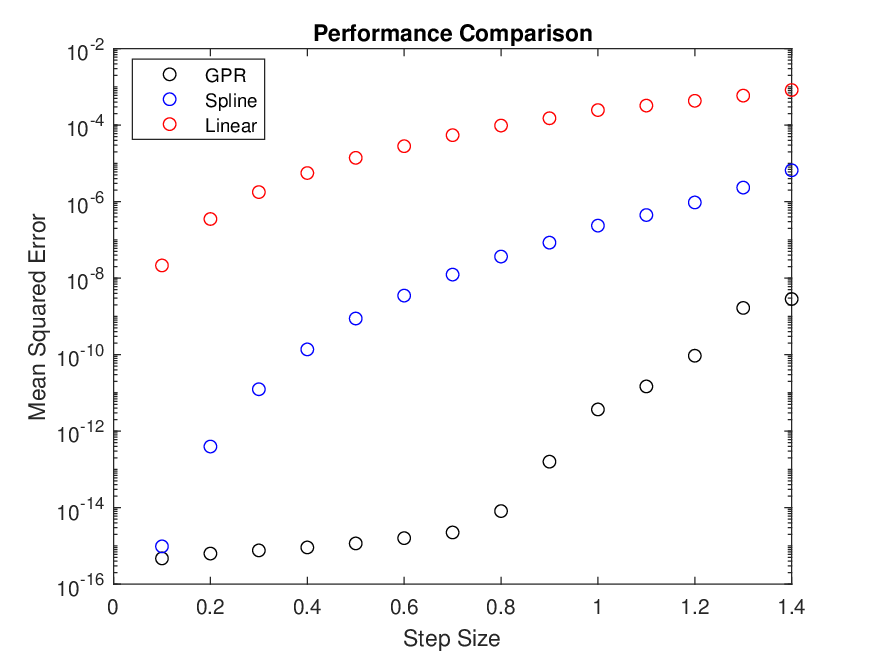}
\caption{MSE Error}
\label{ev:PerformanceComparison1}
\end{subfigure}
\begin{subfigure}{0.38\textwidth}
\centering
\includegraphics[height=5cm,width=6.5cm]{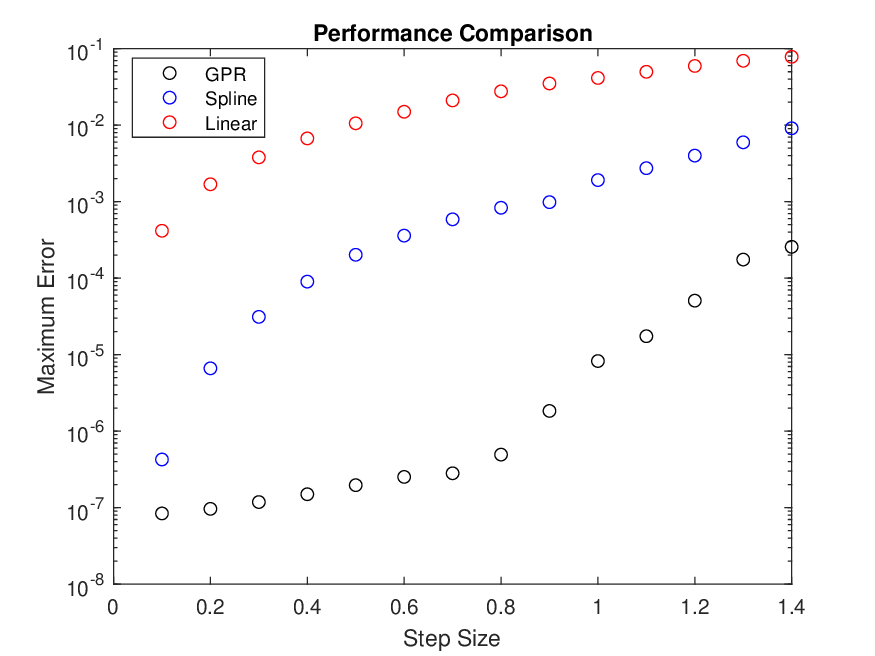}
\caption{MAX Error}
\label{ev:PerformanceComparison2}
\end{subfigure}
\caption{{Performance comparison of spline and GPR in the uniform case}}
\label{ev:PerformanceComparison}
\end{figure}

We analyze the error as a function of the step size. The error is calculated in two different ways: 1) the mean squared error, and 2) the maximum error. A plot of the error as a function of the step size is given in Fig.~\ref{ev:PerformanceComparison}.

In this case, GPR is shown to outperform both of the spline and linear methods. %Must have big enough data set?
\subsection{Case II}
Next, we consider the following non-uniform grid. The grid $G \subset \mathcal{M}$ is: 

$G = [-4.36, -4.04, -2.44, 1.73, 2.05, 4.94, 5.26, 8.14, 9.10]$

\begin{figure}
\centering
\begin{subfigure}{0.3\textwidth}
\includegraphics[height=5cm,width=5.3cm]{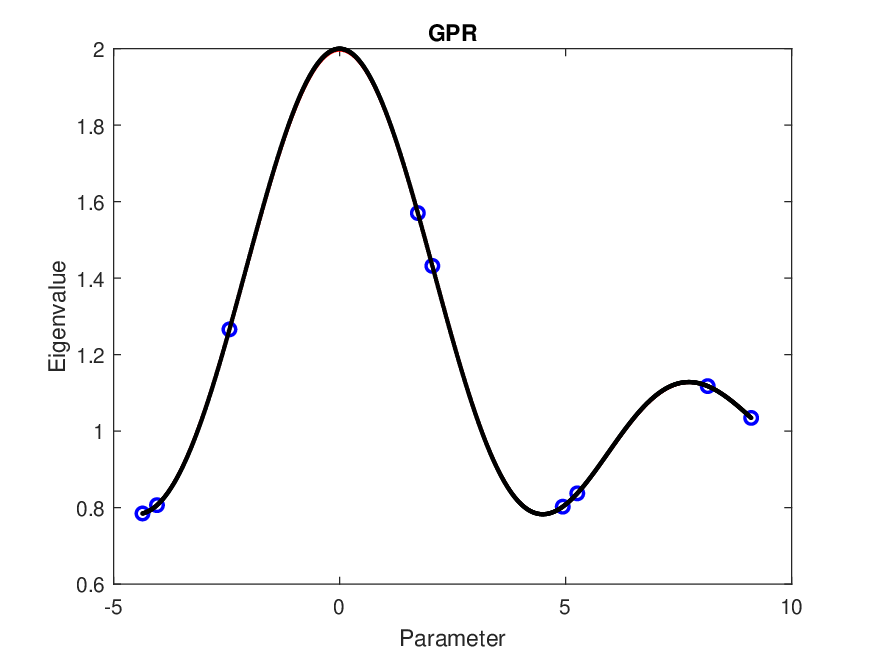}
\caption{GPR}
\label{ev:GPRvsSplines1}
\end{subfigure}
\begin{subfigure}{0.3\textwidth}
\includegraphics[height=5cm,width=5.3cm]{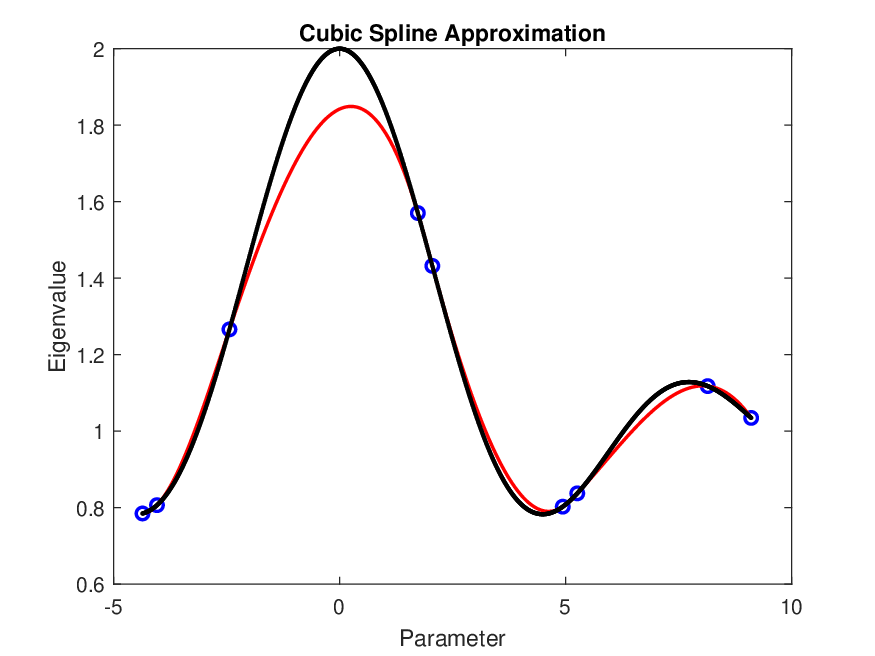}
\caption{Cubic spline interp.}
\label{ev:GPRvsSplines2}
\end{subfigure}
\begin{subfigure}{0.3\textwidth}
\includegraphics[height=5cm,width=5.3cm]{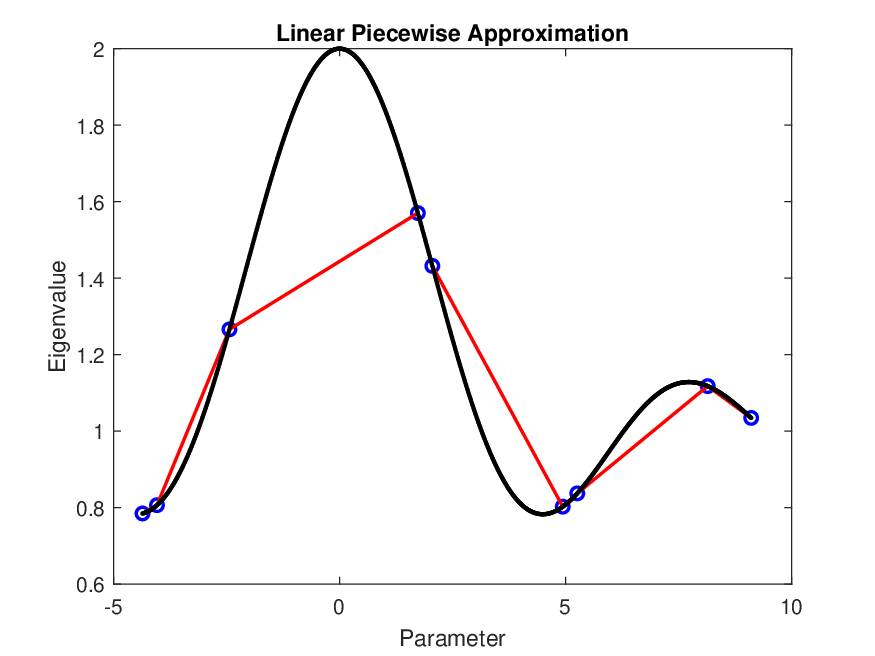}
\caption{Linear piecewise interp.}
\label{ev:GPRvsSplines3}
\end{subfigure}
\caption{Comparison between the function $f(\mu)$ (black) and estimates using three different methods. Blue circles are the training data points.}
\label{ev:GPRvsSplines}
\end{figure}

Due to non-uniformity, some interpolation methods do not work so well. Furthermore, one does not expect that the error using spline methods will approach zero when the size of the larger interval in the grid is not going to zero. GPR, on the other hand, performs fairly well, as can be seen in Fig. \ref{ev:GPRvsSplines}. In particular, it is able to estimate better regions where the second derivative is largest.

\section{Numerical results} \label{se:Num}
In this section, we apply GPR to three different parametric eigenvalue problems. We focus on comparing the numerical results generated using four different covariance functions. The eigenvalue problems considered here have been selected for the following reasons. In the first eigenvalue problem, the eigenvalues cross each other. Such crossings lead to eigenvalue curves with singularities which are common in parametric eigenvalue problems. It is thus instructive to discuss in detail how GPR can be used in such a context. While the curves corresponding to the eigenvalues are continuous but maybe not differentiable, the corresponding eigenvectors present discontinuities when passing across the intersections. The curves considered in the second example, on the other hand, contain no corners/discontinuities. Here, we focus on exploring the performance of our GPR models for a collection of smooth curves with different frequency contents. 
{Finally, the last problem is a two-dimensional eigenvalue problem, an example where we show that our technique is not limited to one-dimensional problems.}

Our analysis is structured as follows. First, let $\bm{\mu_{tr}} = [\mu_1,\dots,\mu_N]$ and $\bm{\lambda_{tr}} = [\lambda_1,\dots,\lambda_N]$ be the training input and output variables, respectively. {In addition to the eigenvalues, also the eigenfunctions $\bm{u_{tr}}$ are computed as output variables.} We generate our training output variables $\bm{\lambda_{tr}}$ {and $\bm{u_{tr}}$} by solving for the eigenvalues {and eigenfunctions} at the training inputs $\bm{\mu_{tr}}$ using the finite element reduced order modeling technique described in Section~\ref{se:ps}. Second, we apply the GPR technique described in Section~\ref{se:GPR} to make predictions at a different collection of test points, $\bm{\mu_{test}}$. We describe in each example how $\bm{\mu_{test}}$ is chosen. Third, we compute the true eigenvalues and eigenvectors at $\bm{\mu_{test}}$. This allows us to compute the error between the predicted and true results at these test points.
{To quantify errors for the eigenvalues, we have used the following error measures: the relative root mean squared error (RRMSE), root mean squared relative error (RMSRE), root mean squared error (RMSE), and relative squared error (RSE). These are defined as follows:
\begin{align*}
  & RRMSE = \sqrt{\frac{\frac{1}{n_t}\sum\limits_{i=1}^{n_t} (\lambda_i-\hat{\lambda}_i)^2}{\sum\limits_{i=1}^{n_t} \lambda_i^2}},\quad {RMSRE = \sqrt{\frac{1}{n_t}\sum\limits_{i=1}^{n_t} \frac{(\lambda_i-\hat{\lambda}_i)^2}{\lambda_i^2} }},\\ 
  &RMSE = \sqrt{\frac{1}{n_t}\sum\limits_{i=1}^{n_t} (\lambda_i-\hat{\lambda}_i)^2}, \quad RSE=\frac{\sum\limits_{i=1}^{n_t} (\lambda_i-\hat{\lambda}_i)^2}{\sum\limits_{i=1}^{n_t} (\lambda_i-\bar{\lambda})^2} \quad \text{with} \quad \bar{\lambda}=\frac{1}{n_t}\sum\limits_{i=1}^{n_t}\lambda_i
\end{align*} where $\lambda_i$ and $\hat{\lambda}_i$ are respectively the FEM and GPR based eigenvalues,  and $n_t$ denotes the number of test points used.
We observe that RRMSE has a different scaling from the other error quantities, since the denominator is not divided by the number of test points $n_t$. The other three errors follow the same pattern so we only report the RMSRE for comparison purposes.}
To compare the eigenvectors, we plot the error between the FEM eigenvectors and GPR-based eigenvectors after normalizing using the $L^2$-norm and choosing the proper sign.
{In analogy of what is done for the eigenvalues, the RMSRE could be computed also for the coefficients of the eigenvectors.}

{Before we analyze the results, we first discuss some technical details regarding our use of the \texttt{fitrgp} library. In all experiments, we set the explicit basis function to be the constant function. Thus, one needs to estimate only one coefficient, $\beta_0$. Furthermore, each kernel contains two hyper-parameters, the length scale $\ell$ and the signal standard deviation $\sigma_f$. Furthermore, we assume our data to be noise-free. In practice, however, we set the noise to be very small, $\sigma_n = 10^{-4}$. In particular, one needs to choose the following setting in \texttt{fitrgp} to achieve this: \texttt{'ConstantSigma', true}. Furthermore, a lower bound has to also be set. In our case, it is \texttt{'SigmaLowerBound',1e-6}. This needs to be always smaller than the chosen $\sigma_n$.}  

{In order to obtain the best fit, we minimize the negative marginal log-likelihood function. To do this, we use the \texttt{'quasinewton'} optimizer. The other available optimizers either a) do not produce better results, and/or b) require additional toolboxes and are more computationally expensive.}

{We initialize our hyperparameters as follows. We set the initial values for $\ell$ and $\sigma_f$  to be $\texttt{mean(std}(\mathcal{M}_{tr}))$ and \texttt{std(y)/$\sqrt(2)$}, respectively. Here, $\mathcal{M}_{tr}$ is the set of training parameters and $y$ is the corresponding output vector.}

{Furthermore, we use for the \texttt{fitMethod} the option: \texttt{'exact'}. All other options (namely, \texttt{'bcd'},\texttt{'sd'},\texttt{'sr'} and \texttt{'fic'}) become noticeably more effective, by design, for large datasets only ($>$ 10000 points). This also applies to \texttt{PredictMethod}, which we choose to also be \texttt{'exact'} instead of \texttt{'sd'},\texttt{'sr'} and \texttt{'fic'}.}

\begin{figure}
\centering
\includegraphics[width=5cm]{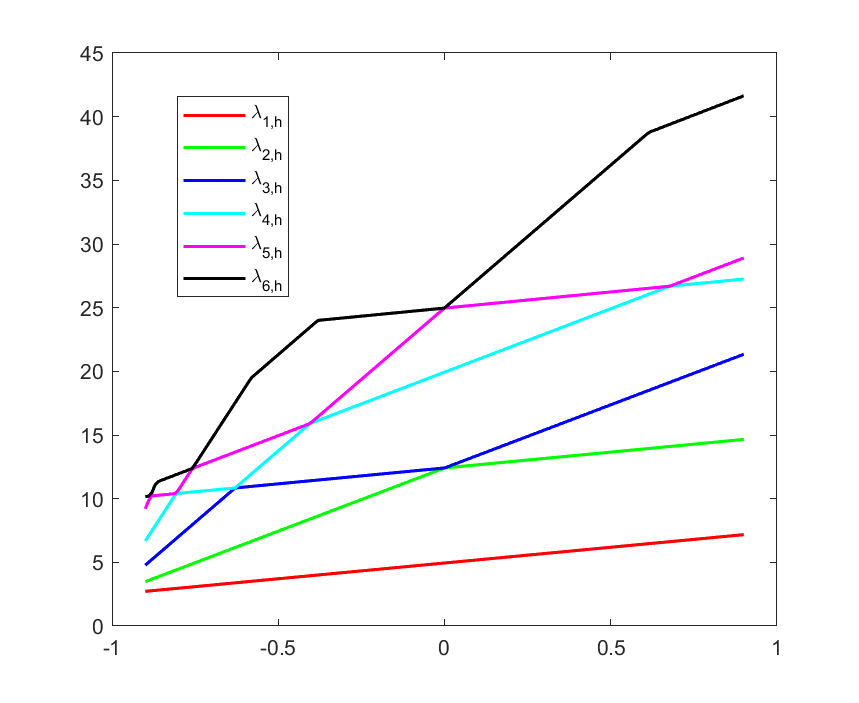}
\caption{First six eigenvalues of the eigenvalue problem \eqref{eq:crs}.}
\label{ev:exp1}
\end{figure}

\begin{example}
Let us consider the following eigenvalue problem 
\begin{equation}
\label{eq:crs}
    \left\{
\begin{array}{ll}
   -\nabla \cdot(A(\mu)\nabla u(\mu))=\lambda(\mu) u & \text{in }\Omega=(-1,1)^2 \\
     u(\mu)=0& \textrm{on }\partial\Omega 
\end{array}
\right.
\end{equation}
where $A(\mu)=\begin{bmatrix}
1 &0 \\
0 & 1+\mu
\end{bmatrix}$. We set our parameter space to be $\mathcal{M}=[-0.9,0.9]$. Furthermore, we use the following notation to represent our data sets. We denote by $a:c:b$ a uniformly sampled grid from the 1-dimensional interval $[a,b]$ with step size $c$. 

The exact solution is known for this problem. The eigenvalues and eigenvectors are given by
\begin{align*}
    &\lambda_{m,n}(\mu)=\frac{\pi^2}{4}(m^2+(1+\mu)n^2) \\
    &u_{m,n}= \cos\left(\frac{m\pi}{2}x\right)\cos\left(\frac{n\pi}{2}y\right),\quad \forall m,n\in \mathbb{N}.
\end{align*}

The first six eigenvalues of the eigenvalue problem~\eqref{eq:crs} are plotted in Fig.~\ref{ev:exp1}. We do not consider the first eigenvalue curve since it is completely separated from the others and thus is trivial. Our training points are $\bm{\mu_{tr}}=-0.9:0.1:0.9$ and the test points are $\bm{\mu_{test}}=-0.9:0.05:0.9$. 
\end{example}

% In Fig.~\ref{ev1:crs}, we consider the first eigenvalue curve which is separated from all of the other eigenvalues. The GPR predictions for this curve first eigenvalue of the eigenvalue problem~ \eqref{eq:crs} using four different covariance functions/kernels. The first eigenvalues at four test points $\mu=-0.75,-0.25,0.25, 0.75$ are reported in Tab.~\ref{crs:ev1}. For all four cases, the eigenvalues \re{are the same up to four decimal places and the RRMSE corresponding to the Mat\'{e}rn kernels are comparable and little better than the squared exponential kernel. The RRMSE is calculated using $40$ random points selected from the parameter space $\mathcal{M}$.}

% Since the eigenvectors are parameter independent and the first eigenvalue is separated from the other so there is only one reduced co-efficient for the eigenvectors. As we are getting the same values for all cases, the first eigenvalues it is not interesting to us so we did not include the GPR for the eigenvectors. 

\begin{figure}
\centering
 \begin{subfigure}{0.32\textwidth}
   \centering
   \includegraphics[height=4cm,width=5cm]{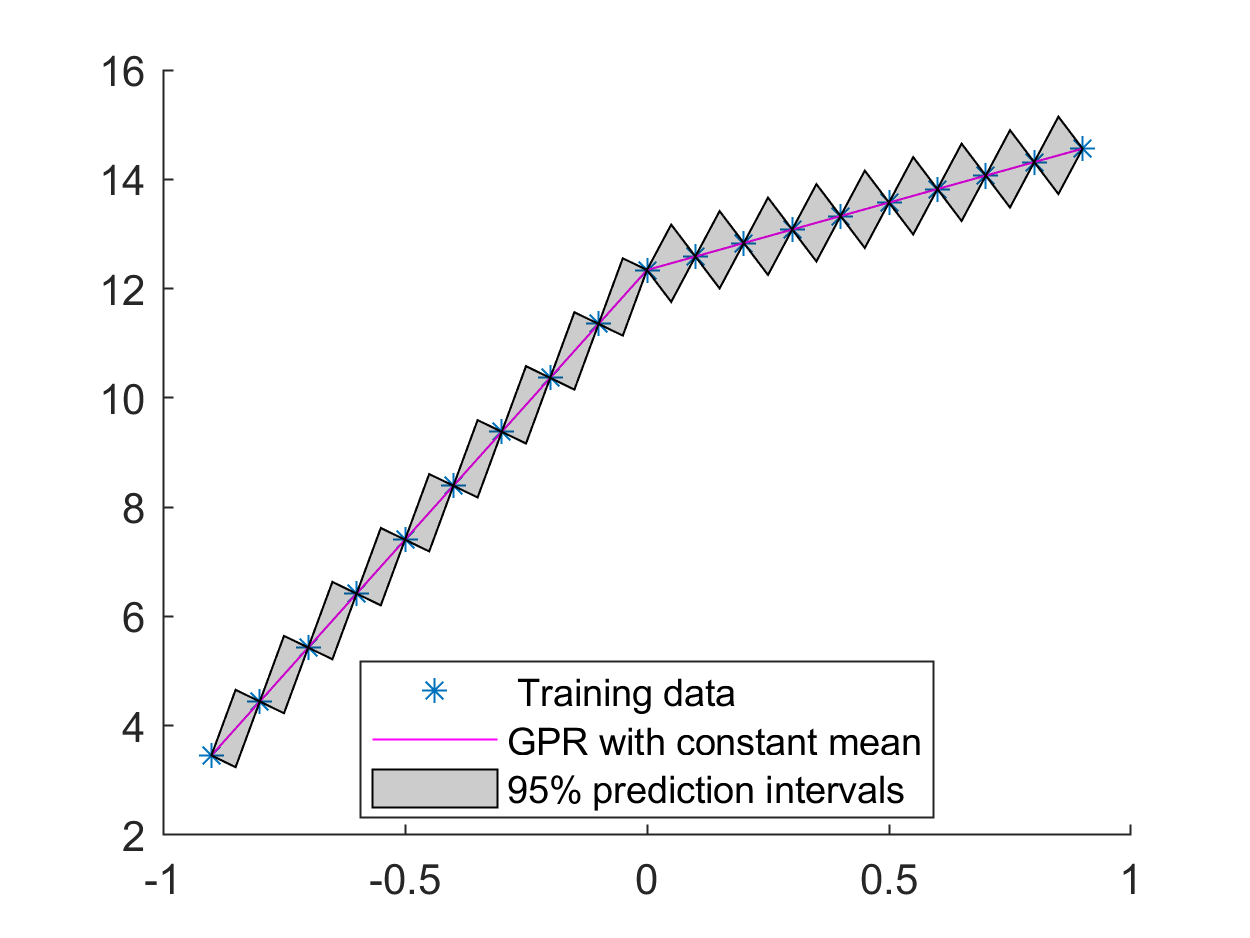}
    \caption{Exp Kernel}
  \end{subfigure}
 \begin{subfigure}{0.32\textwidth}
   \centering
   \includegraphics[height=4cm,width=5cm]{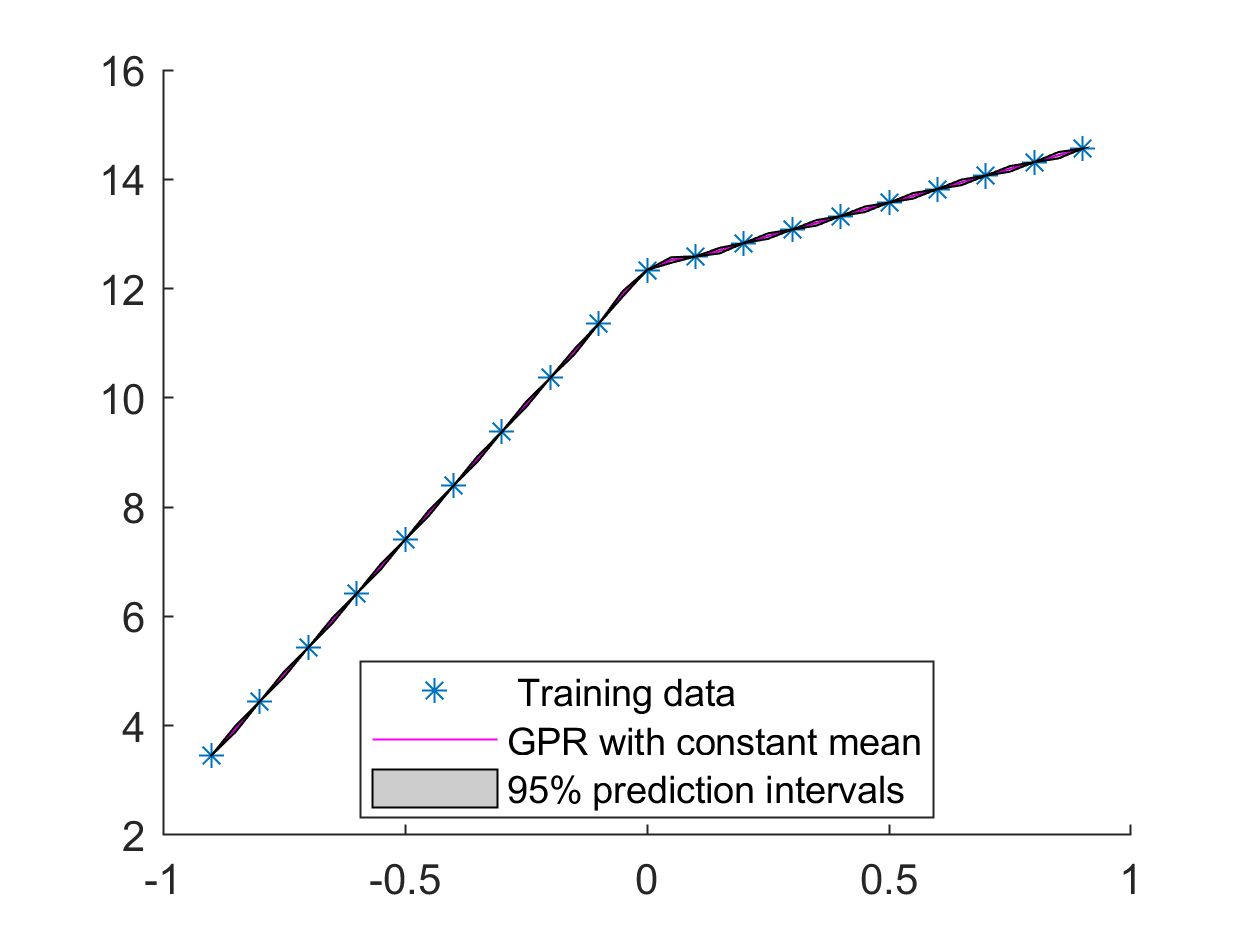}
    \caption{Mat\'{e}rn 3/2 Kernel}
  \end{subfigure}\\
  \begin{subfigure}{0.32\textwidth}
   \centering
   \includegraphics[height=4cm,width=5cm]{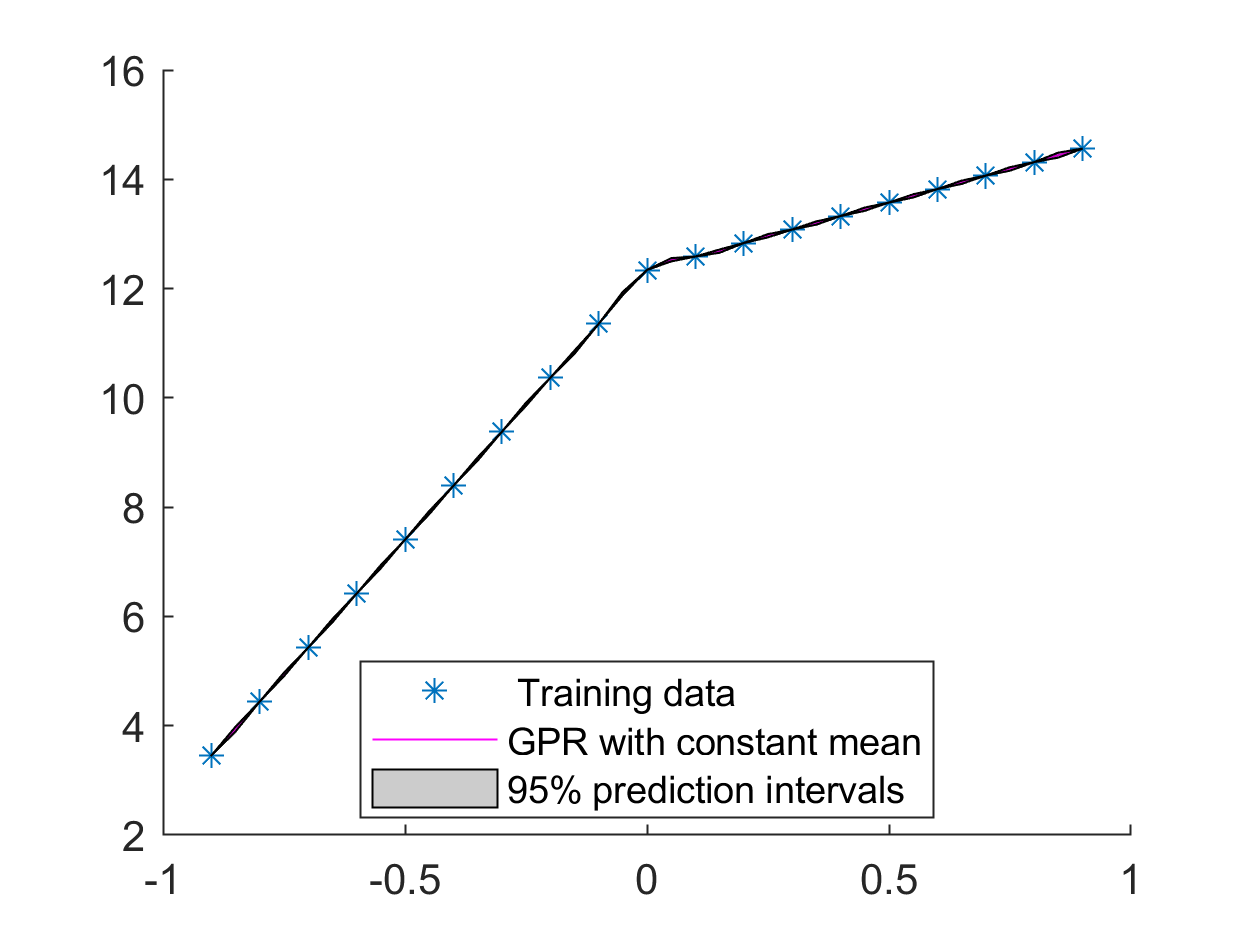}
    \caption{Mat\'{e}rn 5/2 Kernel}
  \end{subfigure}
  \begin{subfigure}{0.32\textwidth}
   \centering
   \includegraphics[height=4cm,width=5cm]{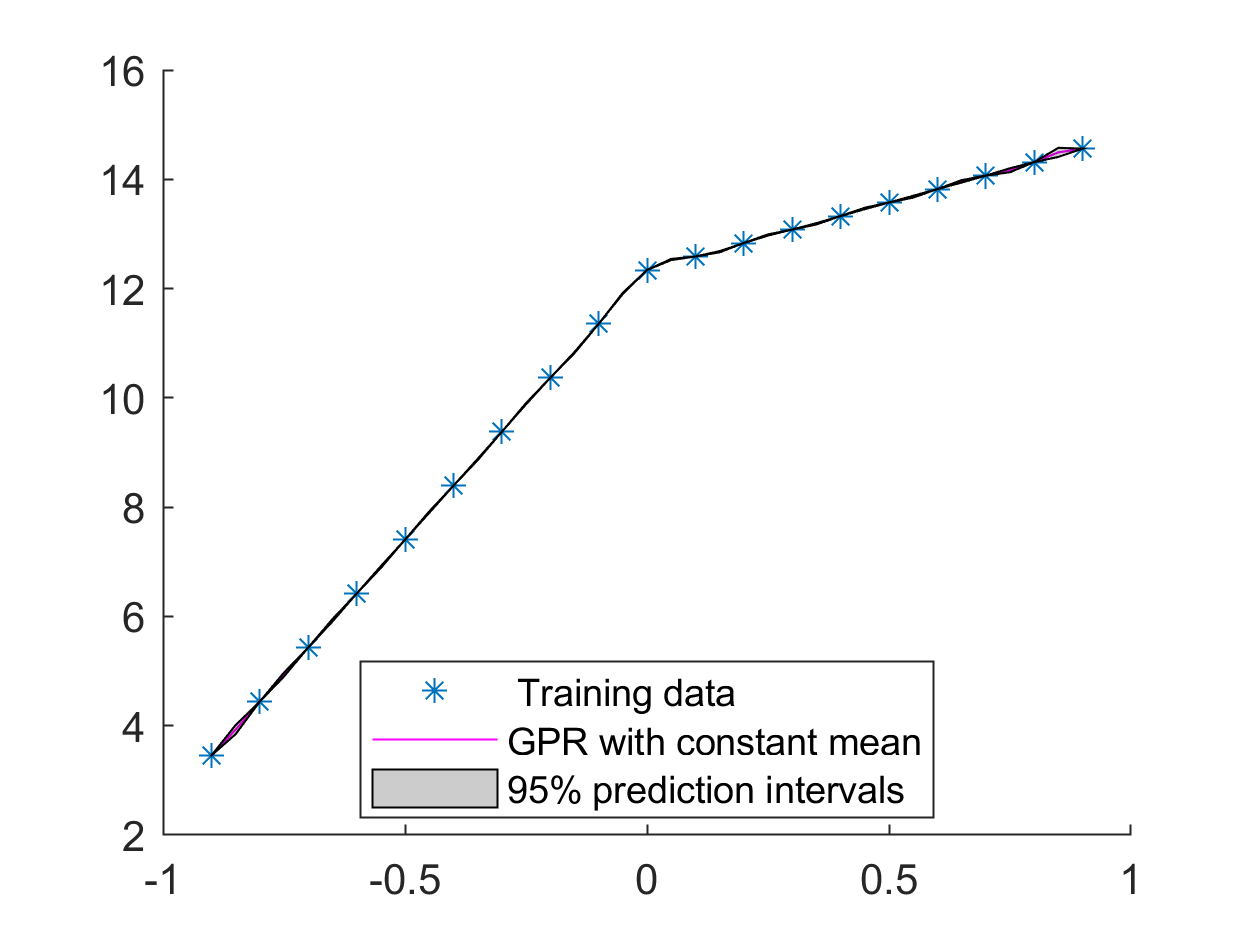}
    \caption{SE Kernel}
  \end{subfigure}
  \caption{GPR corresponding to 2nd eigenvalues of EVP~\eqref{eq:crs} using different kernels.}
  \label{ev2:crs}
\end{figure}

In Fig.~\ref{ev2:crs}, we show the GPR predictions for the second eigenvalue of the EVP~\eqref{eq:crs} using four different covariance functions. In addition, we also show the associated 95 percent confidence intervals. In Tab.~\ref{crs:ev2}, we show a sample of the predictions at the following four test points $\{-0.75,-0.25,0.25, 0.75\}$. To compute the RMSRE listed in Tab.~\ref{crs:ev2}, we use a randomly generated (but fixed) sample set containing 40 points. This avoids selecting test points that might favor our GPR model. Based on the \re{RMSRE} values reported, it can be concluded that the absolute exponential kernel performs best. This kernel's enhanced performance is due to the corner at $\mu = 0$ in the true function. {To check the performance of the GPR with different kernels, we increased the number of test points and calculated the corresponding RMSRE. The RMSRE values reported in Tab.~\ref{crs:ev2_rrse} indicate that the absolute kernel performs better than the other ones also when we increase the number of test points.}

{To analyze the associated eigenvectors, we plot the point-wise error between the eigenvector obtained by FEM and the eigenvector obtained by the GPR at some test points. We need to estimate two reduced coefficients for the second eigenvector of the EVP~\eqref{eq:crs}}. This is because the second eigenvalue curve comes from two different analytical solutions. This indicates that there are only two independent columns of the snapshot matrix since the eigenvectors are parameter-independent. 

\begin{figure}
\centering
 \begin{subfigure}{0.35\textwidth}
   \centering
   \includegraphics[height=4cm,width=5cm]{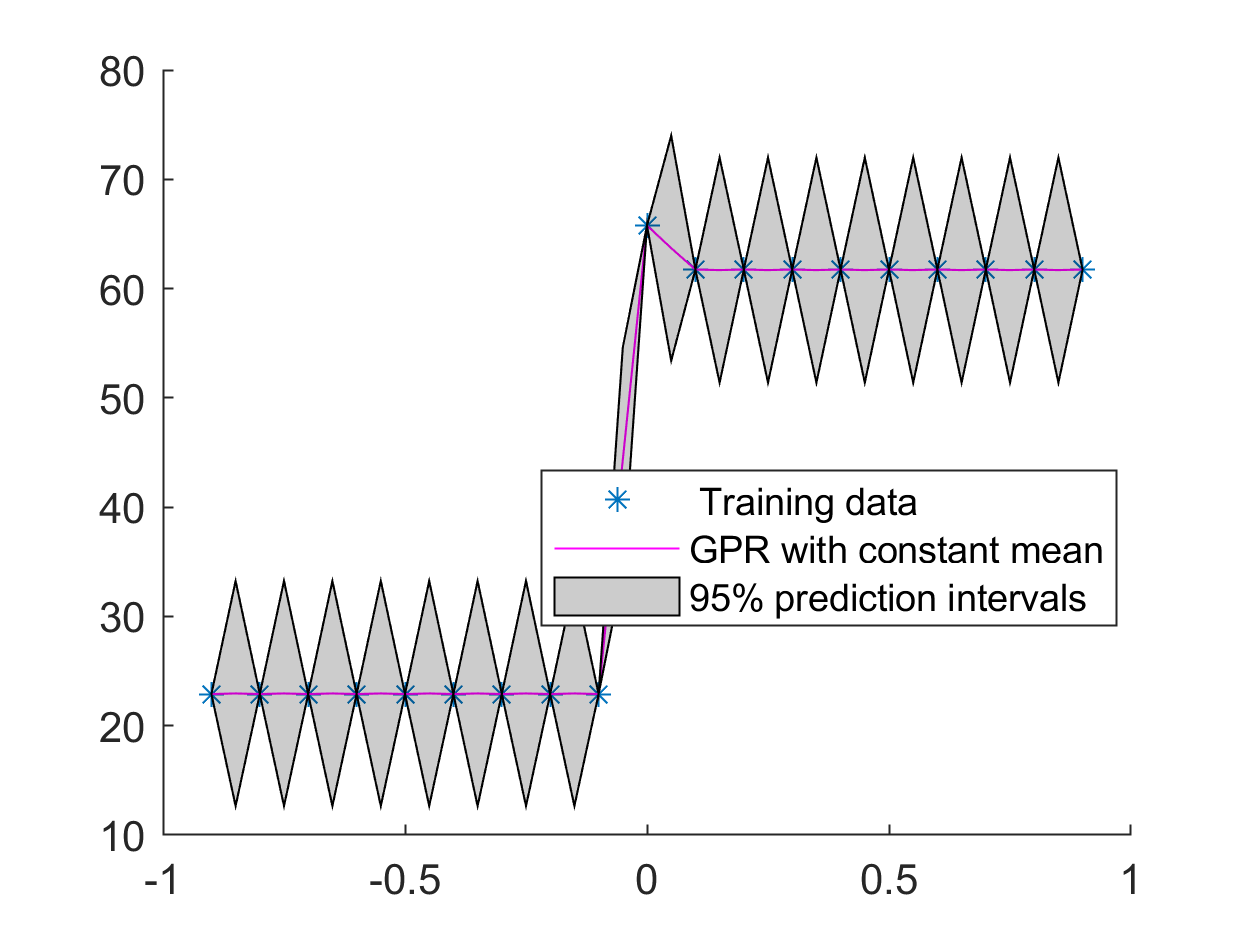}
    \caption{Exp Kernel}
  \end{subfigure}
 \begin{subfigure}{0.35\textwidth}
   \centering
   \includegraphics[height=4cm,width=5cm]{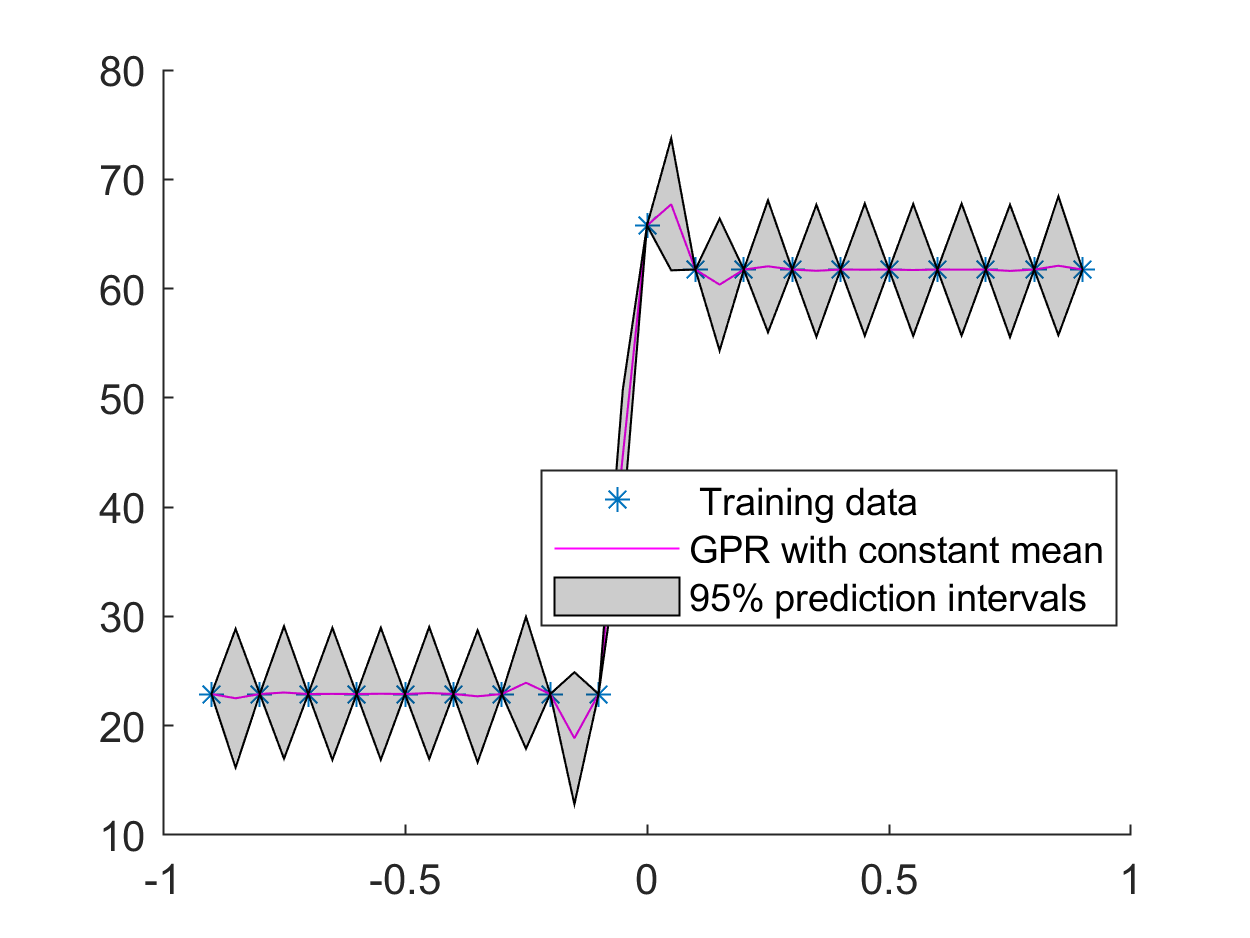}
    \caption{Mat\'{e}rn 3/2 Kernel}
  \end{subfigure}
  \begin{subfigure}{0.35\textwidth}
   \centering
   \includegraphics[height=4cm,width=5cm]{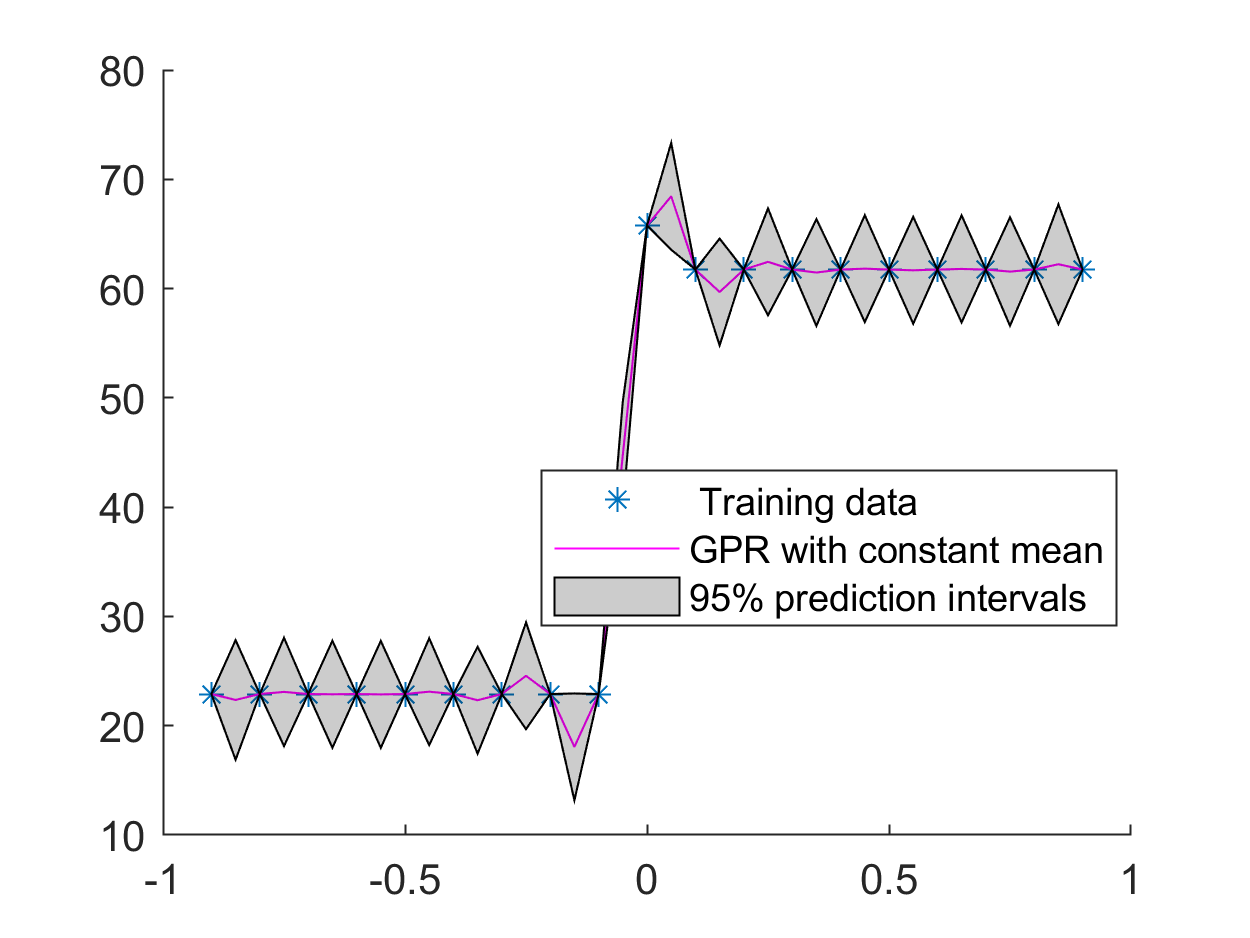}
    \caption{Mat\'{e}rn 5/2 Kernel}
  \end{subfigure}
  \begin{subfigure}{0.35\textwidth}
   \centering
   \includegraphics[height=4cm,width=5cm]{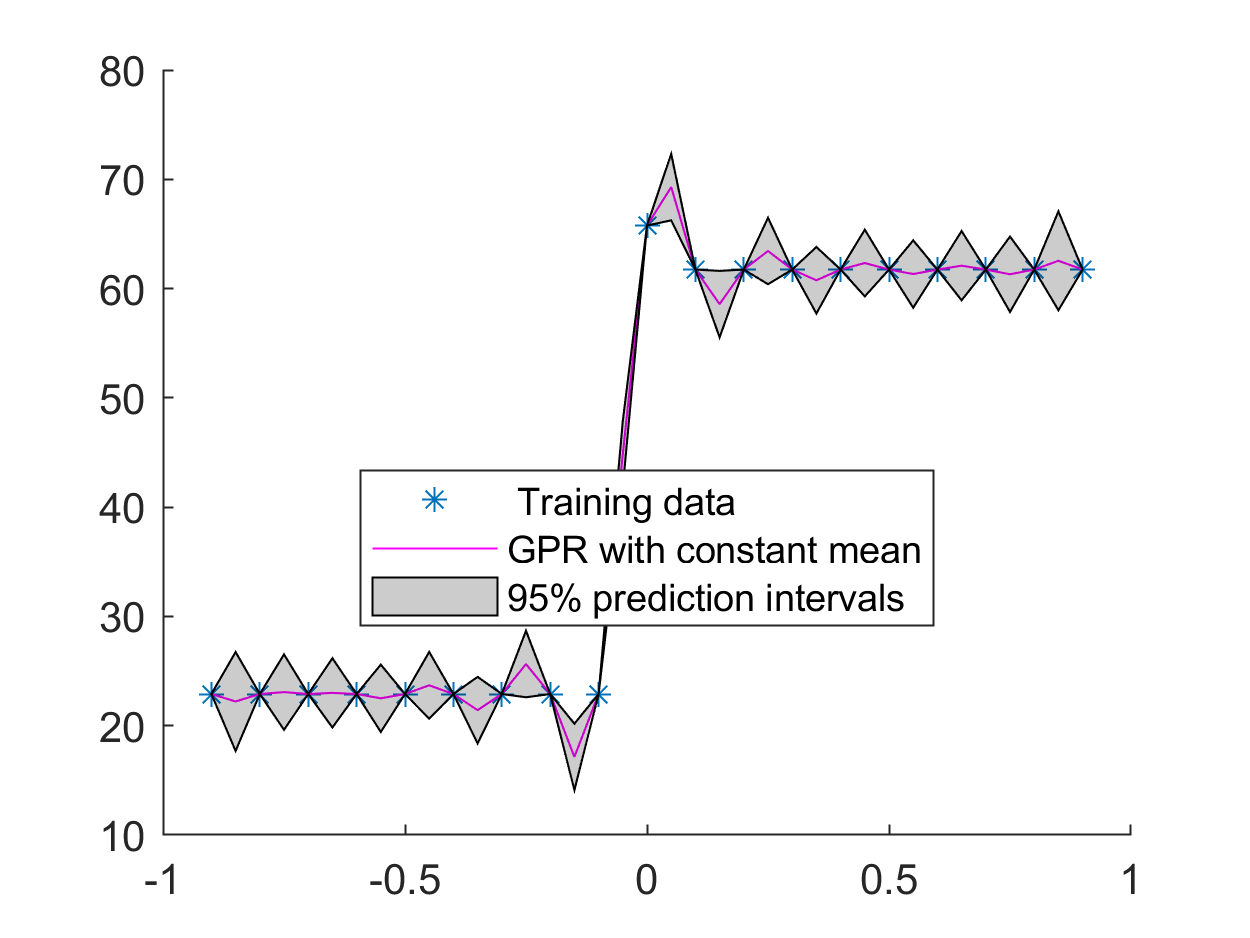}
    \caption{Squared Exponential Kernel}
  \end{subfigure}
  \caption{GPR corresponding to 1st coefficient of the reduced 2nd eigenvector of EVP~\eqref{eq:crs} using different kernels.}
  \label{crs:ev2_coeff1}
\end{figure}

\begin{figure}
\centering
 \begin{subfigure}{0.35\textwidth}
   \centering
   \includegraphics[height=4cm,width=5cm]{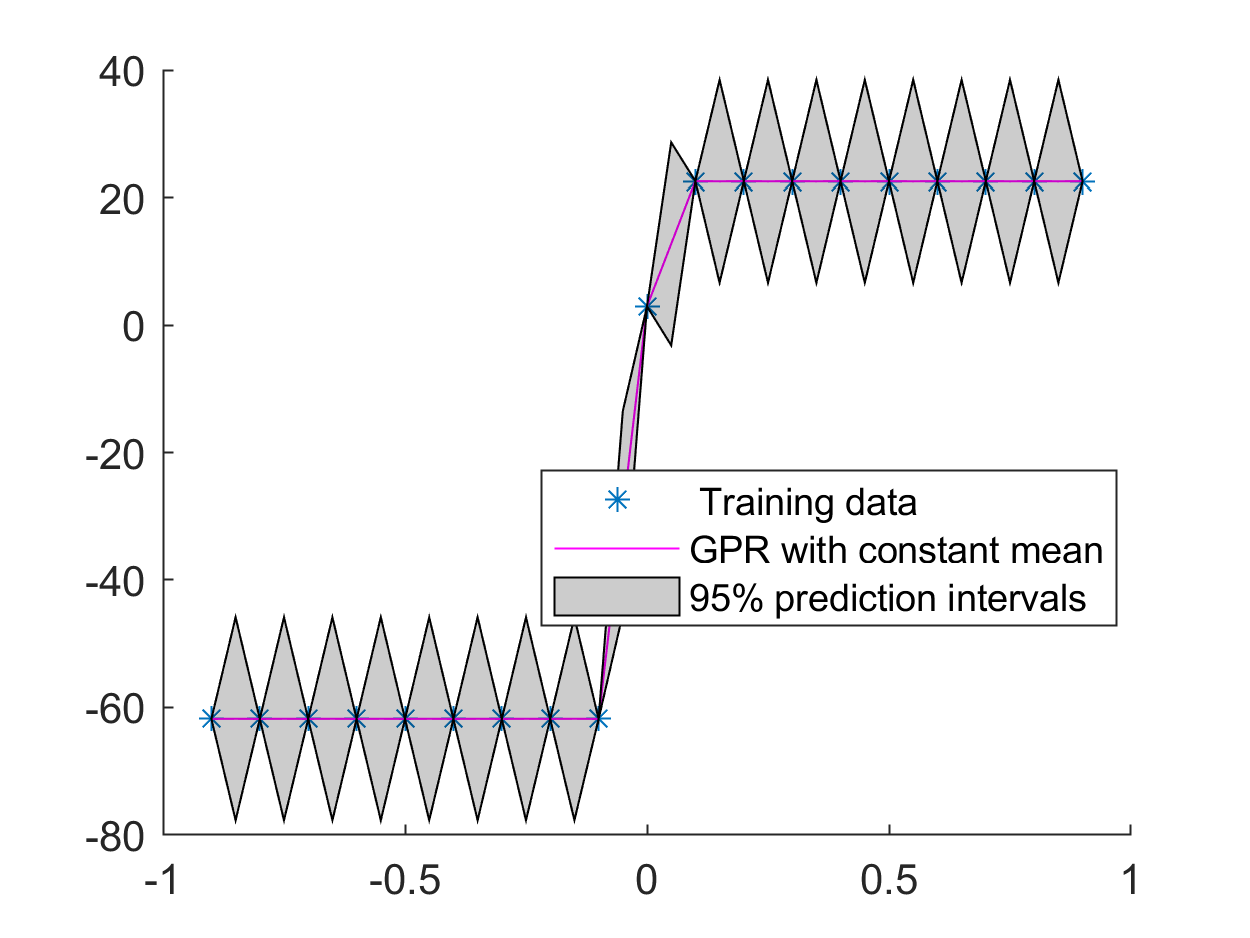}
    \caption{Exp Kernel}
  \end{subfigure}
 \begin{subfigure}{0.35\textwidth}
   \centering
   \includegraphics[height=4cm,width=5cm]{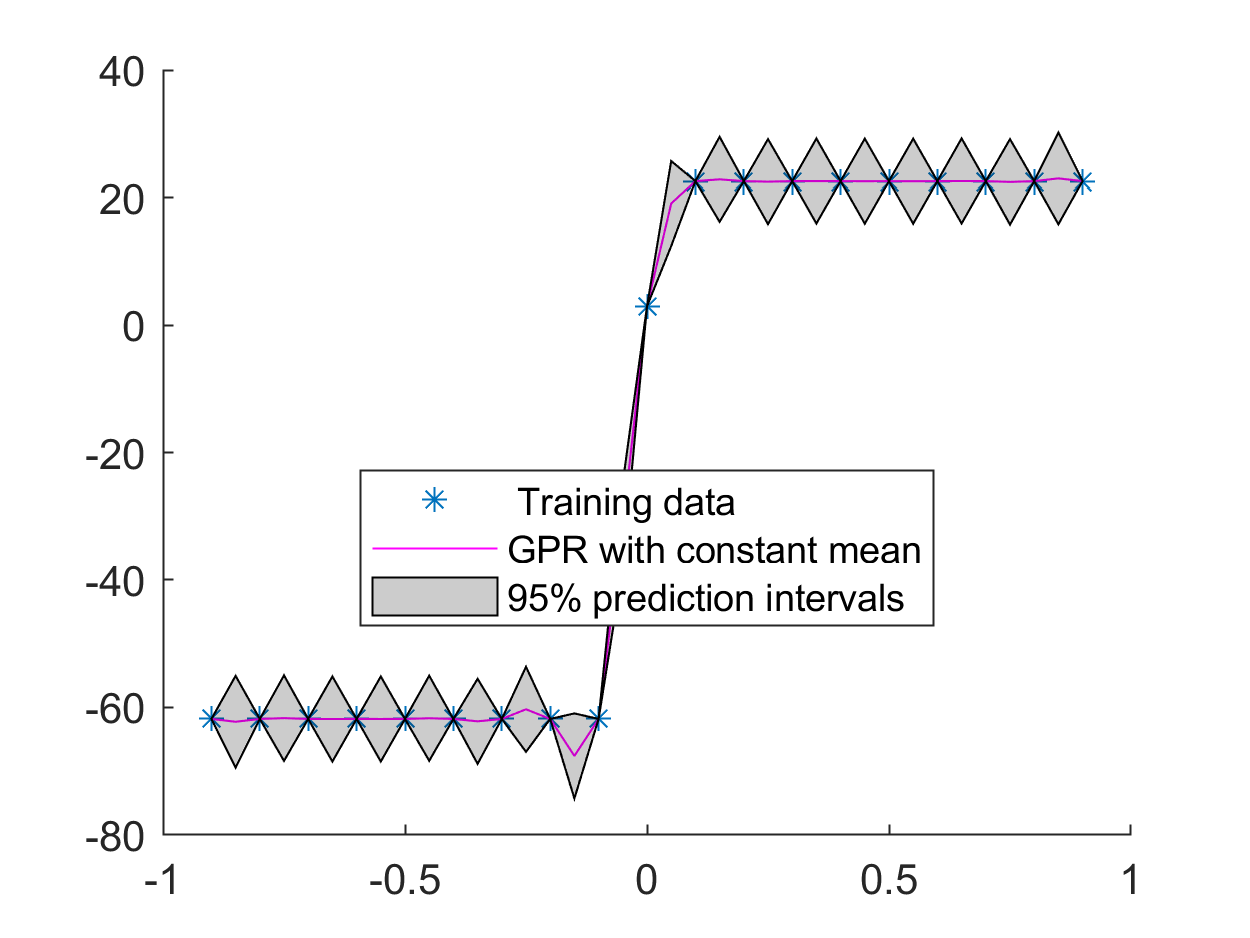}
    \caption{Mat\'{e}rn 3/2 Kernel}
  \end{subfigure}
  \begin{subfigure}{0.35\textwidth}
   \centering
   \includegraphics[height=4cm,width=5cm]{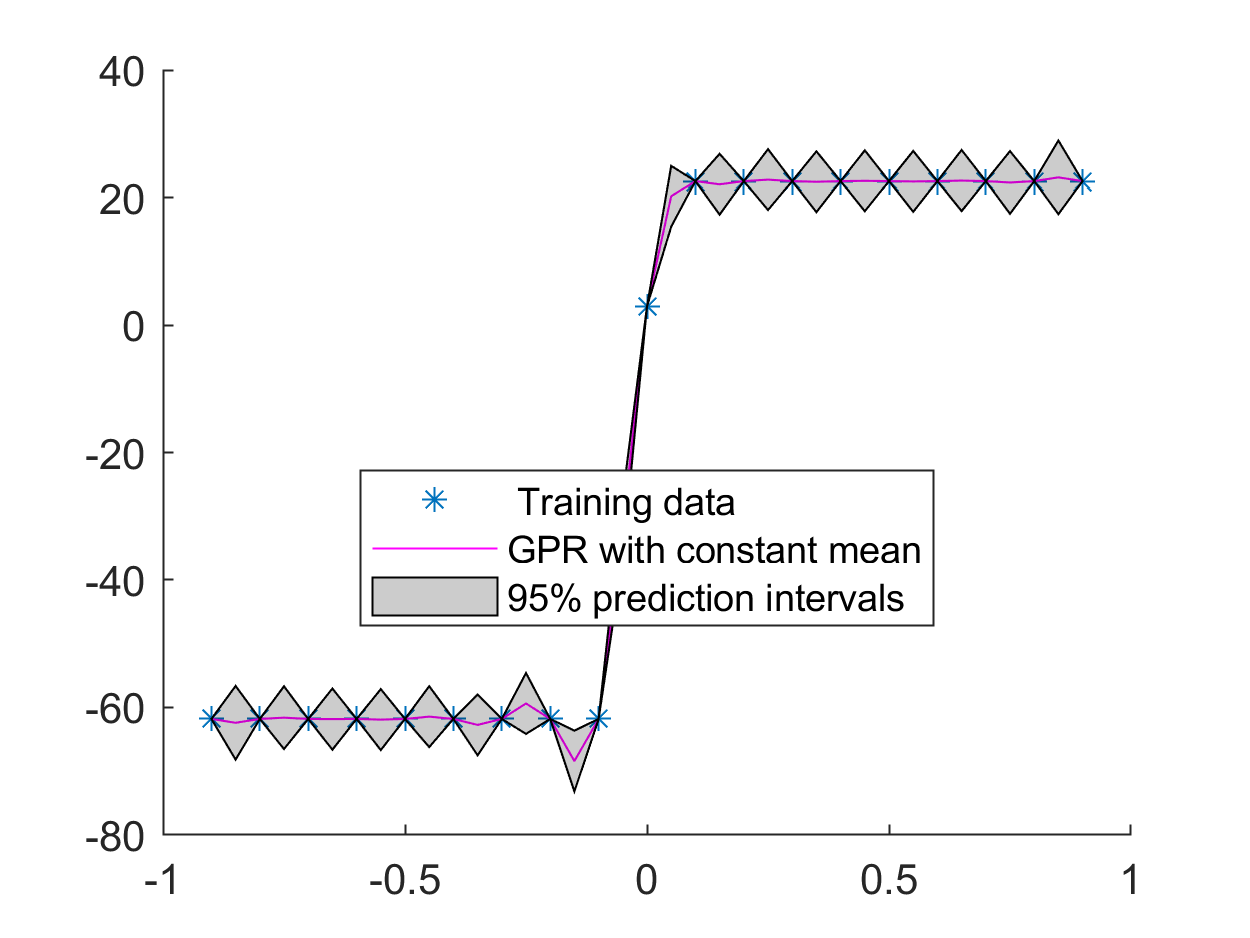}
    \caption{Mat\'{e}rn 5/2 Kernel}
  \end{subfigure}
  \begin{subfigure}{0.35\textwidth}
   \centering
   \includegraphics[height=4cm,width=5cm]{Result2/evct2_dm1_ma5_coef2.png}
    \caption{Squared Exponential Kernel}
  \end{subfigure}
  \caption{GPR corresponding to second coefficient of the reduced 2nd eigenvector of EVP~\eqref{eq:crs} using different kernels.}
  \label{crs:ev2_coeff2}
\end{figure}

In Fig.~\ref{crs:ev2_coeff1}-Fig.~\ref{crs:ev2_coeff2}, we {plot} the GPR for the two coefficients corresponding to the second eigenvector. In both figures, we can see that using the absolute exponential kernel leads to better prediction power. {In the case of the absolute exponential kernel, there is no oscillation near the jump point $\mu=0$. However, for the other kernels, oscillation is observed near the point $\mu=0$.}

\begin{figure}
\centering
  \begin{subfigure}{0.35\textwidth}
   \includegraphics[height=4cm,width=5cm]{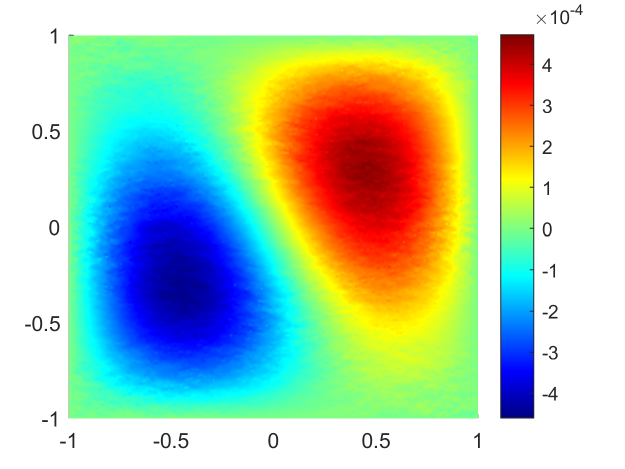}
    \caption{ Exp Kernel}
  \end{subfigure}
  \begin{subfigure}{0.35\textwidth}
   \includegraphics[height=4cm,width=5cm]{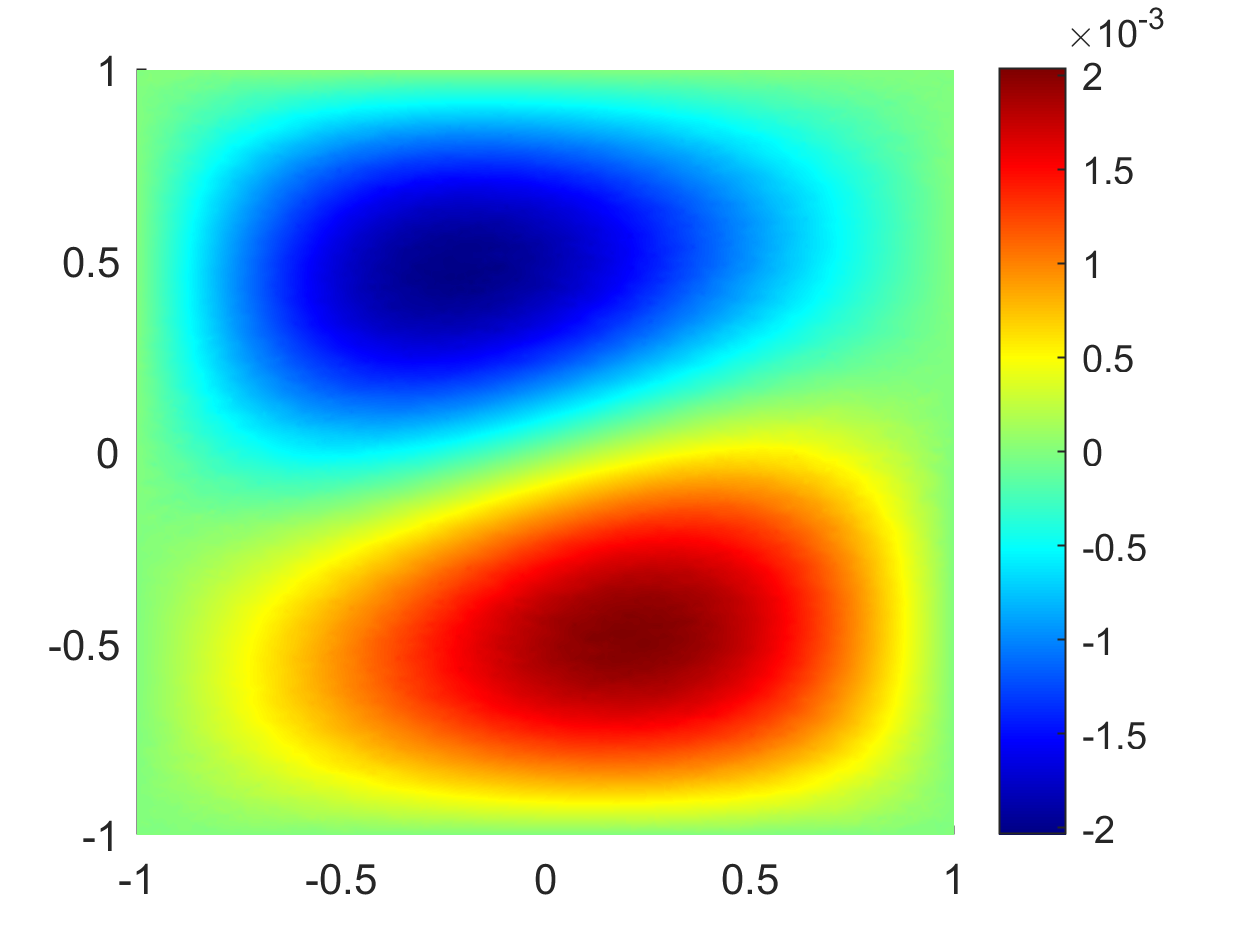}
    \caption{Mat\'{e}rn 3/2 Kernel}
  \end{subfigure}
  \begin{subfigure}{0.35\textwidth}
   \includegraphics[height=4cm,width=5cm]{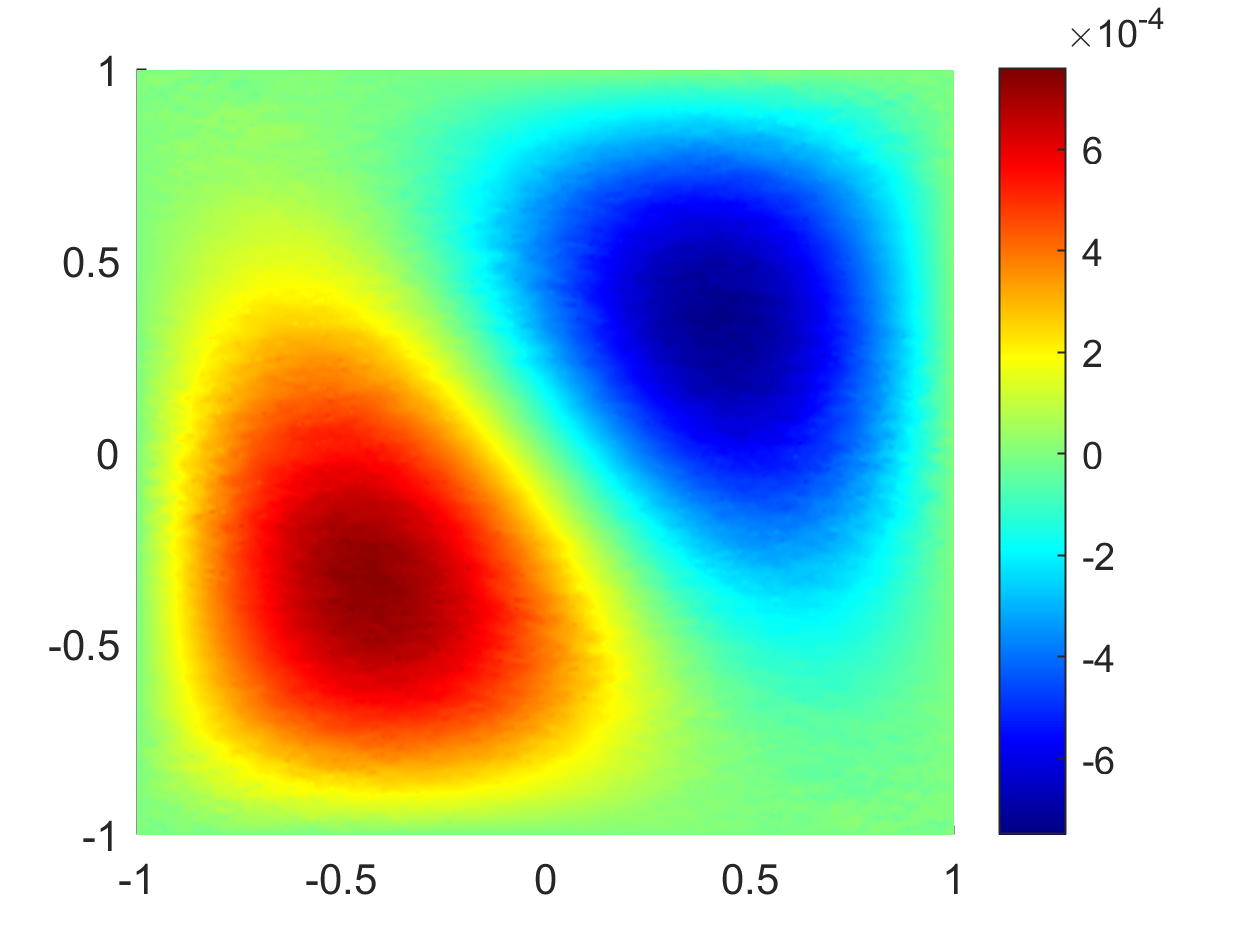}
    \caption{ Err Matrn 5/2 Kernel}
  \end{subfigure}
   \begin{subfigure}{0.35\textwidth}
   \includegraphics[height=4cm,width=5cm]{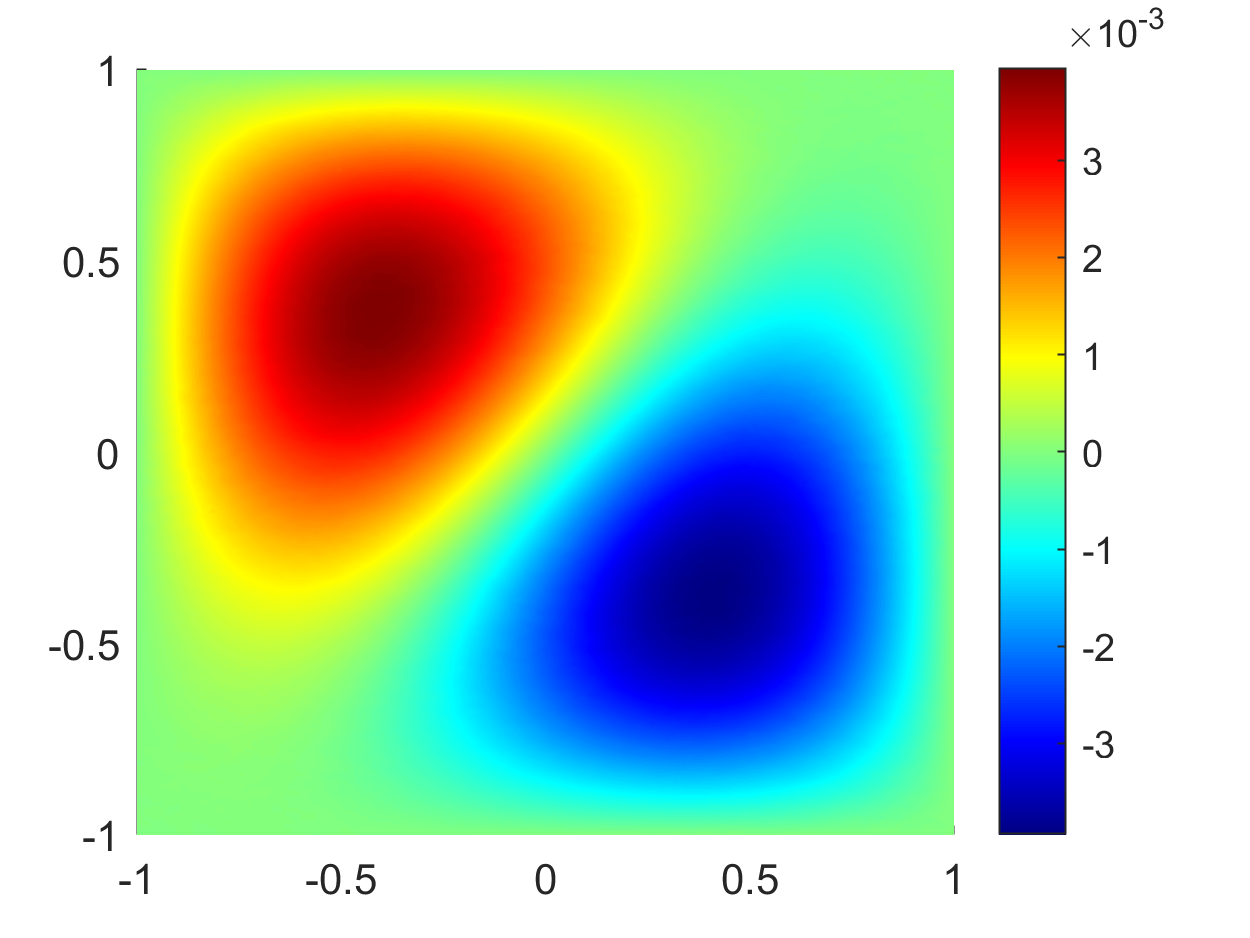}
    \caption{SE Kernel}
  \end{subfigure}
  \caption{Error between the FEM and GPR based 2nd eigenvectors of EVP~\eqref{eq:crs} using different kernels at ${\mu}=-0.75$.}
  \label{crs:evct2}
\end{figure}

\begin{figure}
\centering
 \begin{subfigure}{0.35\textwidth}
   \centering
   \includegraphics[height=4cm,width=5cm]{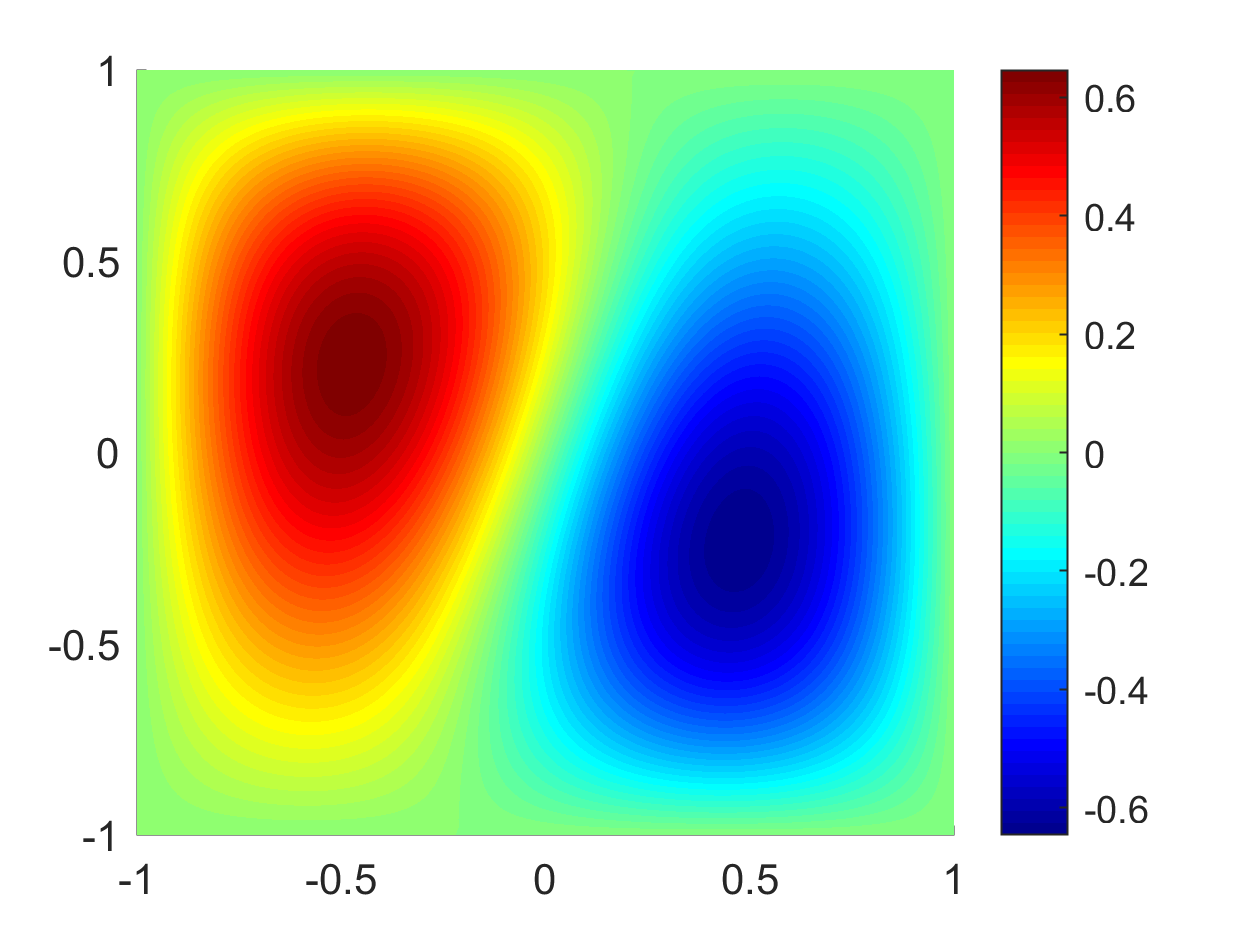}
    \caption{Exp Kernel}
  \end{subfigure}
  \begin{subfigure}{0.35\textwidth}
   \centering
   \includegraphics[height=4cm,width=5cm]{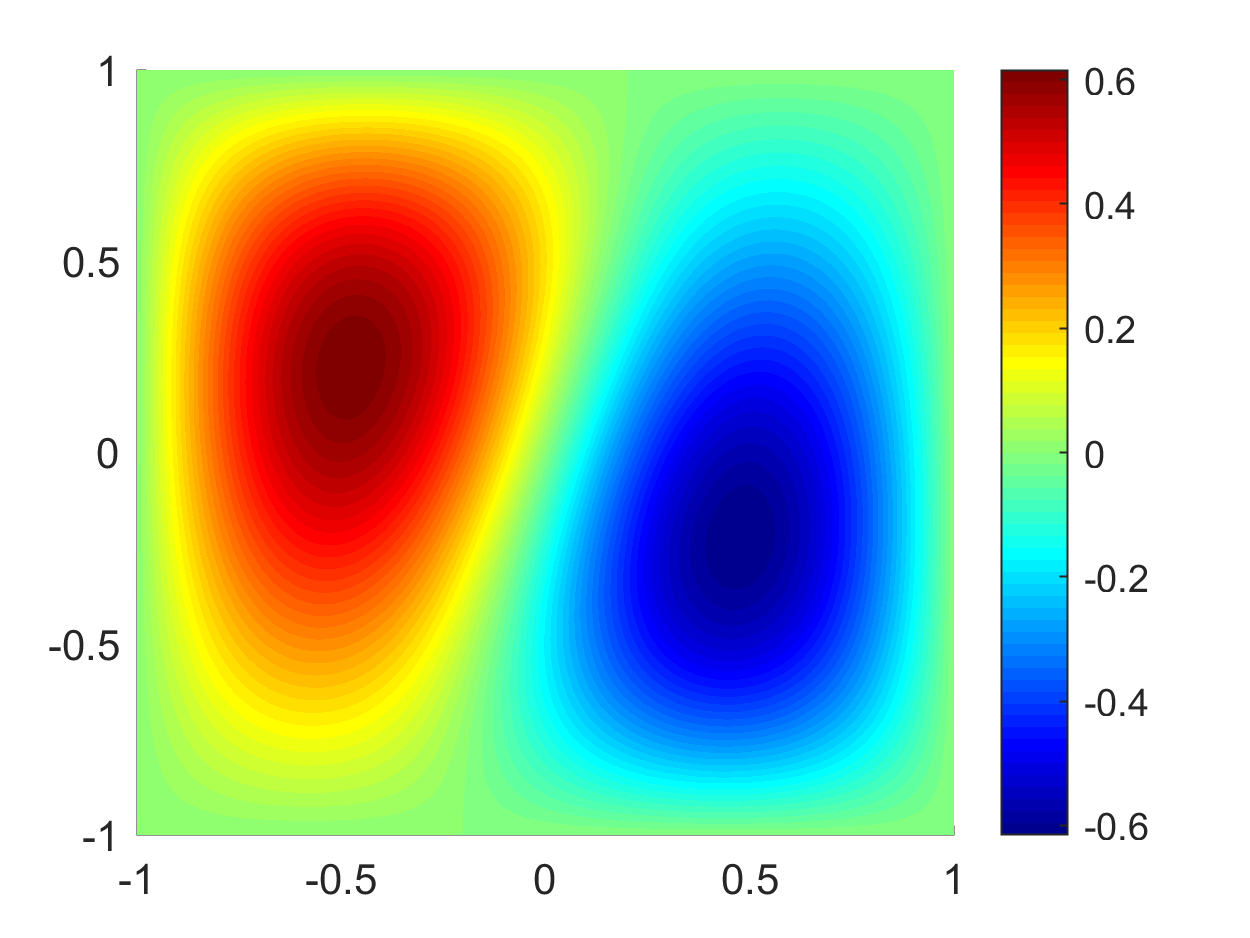}
    \caption{Mat\'{e}rn 3/2 Kernel}
  \end{subfigure}
  \begin{subfigure}{0.35\textwidth}
   \centering
   \includegraphics[height=4cm,width=5cm]{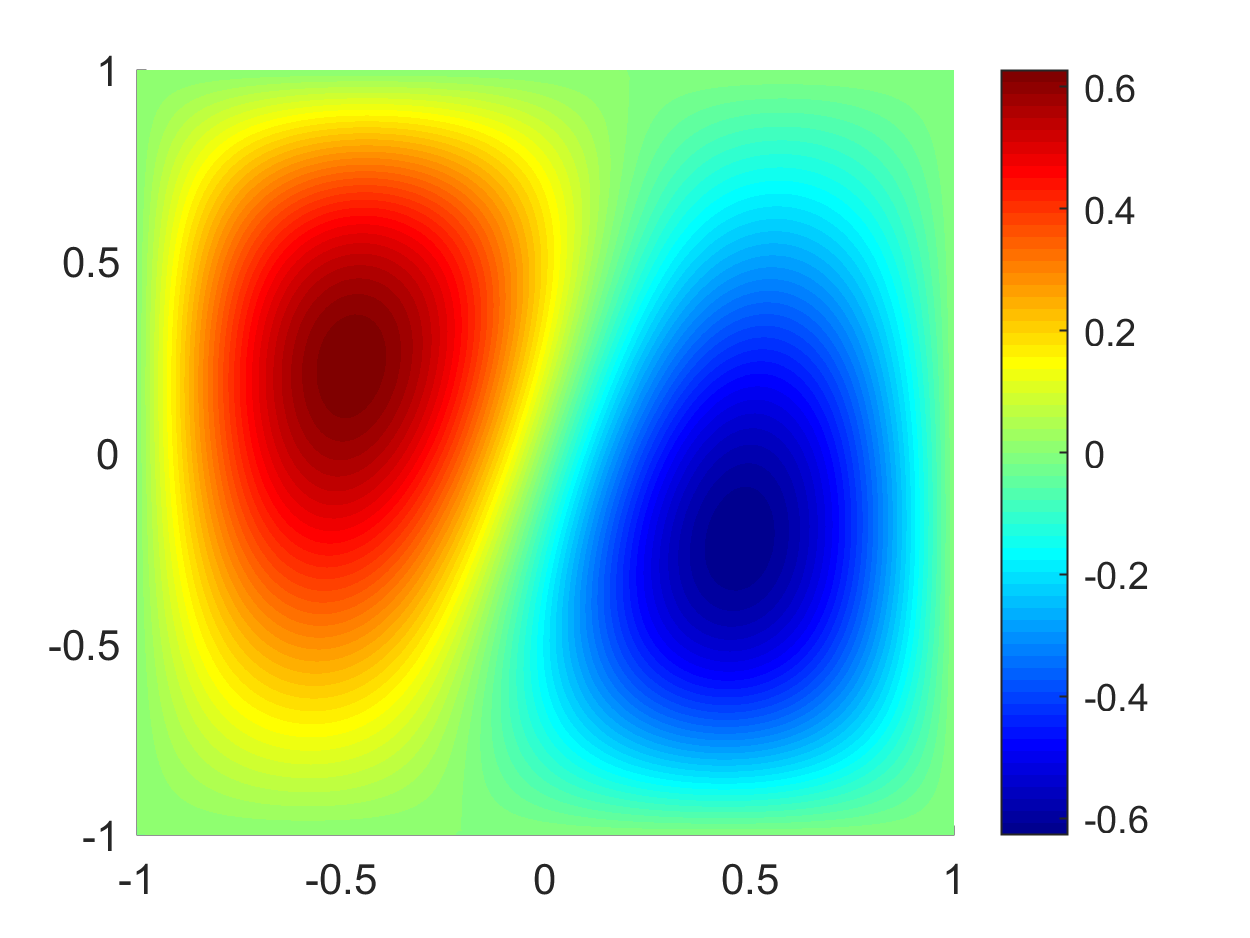}
    \caption{Mat\'{e}rn 5/2 Kernel}
  \end{subfigure}
  \begin{subfigure}{0.35\textwidth}
   \centering
   \includegraphics[height=4cm,width=5cm]{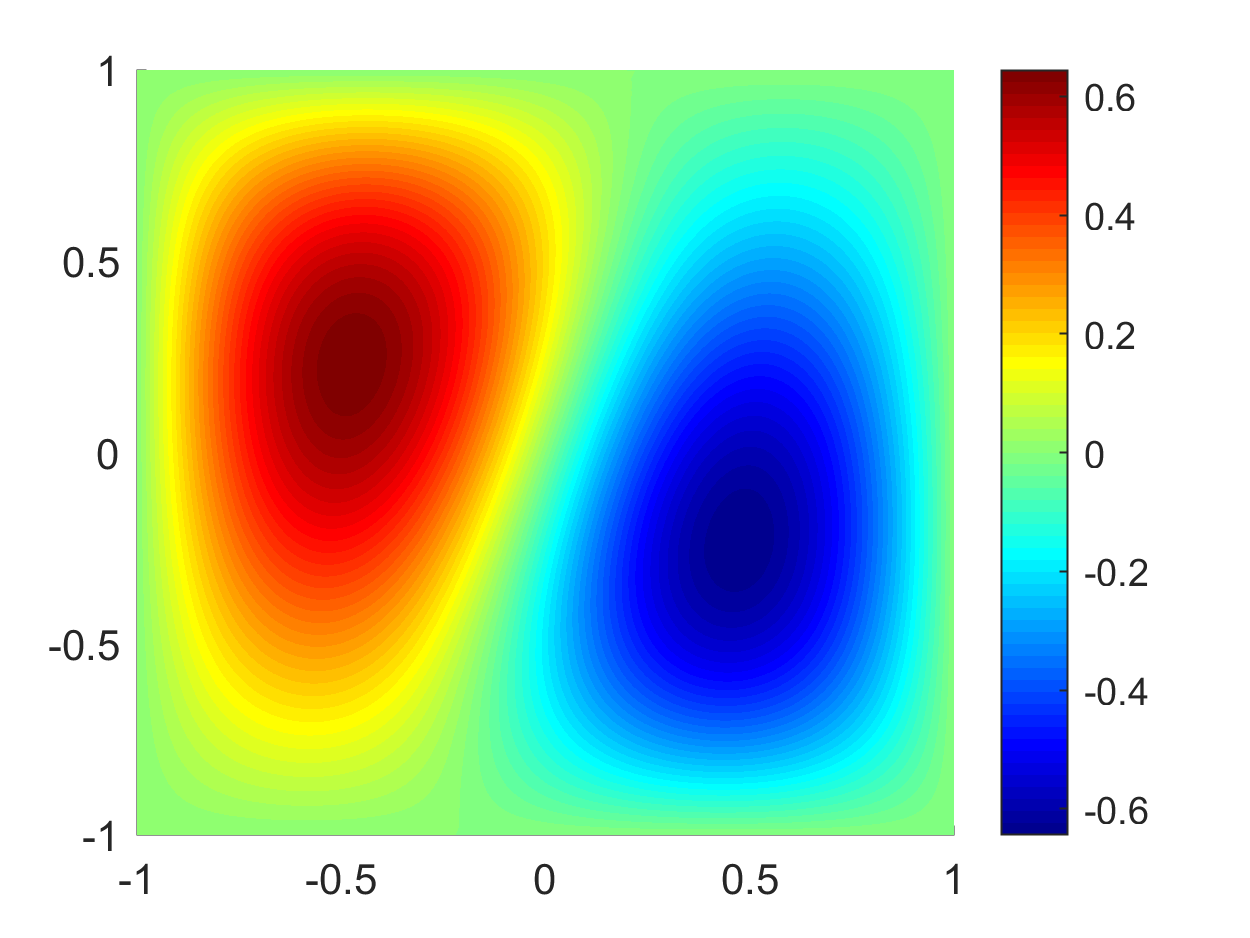}
    \caption{SE Kernel}
  \end{subfigure}
  \caption{Error between the second eigenvector and GPR result for second eigenvector of EVP~\eqref{eq:crs} with exp, Mat\'{e}rn 3/2, Mat\'{e}rn 5/2 and SE kernel at $\mu=-0.05$ when the training sample is $-0.9:0.1:0.9$.}
  \label{nw:fig2}
\end{figure}

In Fig.~\ref{crs:evct2}, we {report} the error between the true second eigenvector (obtained by FEM) and the GPR predictions at $\mu=-0.75$. In this case, the error for the exponential one is of order $10^{-4}$ whereas the error is of order $10^{-3}$ for the squared exponential kernel. Recall that $\mu = 0$ is a corner point. Thus, it is interesting to consider the performance of the four different models near this point. In Fig.~\ref{nw:fig2}, we {have shown} the RRMSE values at $\mu=-0.05$. Note the error is almost identical for all the GPR with different covariance functions at $\mu=-0.05$. This is because we only consider a few samples in the training set. 

\begin{figure}
\centering
 \begin{subfigure}{0.35\textwidth}
   \centering
   \includegraphics[height=4cm,width=5cm]{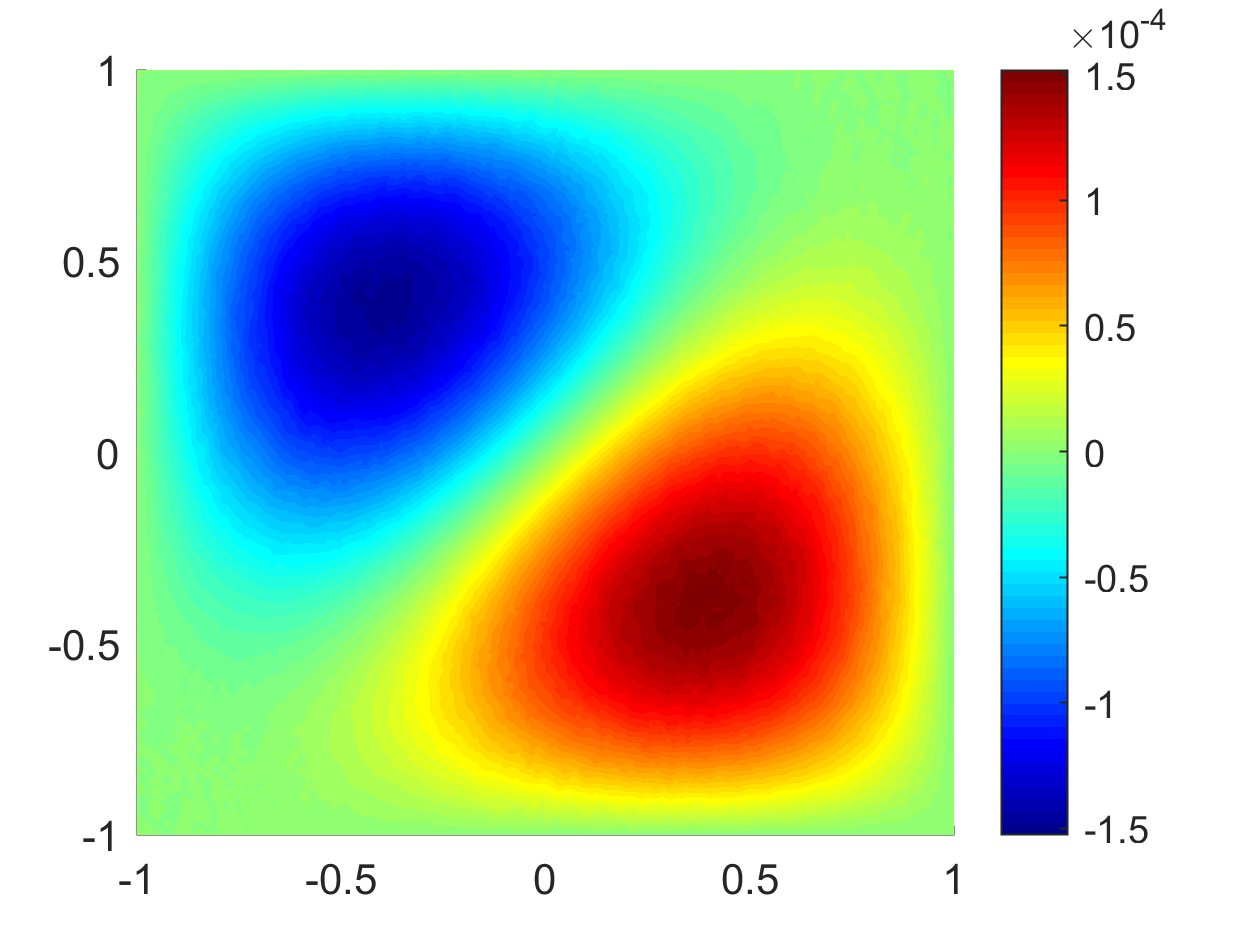}
    \caption{Exp Kernel}
  \end{subfigure}
  \begin{subfigure}{0.35\textwidth}
   \centering
   \includegraphics[height=4cm,width=5cm]{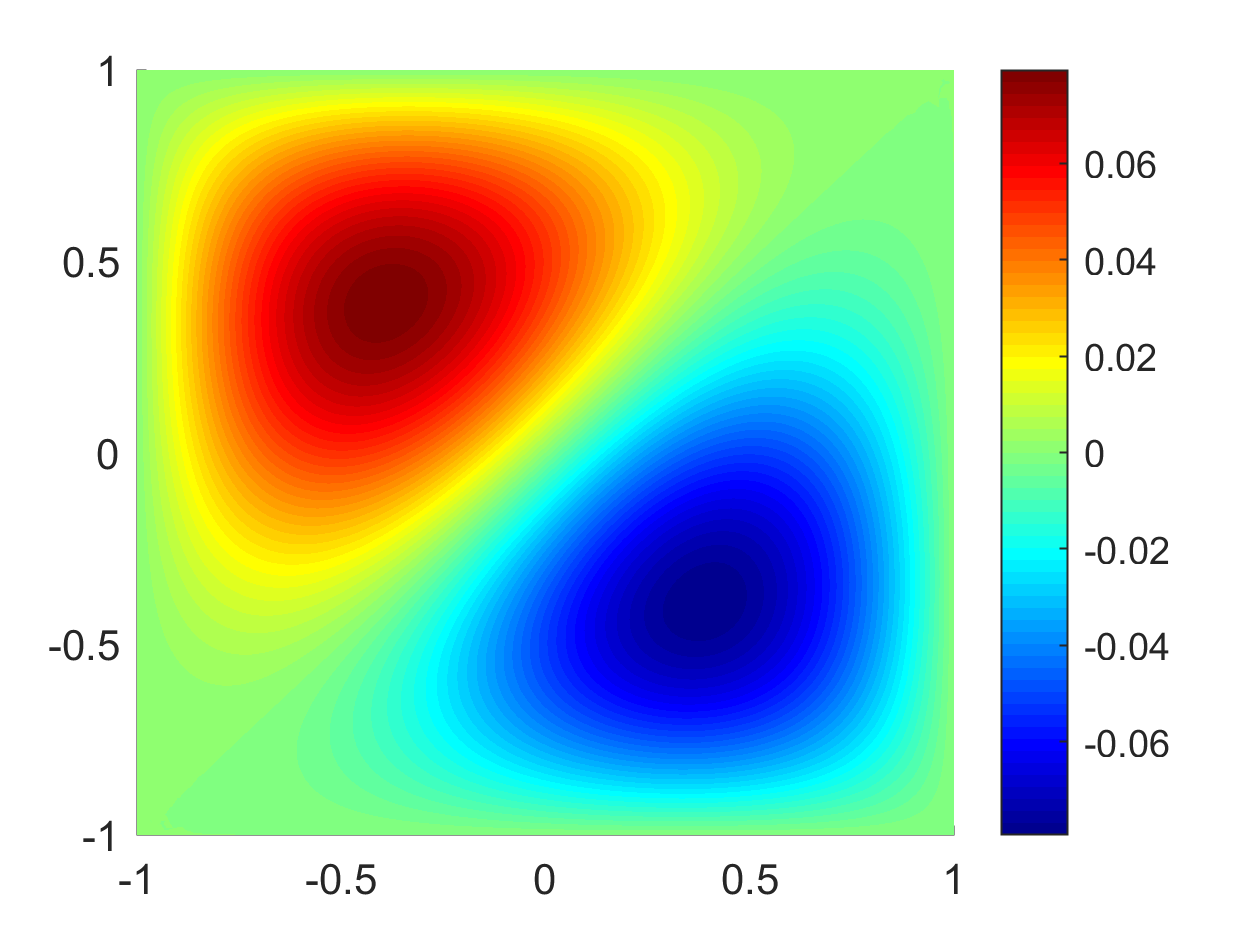}
    \caption{Mat\'{e}rn 3/2 Kernel}
  \end{subfigure}
  \begin{subfigure}{0.35\textwidth}
   \centering
   \includegraphics[height=4cm,width=5cm]{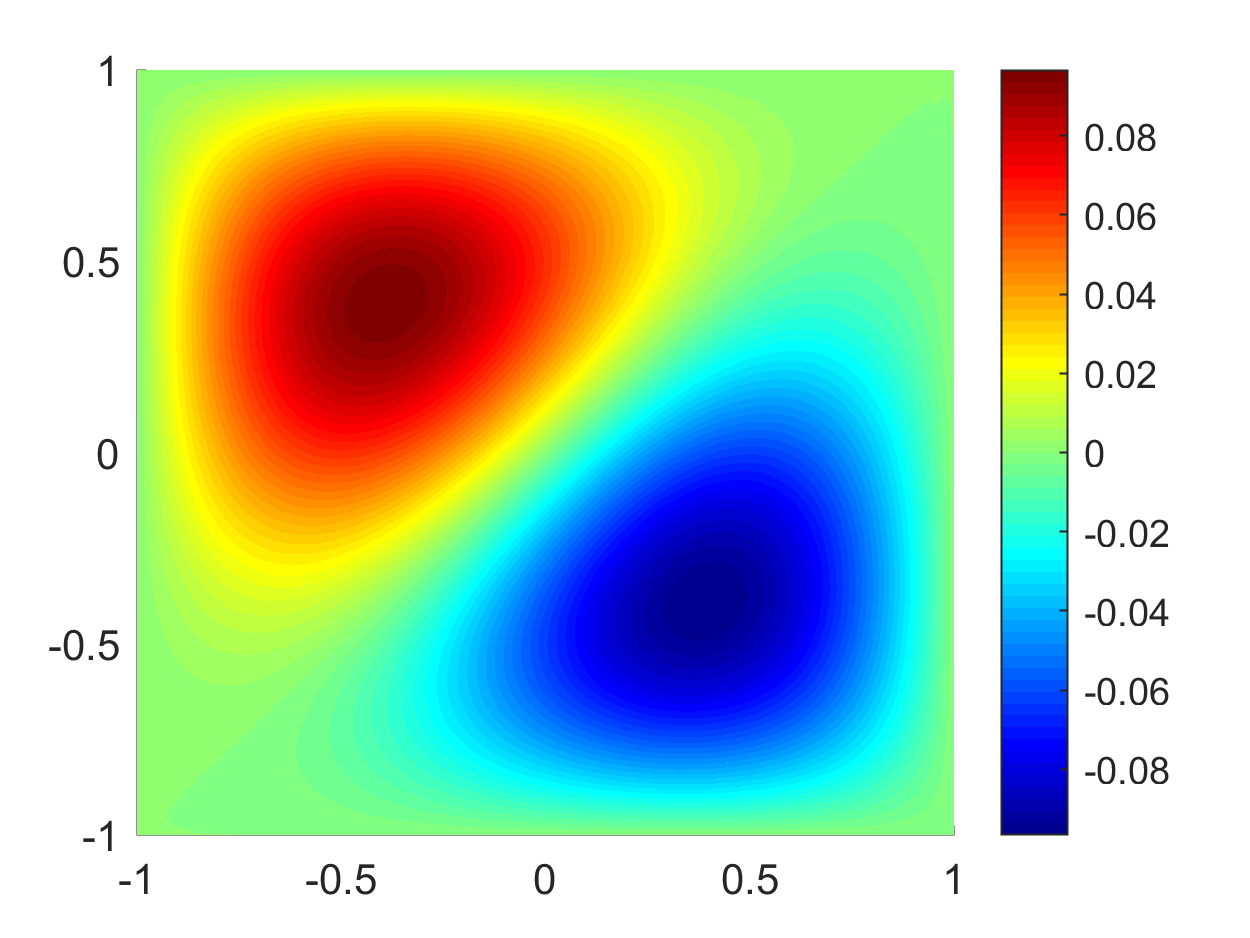}
    \caption{Mat\'{e}rn 5/2 Kernel}
  \end{subfigure}
  \begin{subfigure}{0.35\textwidth}
   \centering
   \includegraphics[height=4cm,width=5cm]{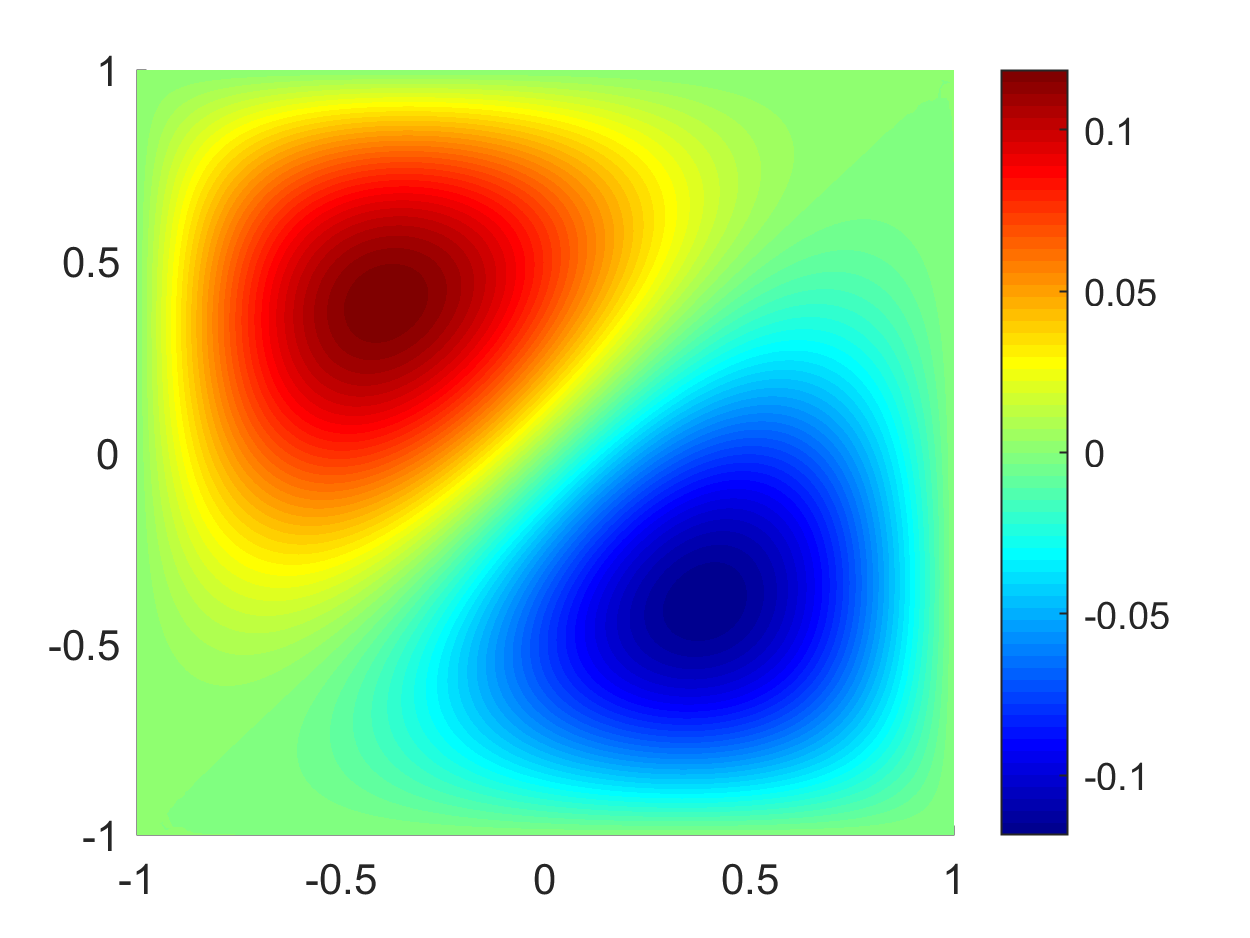}
    \caption{SE Kernel}
  \end{subfigure}
  \caption{Error between the second eigenvector using FEM and GPR of EVP~\eqref{eq:crs} with kernels exp, Mat\'{e}rn 3/2, Mat\'{e}rn 5/2 and SE kernel at $\mu=-0.05$ when sample points are $-0.9:0.04:0.9$.}
  \label{nw:fig4}
\end{figure}

Next, we consider more sample points. In particular, we choose our training set to be $\bm{\mu_{tr}} = -0.9:0.04:0.9$ and use the corresponding solutions as snapshots. The error for the second eigenvector corresponding to the new samples {is} shown in Fig.~\ref{nw:fig4}. In Fig.~\ref{nw:fig4}, one can see that the GPR with an absolute exponential kernel produces the best result among the four GPR models. {The point-wise error for the absolute kernel is of order $10^{-4}$, and of order $10^{-2}$ for  Mat\'{e}rn~3/2, Mat\'{e}rn~5/2 kernels. In Table~\ref{crs:ev2n}, we report the optimal value of the three hyper-parameters of the GPR corresponding to the second eigenvalue and the reduced coefficient of the second eigenvectors. All the parameters lie within the default range $[10^{-3}, 10^3]$ and are well separated from the endpoints of the interval.}

\begin{table}
 \centering
  \begin{tabular}{|c|c|c|c|c|c|c|c|c|} 
 Cov &   Method  &$\mu=-0.75$  &$\mu=-0.25$  & $\mu=0.25$& $\mu=0.75$&RMSRE \\
 \hline
 &        FEM&4.93615339& 9.87196392& 12.95691963& 14.19104350 &\\
 \hline
Exp      &GPR&4.93620186& 9.87195587& 12.95687818& 14.19098776 & ${7.3\times 10^{-7}}$ \\% 
& Rel. Err. & $9.8 \times 10^{-6}$& $8.1 \times 10^{-7}$&$3.1 \times 10^{-6}$&$3.9 \times 10^{-6}$&  \\
\hline
Mat\'{e}rn 3/2   &GPR  & 4.93528087 &9.87007887& 12.95503363 &14.19140987 &{ $1.3 \times 10^{-4}$} \\
& Rel. Err. & $1.7 \times 10^{-4}$& $1.9 \times 10^{-4}$&$1.4 \times 10^{-4}$&$2.5 \times 10^{-5}$& \\
\hline
Mat\'{e}rn 5/2   & GPR &4.92953921 &9.88022559& 12.96546014& 14.19445359& {$1.7 \times 10^{-4}$} \\%& $\re{9.6 \times 10^{-6}}$\\
& Rel. Err. & $1.3 \times 10^{-3}$& $8.3 \times 10^{-4}$&$6.5 \times 10^{-4}$&$2.4 \times 10^{-4}$& \\
\hline
SE &GPR& 4.94133382 &9.94586709& 13.04302796 &14.23611366 &{$3.0 \times 10^{-4}$} \\%&$\re{3.3 \times 10^{-5}}$\\
& Rel. Err. & $1.0 \times 10^{-3}$& $7.5 \times 10^{-3}$&$6.6 \times 10^{-3}$&$3.2 \times 10^{-4}$& \\
 \end{tabular}
\caption{{The values} of RMSRE corresponding to 2nd eigenvalues of EVP~\eqref{eq:crs} using GPR model with different covariance functions.}
	 	\label{crs:ev2}
 \end{table}

\begin{table}
 \centering
  \begin{tabular}{|c|c|c|c|c|c|c|c|c|} 
$\delta \mu$& Exponential  &Mat\'{e}rn 3/2 &Mat\'{e}rn 5/2  & Squared Exp.\\
 \hline
 %EV2 & $0.1$ & $7.3\times 10^{-7}$ & $1.3\times 10^{-4}$ & $1.7\times 10^{-4}$ & $3.0\times 10^{-5}$\\
      $0.06$& $6.0\times 10^{-8}$ & $7.9\times 10^{-5}$ & $9.0\times 10^{-5}$ & $1.1\times 10^{-4}$\\
      $0.04$& $1.5\times 10^{-8}$ & $1.3\times 10^{-5}$ & $1.3\times 10^{-5}$ & $2.3\times 10^{-5}$\\
    $0.02$& $6.3\times 10^{-10}$ & $6.8\times 10^{-6}$ & $1.0\times 10^{-5}$ & $1.4\times 10^{-5}$\\
 \end{tabular}
\caption{{The values of RMSRE for the second eigenvalue of the EVP~\eqref{eq:crs} using GPR model with different covariance functions when training samples are $-0.9:\delta \mu:0.9$.}}
\label{crs:ev2_rrse}
 \end{table}

\begin{table}
 \centering
 \begin{tabular}{|c|c|c|c|c|c|c|c|} 
 Name &Covariance  & $\ell$& $\sigma_f$ & $\beta_0$ & Log-Likelihood\\
 \hline
 2nd eigenvalue
     &Exponential & 10.90 &  5.33 &9.13 &  -22.68 \\%&$9.3\times 10^{-11}$ \\
     & Matern 3/2 & 3.17  &  9.76 &6.67 &  0.01 \\%&$9.6 \times 10^{-6}$  \\
     & Matern 5/2 & 0.69  &  4.34 &9.35 & -4.84 \\%&$1.2\times 10^{-5}$ \\
     &Squared Exp & 0.15  &  2.73 &10.19& -24.80 \\%&$3.5\times 10^{-5}$   \\
\hline
 coefficient-1  
      & Exponential&0.69 &  19.57 &  43.05 & -71.01\\%&    0.1764\\
     & Matern 3/2 & 0.19 &  17.99 &  43.38  & -72.75\\% &    0.1909 \\   
     & Matern 5/2 & 0.14 &  17.52 &  43.42 &-73.73\\% &    0.1918\\  
     &Squared Exp & 0.09 &  16.30 &  43.46 & -75.71\\ %&    0.1925\\
 \hline
coefficient-2  
      & Exponential&1.31 &  41.69 &  -19.10 & -80.20\\ %&   0.0853 \\
     & Matern 3/2 &  0.30 &  38.15 &  -18.69 & -79.85\\%&   0.0870 \\  
     & Matern 5/2 & 0.20  &  36.05 & -18.61  & -81.14\\% &   0.0866\\ 
    &Squared Exp &0.11   &  33.53  &  -18.71 & -84.96\\% &   0.0861\\
 \end{tabular}
 \caption{{Optimum hyperparameters corresponding to the second eigenpair of the EVP~\eqref{eq:crs} using different kernels while. We fixed the noise SD $\sigma_n=0.0001$ and the training points as $-0.9:0.1:0.9$.}}
\label{crs:ev2n}
 \end{table}

% \begin{figure}
%  \begin{subfigure}{0.23\textwidth}
%    \includegraphics[height=3cm,width=4cm]{rev_result/evct2_ex_m05.png}
%     \caption{Exp Kernel}
%   \end{subfigure}
%   \begin{subfigure}{0.23\textwidth}
%    \includegraphics[height=3cm,width=4cm]{rev_result/evct2_ma3_m05.png}
%     \caption{Mat\'{e}rn 3/2 Kernel}
%   \end{subfigure}
%   \begin{subfigure}{0.23\textwidth}
%    \includegraphics[height=3cm,width=4cm]{rev_result/evct2_ma5_m05.png}
%     \caption{Mat\'{e}rn 5/2 Kernel}
%   \end{subfigure}
%   \begin{subfigure}{0.23\textwidth}
%    \includegraphics[height=3cm,width=4cm]{rev_result/evct2_se_m05.png}
%     \caption{SE Kernel}
%   \end{subfigure}
%   \caption{\re{Second eigenvector and GPR result for second eigenvector with exp, Mat\'{e}rn 3/2, Mat\'{e}rn 5/2 and SE kernel at $\mu=-0.05$ when sample points are $-0.9:0.04:0.9$.}}
%   \label{nw:fig3}
% \end{figure}

\begin{figure}
\centering
 \begin{subfigure}{0.35\textwidth}
   \centering
   \includegraphics[height=4cm,width=5cm]{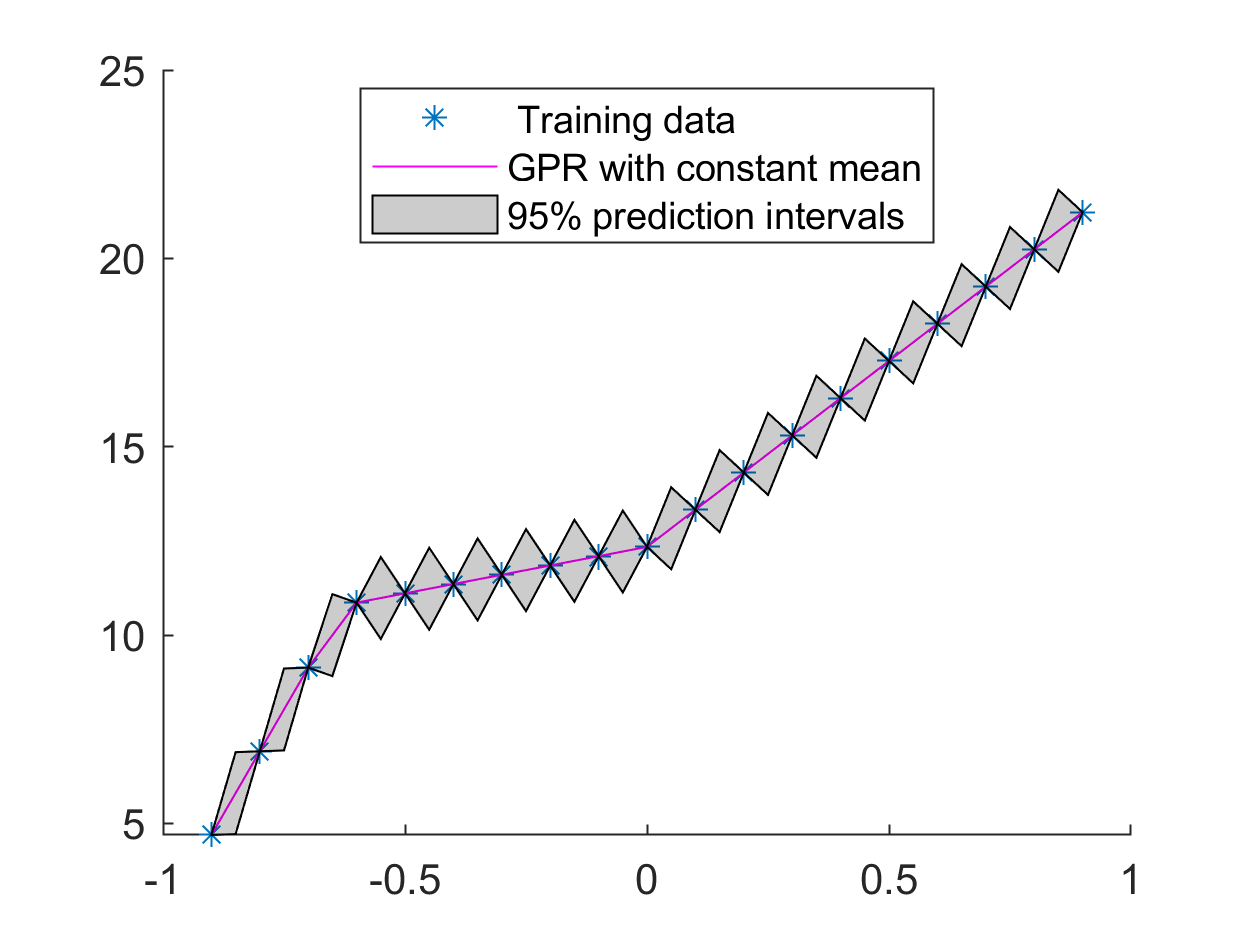}
    \caption{Exp Kernel}
  \end{subfigure}
 \begin{subfigure}{0.35\textwidth}
   \centering
   \includegraphics[height=4cm,width=5cm]{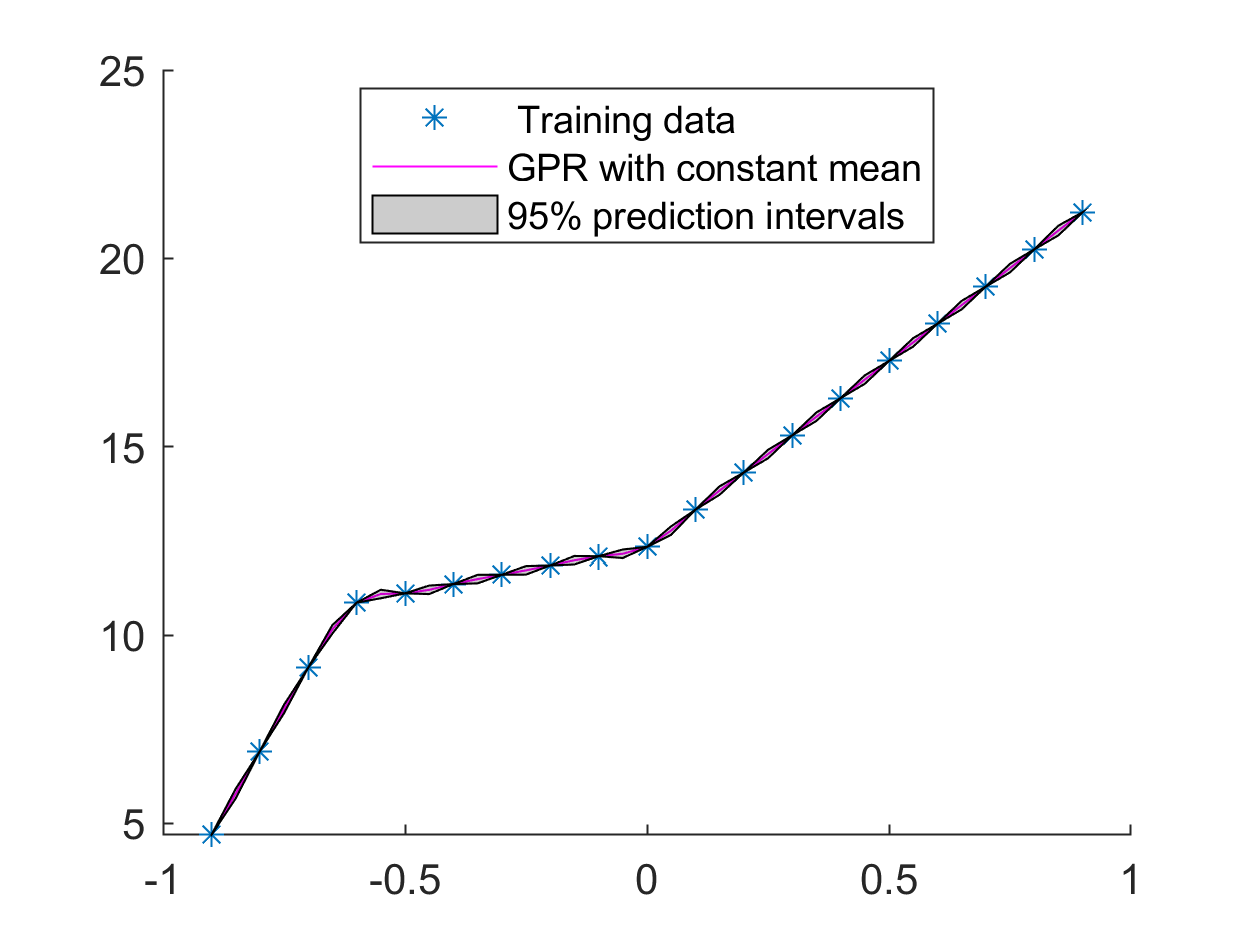}
    \caption{Mat\'{e}rn 3/2 Kernel}
  \end{subfigure}
  \begin{subfigure}{0.35\textwidth}
   \centering
   \includegraphics[height=4cm,width=5cm]{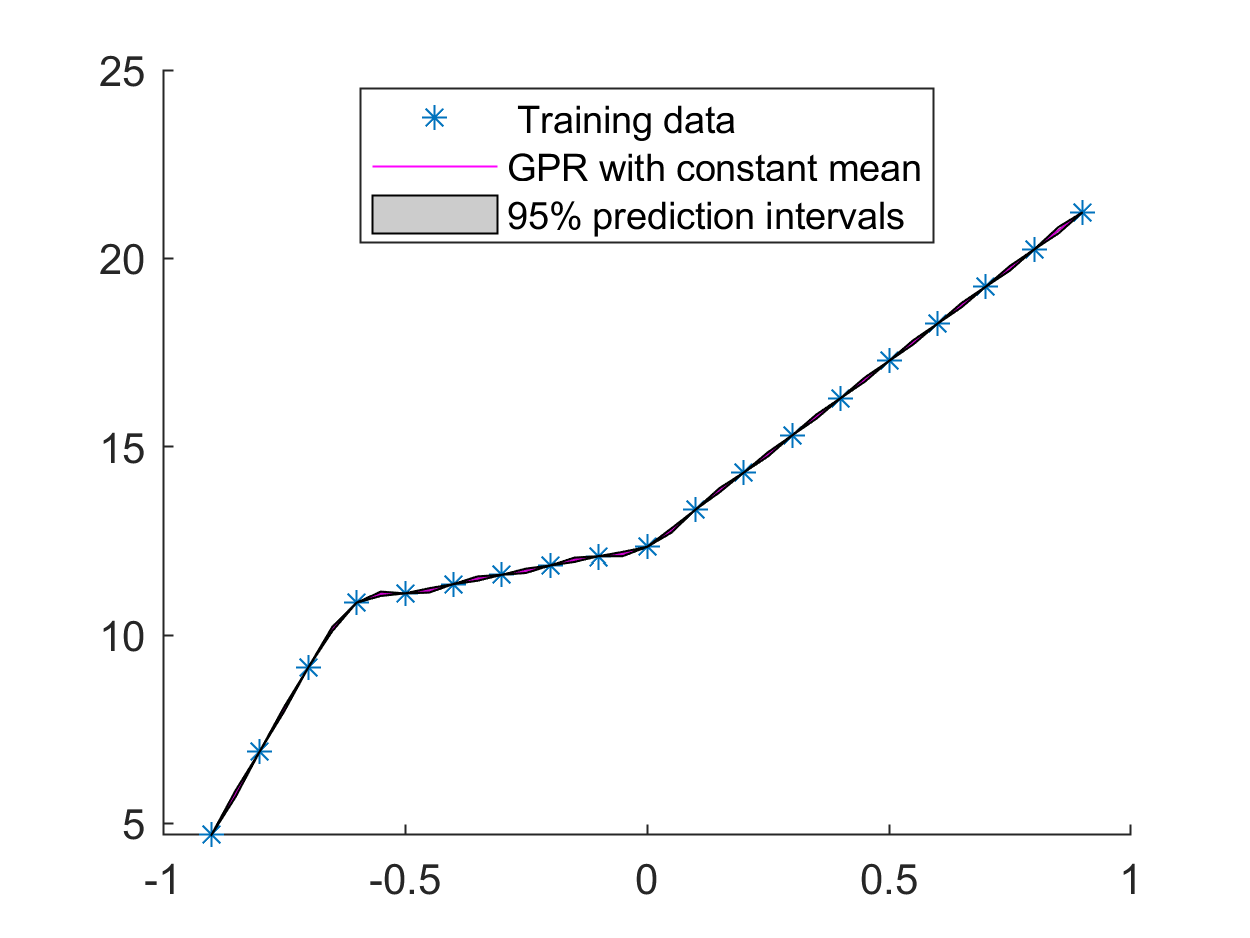}
    \caption{Mat\'{e}rn 5/2 Kernel}
  \end{subfigure}
  \begin{subfigure}{0.35\textwidth}
   \centering
   \includegraphics[height=4cm,width=5cm]{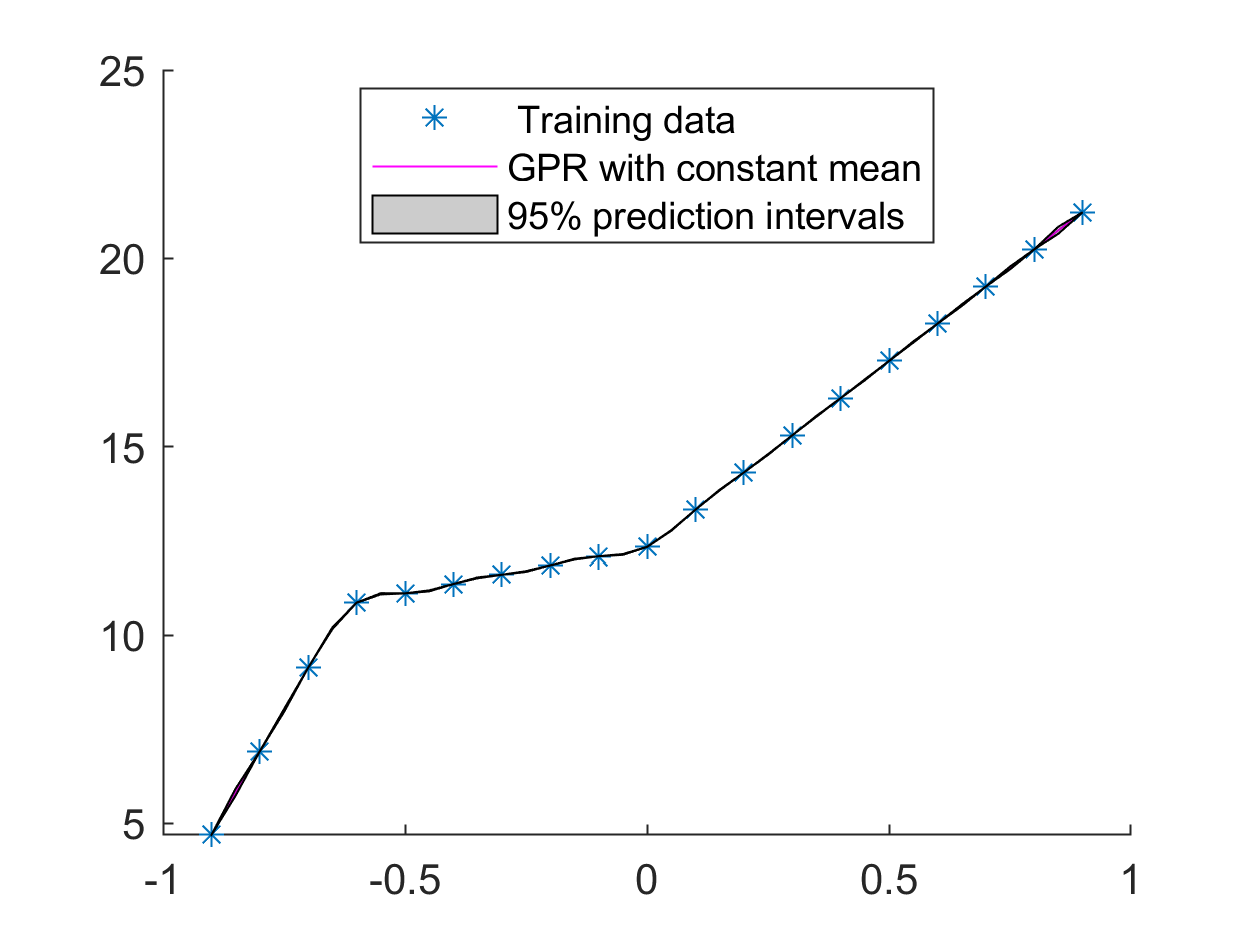}
    \caption{SE Kernel}
  \end{subfigure}
  \caption{GPR corresponding to 3rd eigenvalues of EVP~\eqref{eq:crs} using different kernels.}
  \label{ev3:crs}
\end{figure}
%%%%%%%%%%%%%%%%%%%%%%%%%%%% NEW
% \begin{table}
%  \centering
%   \begin{tabular}{|c|c|c|c|c|c|c|}  
%  Cov &   Method  &$\mu=-0.75$  &$\mu=-0.25$  & $\mu=0.25$& $\mu=0.75$&RRMSE\\
%  \hline
%  &  FEM&8.02390130 &11.72278922&14.80775895 &19.74354759 &\\
%  \hline
% Exp      &DD& 8.02357199& 11.72280632& 14.80773500& 19.74345655 & $6.4 \times 10^{-4}$\\
% %$\re{1.0 \times 10^{-3}}$\\
% & Rel. Err. & $4.1 \times 10^{-5}$& $1.4 \times 10^{-6}$&$1.6 \times 10^{-6}$&$4.6 \times 10^{-6}$ & \\
% \hline
% Mat\'{e}rn 3/2   &  &8.04670722& 11.72157122& 14.80677595 &19.74230459&$3.7 \times 10^{-4}$ \\
% %$\re{6.7 \times 10^{-4}}$\\
% & Rel. Err. & $2.8 \times 10^{-3}$& $1.0 \times 10^{-4}$&$6.6 \times 10^{-5}$&$6.2 \times 10^{-5}$& \\
% \hline
% Mat\'{e}rn 5/2   &  &8.08692759& 11.72891380& 14.81734316& 19.74572554 &$5.7 \times 10^{-4}$ \\
% %\re{6.9 \times 10^{-4}}$\\
% & Rel. Err. & $7.9 \times 10^{-3}$& $5.2 \times 10^{-4}$&$6.4 \times 10^{-4}$&$1.1 \times 10^{-4}$& \\
% \hline
% SE       &DD& 8.18813426& 11.65039408& 14.76794206& 19.74988015  &$1.5 \times 10^{-3}$ \\
% %$\re{8.4 \times 10^{-4}}$ \\
% & Rel. Err. & $2.0 \times 10^{-2}$& $6.2 \times 10^{-3}$&$2.7 \times 10^{-3}$&$3.2 \times 10^{-4}$&\\
%  \end{tabular}
% \caption{Comparison of 3rd eigenvalues of EVP~\eqref{eq:crs} using GPR model with different covariance functions with training sample points $-0.9:0.1:0.9$}.
% \label{crs:ev3}
%  \end{table}
\begin{table}
 \centering
  \begin{tabular}{|c|c|c|c|c|c|c|c|c|} 
  $\delta \mu$& Exponential  &Mat\'{e}rn 3/2 &Mat\'{e}rn 5/2  & Squared Exp.\\
    \hline
     $0.1$ & $1.7\times 10^{-4}$ & $4.7\times 10^{-5}$ & $4.4\times 10^{-5}$ & $4.6\times 10^{-5}$\\
      $0.06$& $5.9\times 10^{-5}$ & $2.1\times 10^{-5}$ & $1.9\times 10^{-5}$ & $2.0\times 10^{-5}$\\
      $0.04$& $1.7\times 10^{-6}$ & $2.3\times 10^{-6}$ & $3.0\times 10^{-6}$ & $5.6\times 10^{-6}$\\
     $0.02$& $1.5\times 10^{-16}$ & $4.0\times 10^{-7}$ & $4.6\times 10^{-7}$ & $5.5\times 10^{-7}$\\
 \end{tabular}
\caption{{The values of RMSRE for the third eigenvalue of the EVP~\eqref{eq:crs} using GPR model with different covariance functions when training samples are $-0.9:\delta \mu:0.9$.}}
\label{crs:ev3_new}
 \end{table}

\begin{table}
 \centering
 \begin{tabular}{|c|c|c|c|c|c|c|c|} 
 Name &Covariance  & $\ell$& $\sigma_f$ & $\beta_0$& Likelihood\\
 \hline
 3rd eigenvalue
     &Exponential &  9.77&   7.76&  13.01 & -30.79 \\ %&$1.6\times 10^{-4}$ \\
     & Matern 3/2 &  2.50&   15.15&    9.83 & -14.43\\ % &$5.8 \times 10^{-5}$  \\
     & Matern 5/2 &  0.63 &  6.7601 &   12.52 & -16.28\\ % &$5.7\times 10^{-5}$ \\
     &Squared Exp &  0.16  &  4.01 & 13.24  & -29.78\\ %&$5.9\times 10^{-5}$   \\
\hline
 Coefficient-1  
      & Exponential&0.96 &   31.95 &   34.31 &-77.67 \\%&    0.1302\\
     & Matern 3/2 & 0.22  &  28.77 &   35.26 & -79.12 \\% &   0.1384 \\  
     & Matern 5/2 &0.16 &  27.72 &  35.38  & -80.34\\% &   0.1386\\   
     &Squared Exp &0.10 &  25.75 & 35.47   & -83.06\\% &    0.1386\\
 \hline %%%%%%%%%%%%%%%%%%
 Coefficient-2  
      & Exponential&0.41 &   30.46 &  -13.24 &-83.18 \\%&    0.1214 \\
     & Matern 3/2 &0.18  &  29.95 &  -16.2  & -82.94 \\%&   0.1221 \\   
     & Matern 5/2 & 0.14 &   29.19 &  -16.77 &-83.50 \\%&  0.1226\\  
    &Squared Exp & 0.10 & 29.25  &  -17.13  & -84.75 \\% & 0.1250\\
    \hline %%%%%%%%%%%%%%%%%%
 Coefficient-3   
      & Exponential&0.58 &   27.77 &  18.55 & -79.00\\%&    0.0575\
     & Matern 3/2 &  0.17 &   24.16 &   13.45 & -79.54\\ %&  0.0573 \\  
     & Matern 5/2 &0.13  & 23.30&  12.84   & -80.20 \\ %   0.0582\\   
    &Squared Exp & 0.09 &   22.42 &   12.17 & -81.61  \\% &    0.0601\\
 \end{tabular}
 \caption{{Optimum hyperparameters corresponding to the third eigen-pair of the EVP~\eqref{eq:crs} using different kernels while We fixed the noise SD $\sigma_n=0.0001$ and Training points are $-0.9:0.1:0.9$.}}
\label{crs:ev3n}
 \end{table}

\begin{figure}
\centering
 \begin{subfigure}{0.35\textwidth}
   \centering
   \includegraphics[height=4cm,width=5cm]{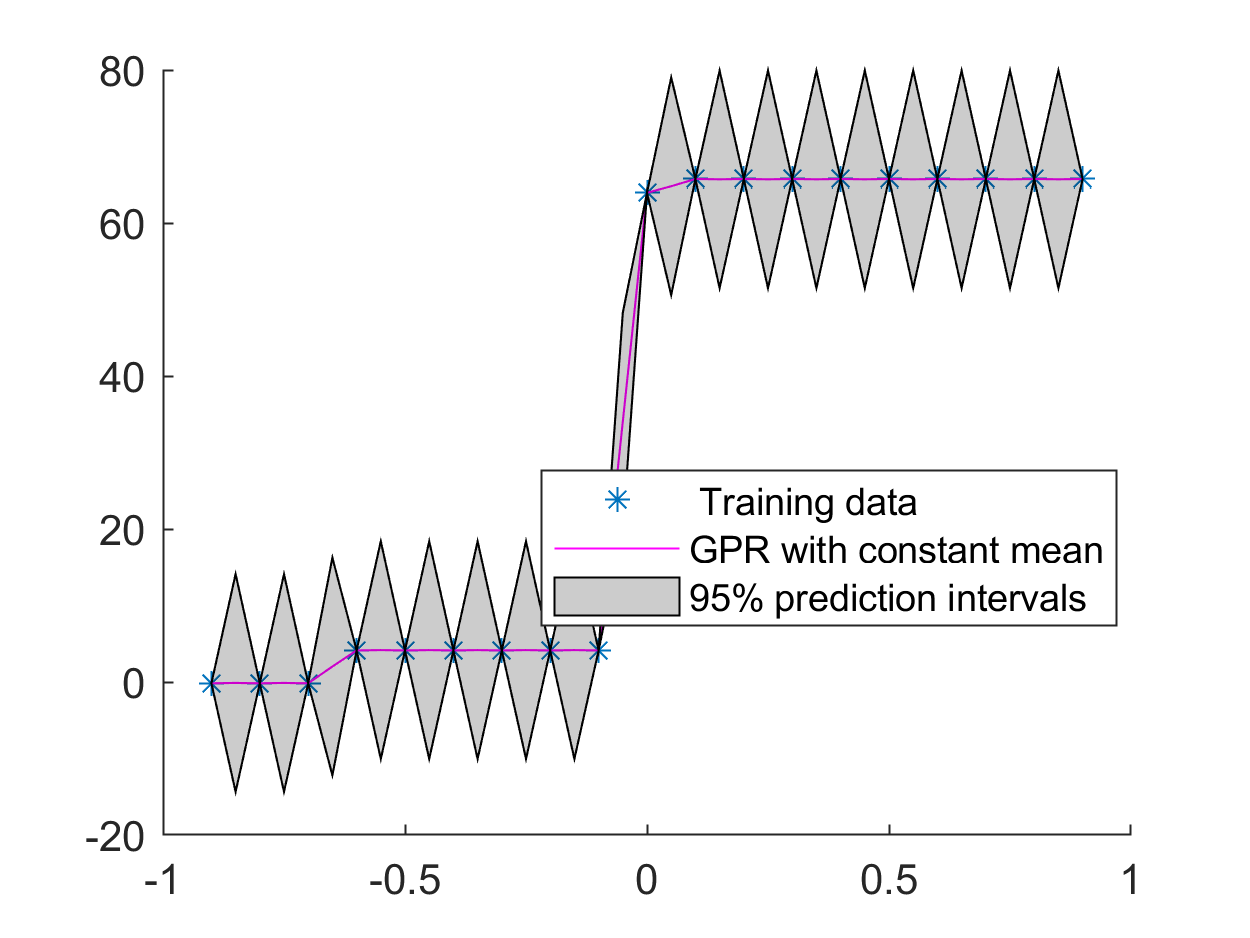}
    \caption{Exp Kernel}
  \end{subfigure}
 \begin{subfigure}{0.35\textwidth}
   \centering
   \includegraphics[height=4cm,width=5cm]{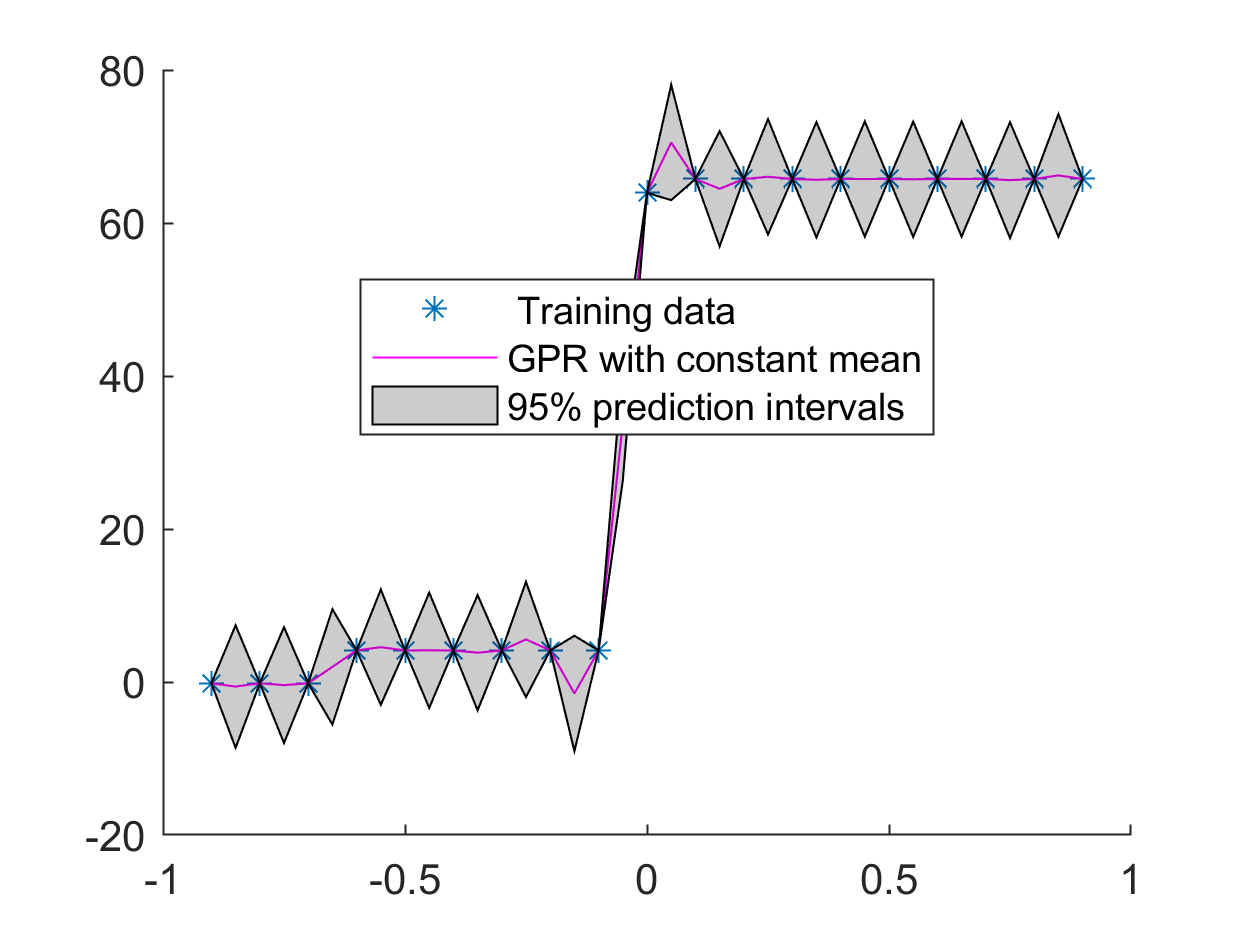}
    \caption{Mat\'{e}rn 3/2 Kernel}
  \end{subfigure}
  \begin{subfigure}{0.35\textwidth}
   \centering
   \includegraphics[height=4cm,width=5cm]{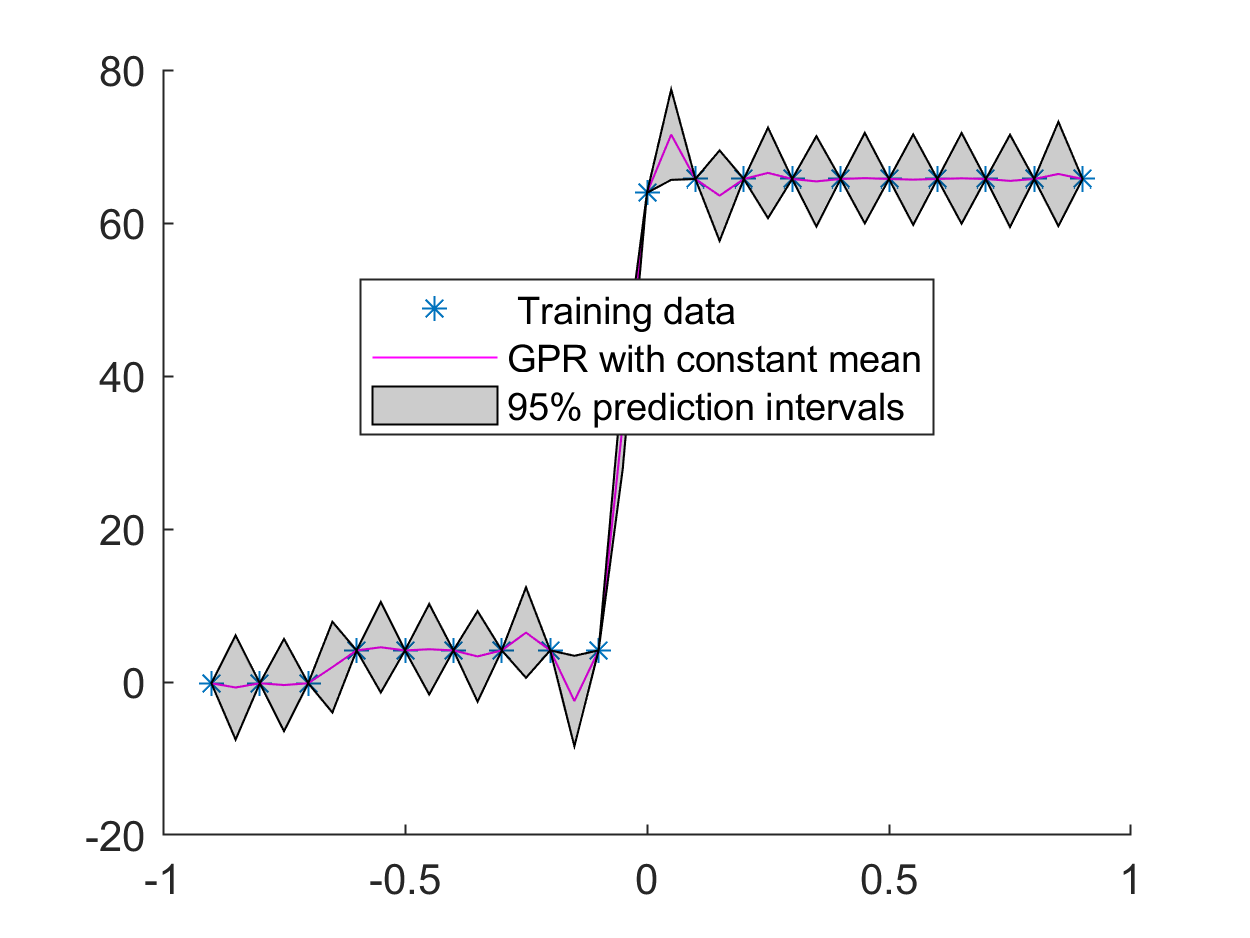}
    \caption{Mat\'{e}rn 5/2 Kernel}
  \end{subfigure}
  \begin{subfigure}{0.35\textwidth}
   \centering
   \includegraphics[height=4cm,width=5cm]{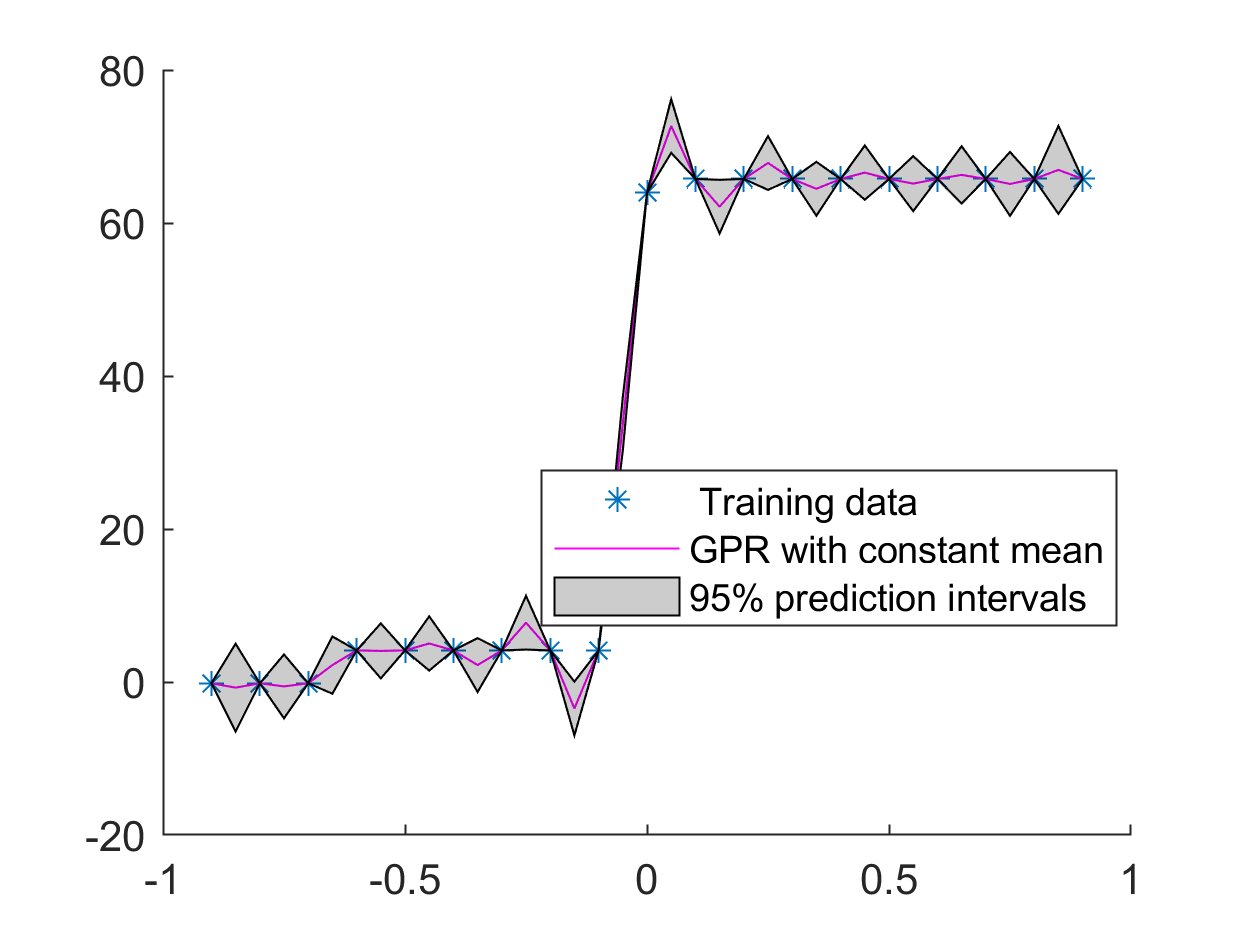}
    \caption{SE Kernel}
  \end{subfigure}
  \caption{GPR corresponding to 1st coefficient of the reduced 3rd eigenvector of EVP~\eqref{eq:crs} using different kernels and the training points are $-0.9:0.1:0.9$.}
  \label{crs:ev3_coeff1}
\end{figure}

\begin{figure}
\centering
 \begin{subfigure}{0.35\textwidth}
   \centering
   \includegraphics[height=4cm,width=5cm]{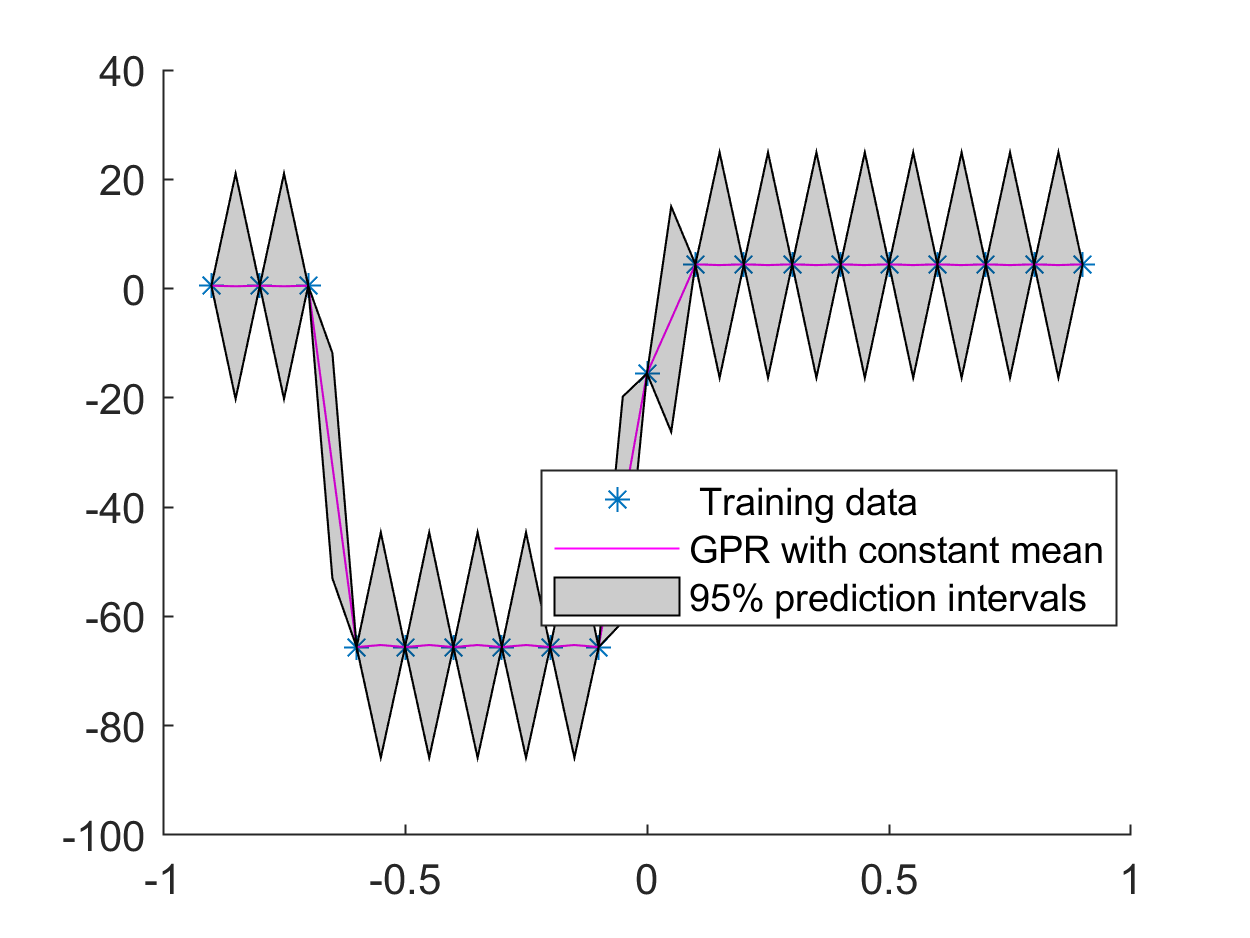}
    \caption{Exp Kernel}
  \end{subfigure}
 \begin{subfigure}{0.35\textwidth}
   \centering
   \includegraphics[height=4cm,width=5cm]{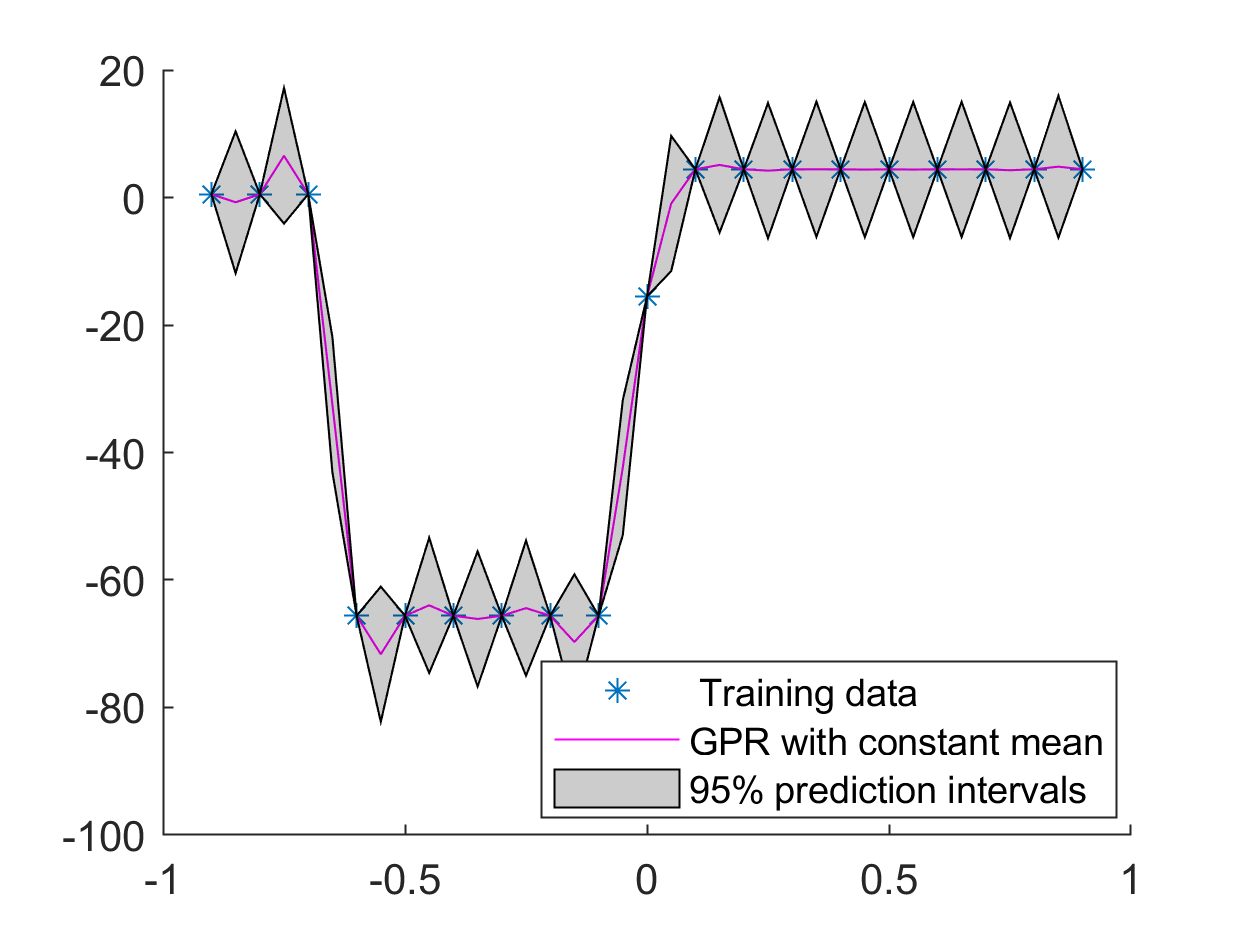}
    \caption{Mat\'{e}rn 3/2 Kernel}
  \end{subfigure}
  \begin{subfigure}{0.35\textwidth}
   \centering
   \includegraphics[height=4cm,width=5cm]{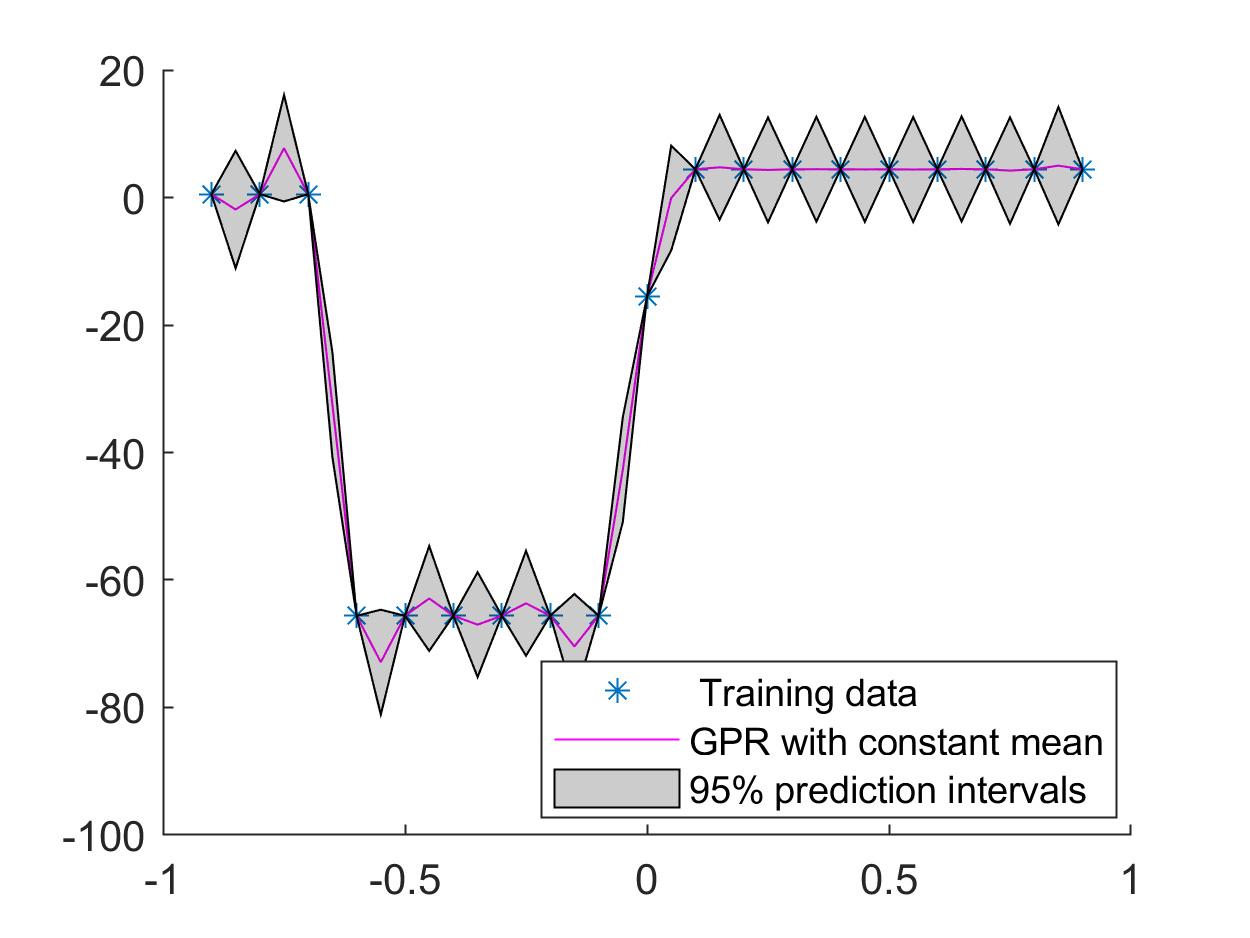}
    \caption{Mat\'{e}rn 5/2 Kernel}
  \end{subfigure}
  \begin{subfigure}{0.35\textwidth}
   \centering
   \includegraphics[height=4cm,width=5cm]{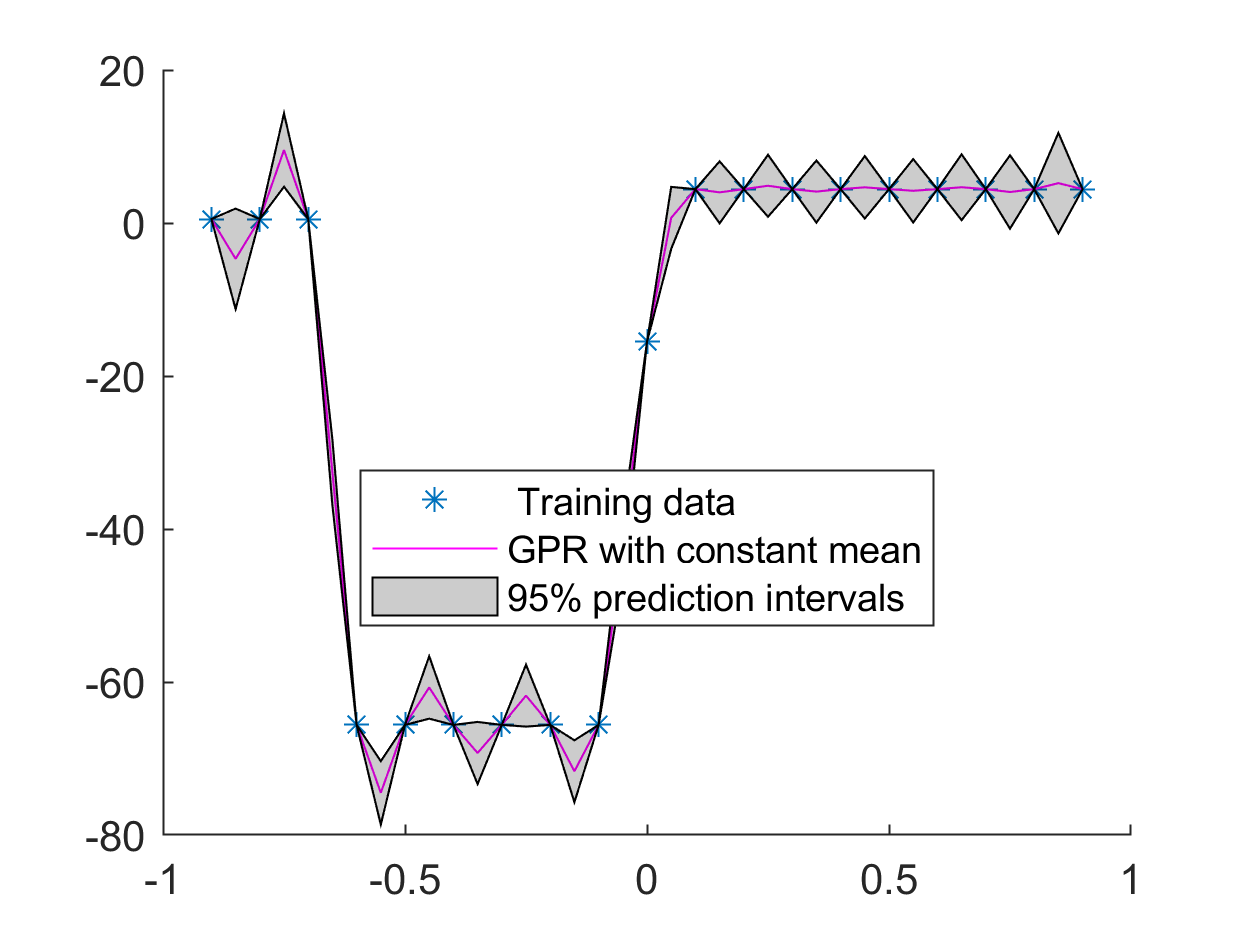}
    \caption{SE Kernel}
  \end{subfigure}
  \caption{GPR corresponding to 1st coefficient of the reduced 3rd eigenvector of EVP~\eqref{eq:crs} using different kernels and the training points are $-0.9:0.1:0.9$.}
  \label{crs:ev3_coeff2}
\end{figure}

\begin{figure}
\centering
 \begin{subfigure}{0.35\textwidth}
   \centering
   \includegraphics[height=4cm,width=5cm]{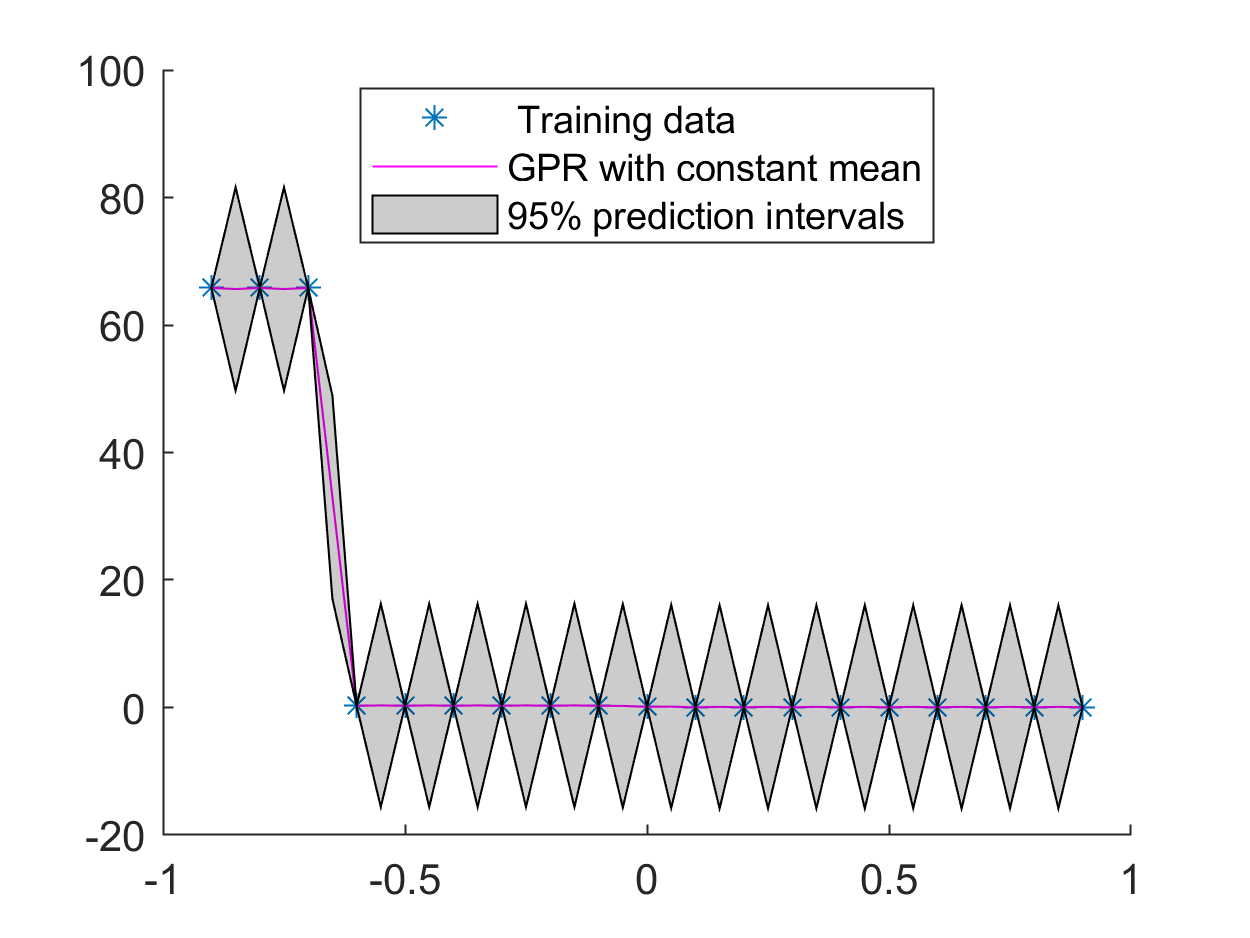}
    \caption{Exp Kernel}
  \end{subfigure}
 \begin{subfigure}{0.35\textwidth}
   \centering
   \includegraphics[height=4cm,width=5cm]{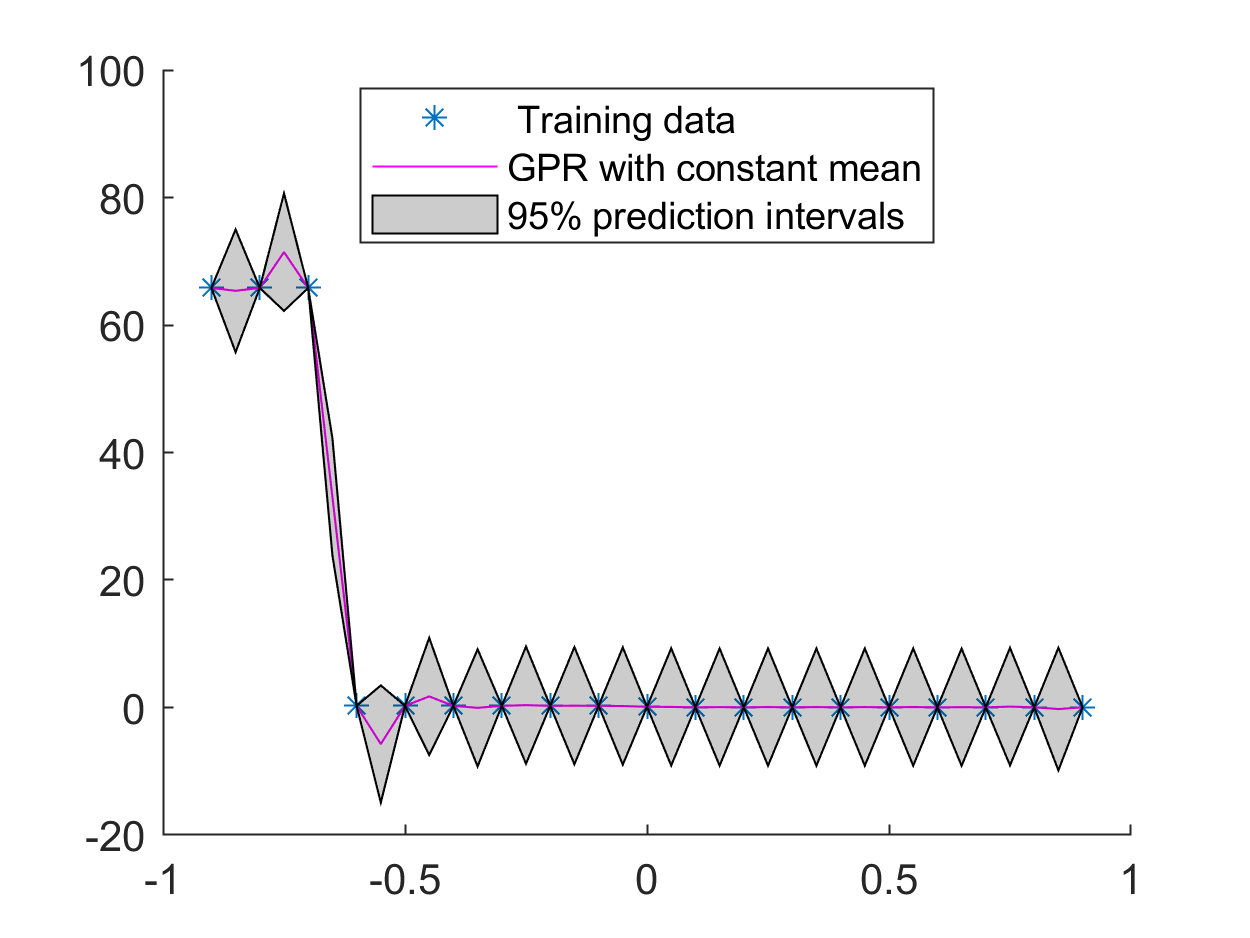}
    \caption{Mat\'{e}rn 3/2 Kernel}
  \end{subfigure}
  \begin{subfigure}{0.35\textwidth}
   \centering
   \includegraphics[height=4cm,width=5cm]{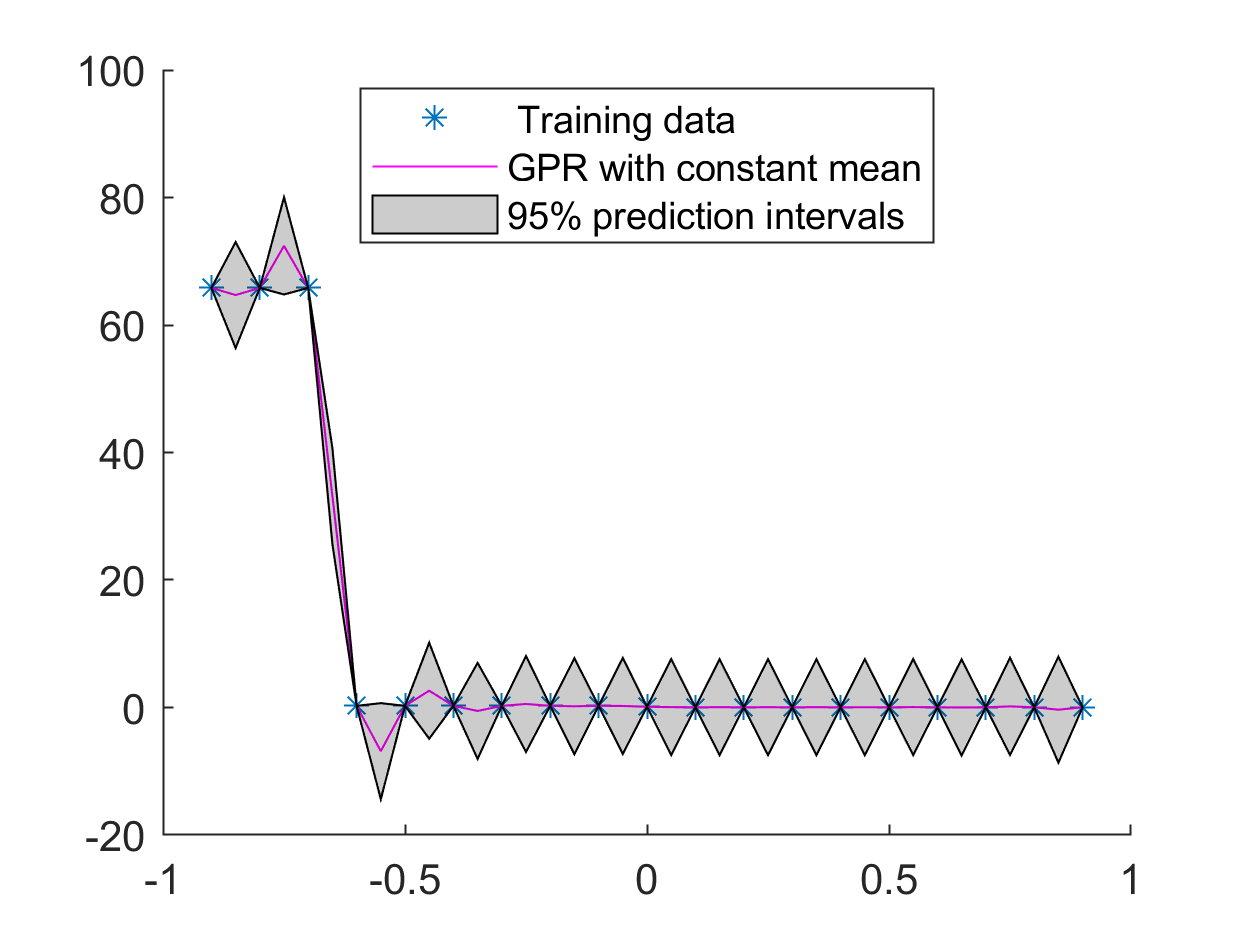}
    \caption{Mat\'{e}rn 5/2 Kernel}
  \end{subfigure}
  \begin{subfigure}{0.35\textwidth}
   \centering
   \includegraphics[height=4cm,width=5cm]{Result2/evct3_dm1_ma5_coef3.png}
    \caption{SE Kernel}
  \end{subfigure}
  \caption{GPR corresponding to 1st coefficient of the reduced 3rd eigenvector of EVP~\eqref{eq:crs} using different kernels and the training points are $-0.9:0.1:0.9$.}
  \label{crs:ev3_coeff3}
\end{figure}

\begin{figure}
\centering
 \begin{subfigure}{0.35\textwidth}
   \centering
   \includegraphics[height=4cm,width=5cm]{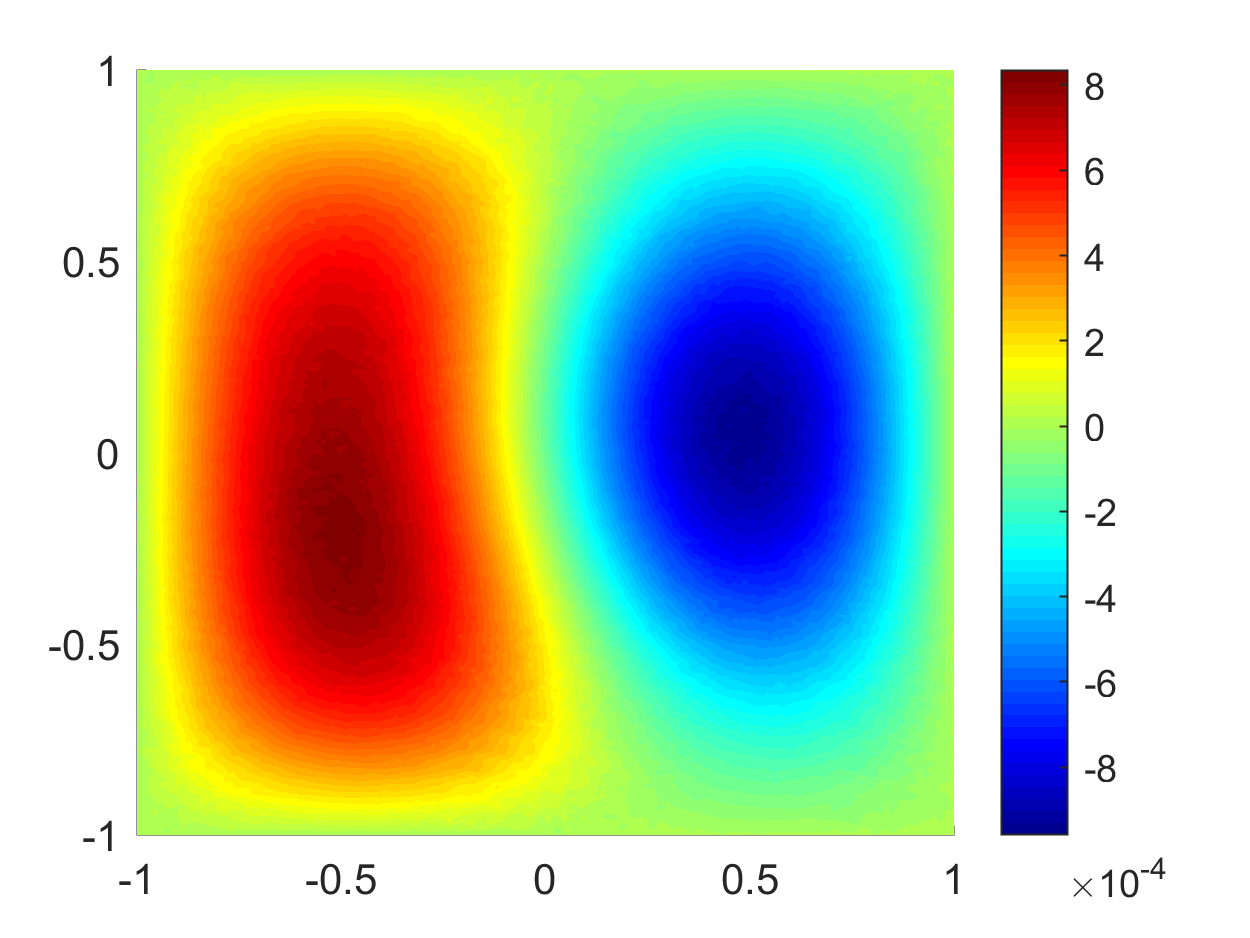}
    \caption{Exp Kernel}
  \end{subfigure}
  \begin{subfigure}{0.35\textwidth}
   \centering
   \includegraphics[height=4cm,width=5cm]{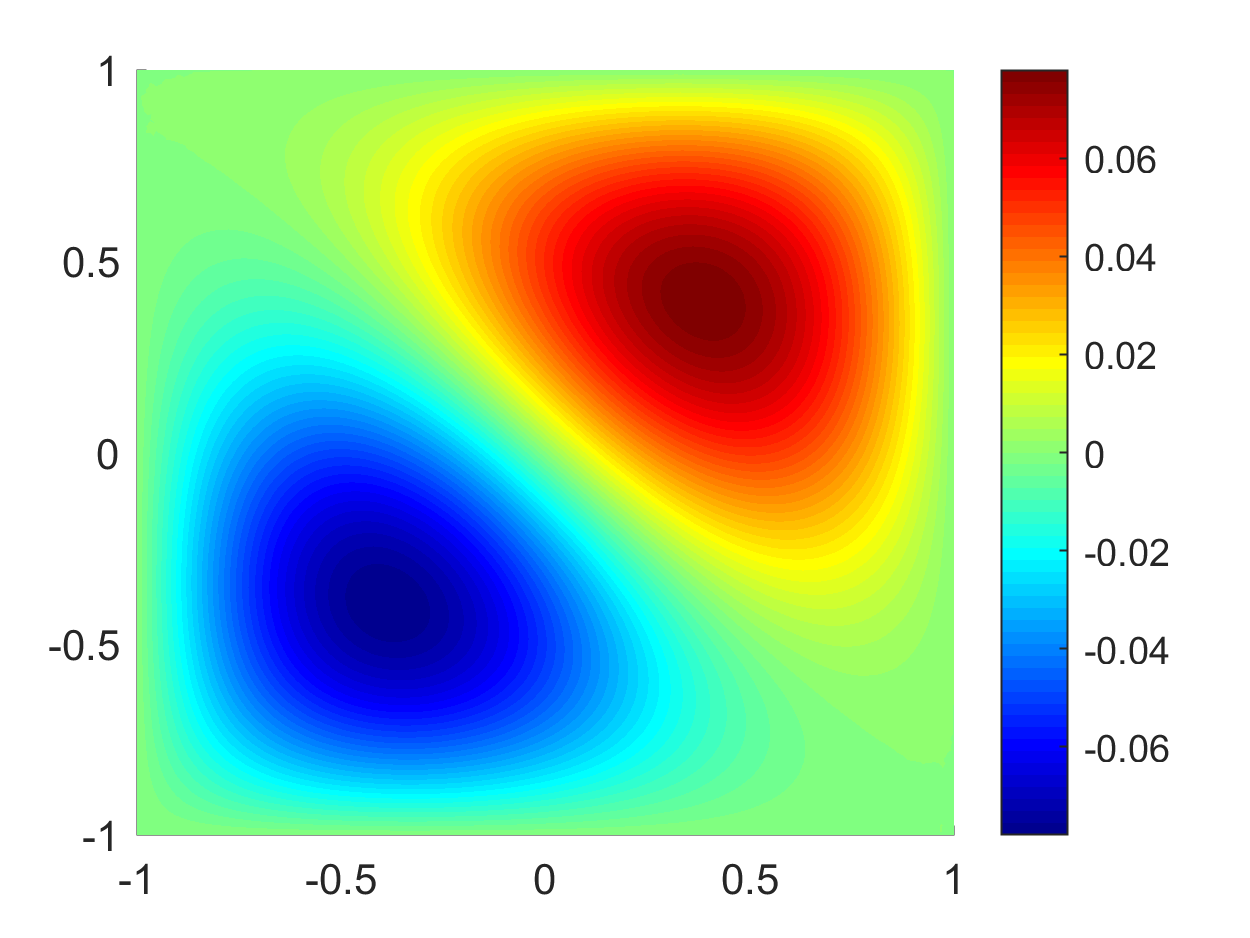}
    \caption{Mat\'{e}rn 3/2 Kernel}
  \end{subfigure}
  \begin{subfigure}{0.35\textwidth}
   \centering
   \includegraphics[height=4cm,width=5cm]{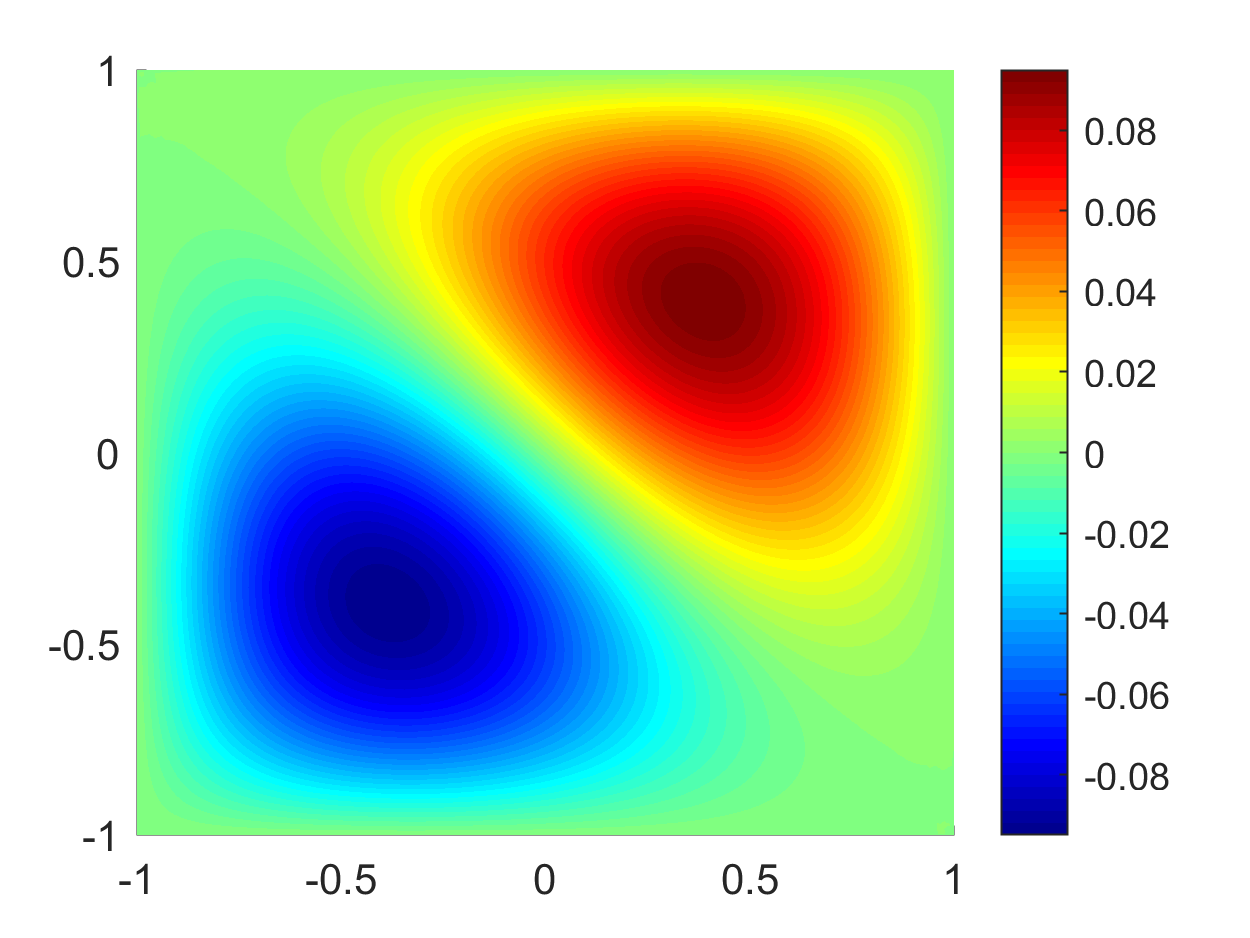}
    \caption{Mat\'{e}rn 5/2 Kernel}
  \end{subfigure}
  \begin{subfigure}{0.35\textwidth}
   \centering
   \includegraphics[height=4cm,width=5cm]{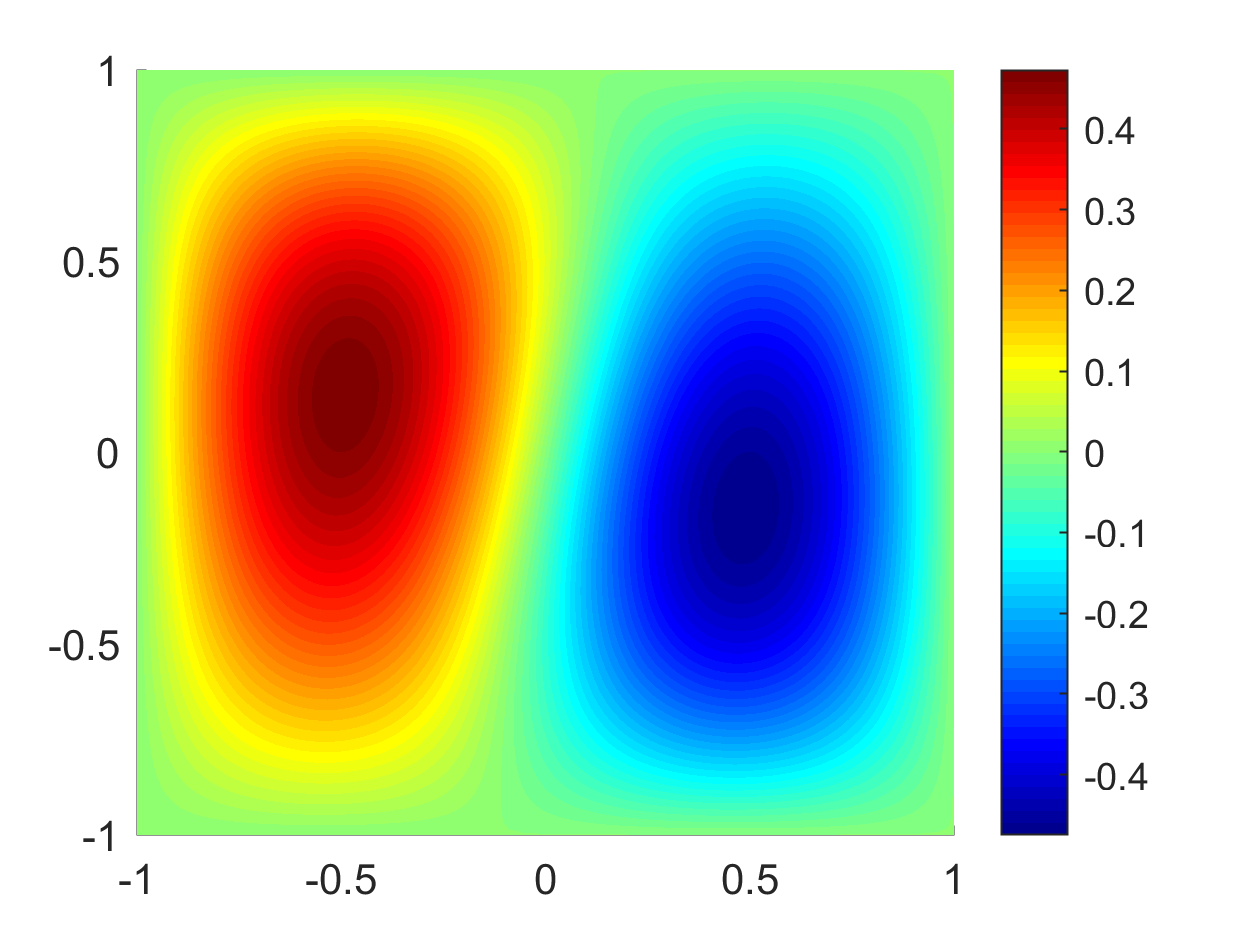}
    \caption{SE Kernel}
  \end{subfigure}
  \caption{Error between the third eigenvector using FEM and GPR of EVP~\eqref{eq:crs} with kernels exp, Mat\'{e}rn 3/2, Mat\'{e}rn 5/2 and SE kernel at $\mu=-0.05$ when sample points are $-0.9:0.04:0.9$.}
  \label{evct3:mu04}
\end{figure}

\begin{figure}
\centering
 % \begin{subfigure}{0.35\textwidth}
   % \centering
   \includegraphics[height=4cm,width=5cm]{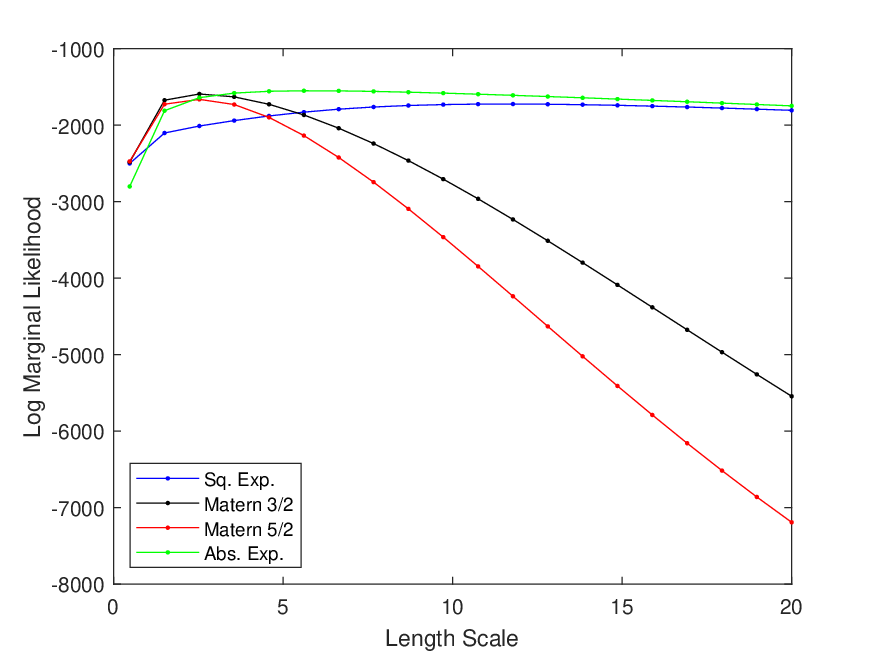}
    \caption{{A sample of the log marginal likelihood function. It is plotted as a function of the length scale $\ell$, and is evaluated at a fixed signal variance $\sigma_f$.}}
    \label{loglikelihoodEx}
  % \end{subfigure}
  % \begin{subfigure}{0.35\textwidth}
  %  \centering
  %  \includegraphics[height=4cm,width=5cm]{rev_result_crs/ev3_err_dmu04_m05_cov2.png}
  %   \caption{Mat\'{e}rn 3/2 Kernel}
  % \end{subfigure}
\end{figure}

{The GPRs corresponding to the third eigenvalue of the EVP~~\eqref{eq:crs} are presented in Fig~\ref{ev3:crs}. The RMSRE for four sets of training points are calculated and reported in Tab.~\ref{crs:ev3_new}. From Tab.~\ref{crs:ev3_new}, one can see that the GPR with absolute kernel performs better when the number of test points is increased. The GPRs corresponding to the three reduced coefficients of the third eigenvalues of the EVP~\eqref{eq:crs} are presented in Figs.~\ref{crs:ev3_coeff1}--~\ref{crs:ev3_coeff3}. In this case, the oscillation is also observed near the two jumps for all kernels except the exponential kernel.
In Fig.~\ref{evct3:mu04}, we report the error between the third eigenvalue computed by FEM and their GPR approximation at $\mu=-0.05$ with training set $-0.9:0.04:0.9$. Again, the point-wise error is of order $10^{-4}$ for the absolute exponential kernel. In Table~\ref{crs:ev3n}, we report the optimal value of the three hyperparameters of the GPR corresponding to the third eigenvalue and third eigenvector. Also, in this case, the optimal values are well within the range interval end separated from the endpoints.} {It is important to note that the log-likelihood function, in all of the cases discussed above, has a unique global extreme point. In Fig.~\ref{loglikelihoodEx}, we show an example of this solution.}

%%%%%%%%%%%%%%%%%%%%%%%%%%%%%%%%%%%%%%%%%%%%%%%%%%%%%%%%%%%%%%%%%%%%%%%%%%%%%%%
% %%%%% EXAMPLE 2.3- NON AFFINE PROBLEM 
%%%%%%%%%%%%%%%%%%%%%%%%%%%%%%%%%%%%%%%%%%%%%%%%%%%%%%%%%%%%%%%%%%%%%%%%%%%%%

\begin{figure}[H]
\centering
\includegraphics[height=4cm,width=6.5cm]{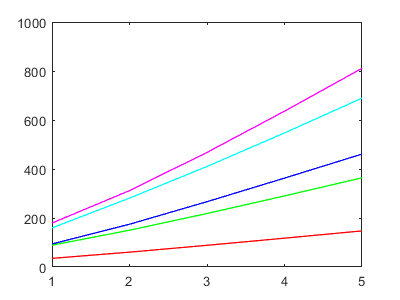}
\caption{First five eigenvalues of the EVP~\eqref{eq:nafn1}.}
\label{ev1:nafn1}
\end{figure}

\begin{example}
{We} consider the following eigenvalue problem
\begin{equation}
\label{eq:nafn1}
    \left\{
\aligned
&-\Delta u(\mu)=\lambda(\mu)e^{-\mu(x^2+y^2)}u(\mu)&&\text{in }\Omega=(0,1)^2\\
&u(\mu)=0&&\text{on }\partial\Omega,
\endaligned
\right.
\end{equation}
where the parameter $\mu$ belongs to the parameter domain $\mathcal{M}=[1,8]$. The parameters used for training are $\bm{\mu_{tr}} = 1:0.4:8$ and test points are $\bm{\mu_{test}} = 1:0.2:8$ for plotting the GPR. Note that this eigenvalue problem depends on the parameter in a non-affine way. The first five eigenvalues are shown in Fig.~\ref{ev1:nafn1}.
\end{example}

\begin{figure}
\centering
 \begin{subfigure}{0.35\textwidth}
   \centering
   \includegraphics[height=4cm,width=5cm]{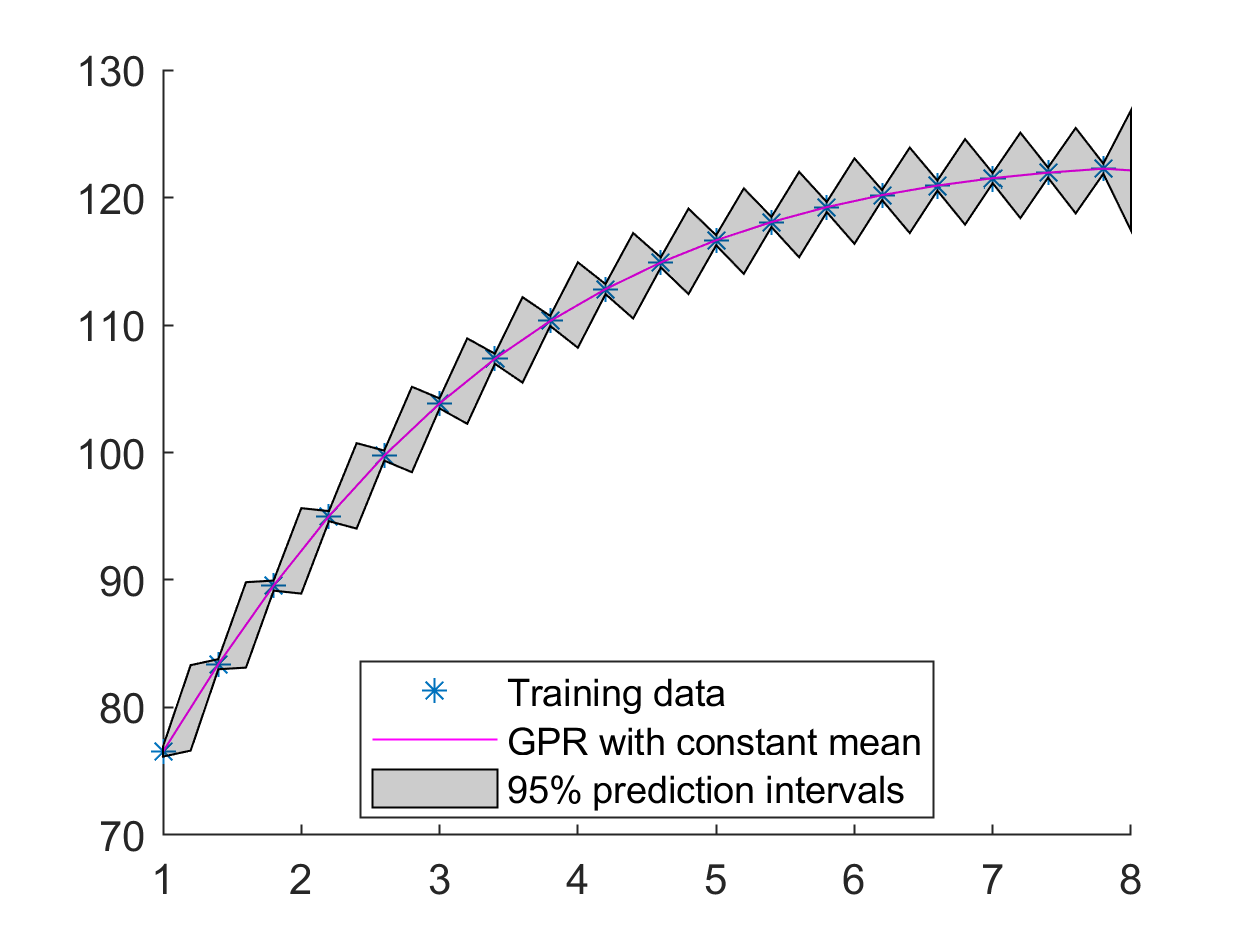}
    \caption{Exp Kernel}
  \end{subfigure}
 \begin{subfigure}{0.35\textwidth}
   \centering
\includegraphics[height=4cm,width=5cm]{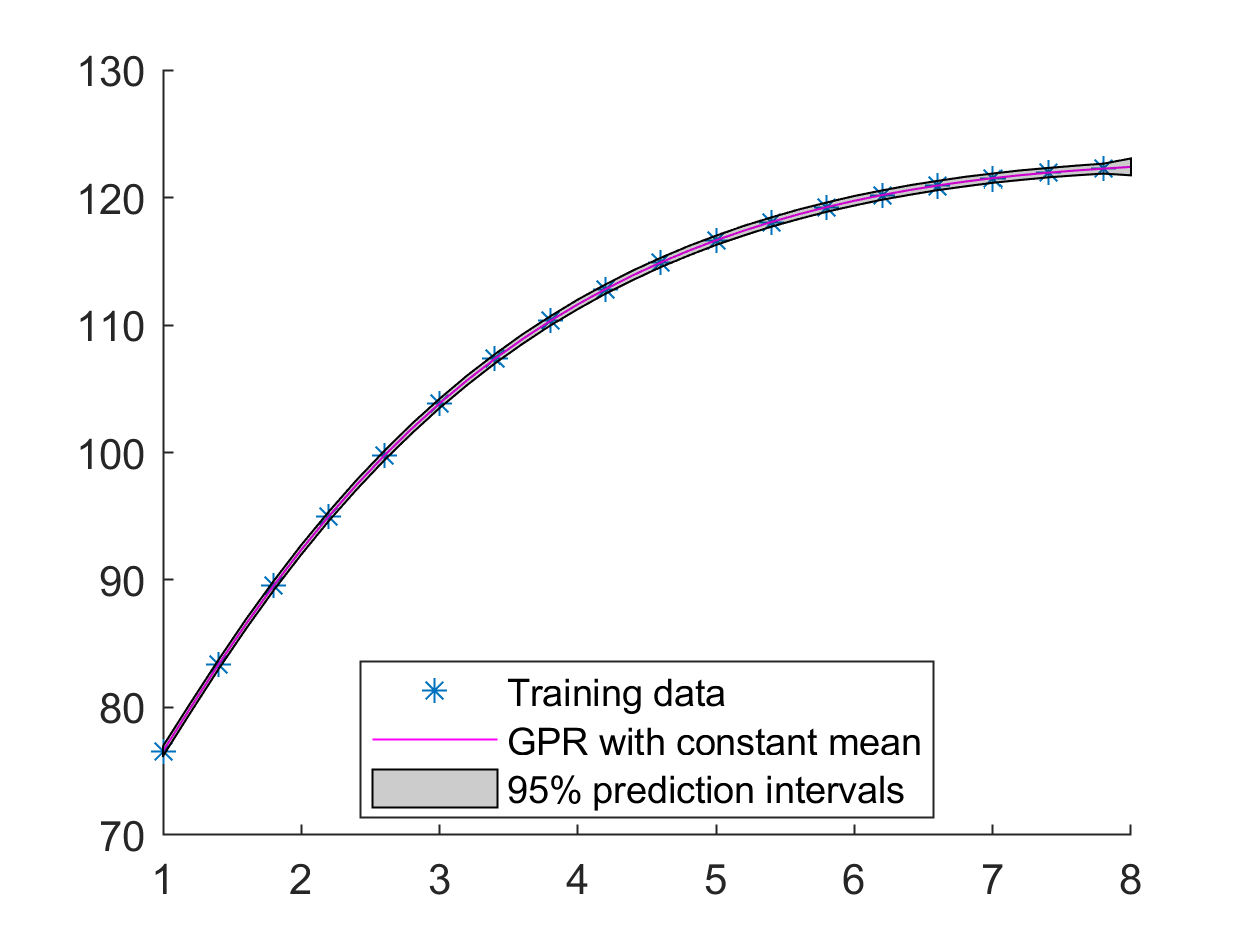}
    \caption{Mat\'{e}rn 3/2 Kernel}
  \end{subfigure}
  \begin{subfigure}{0.35\textwidth}
   \centering
\includegraphics[height=4cm,width=5cm]{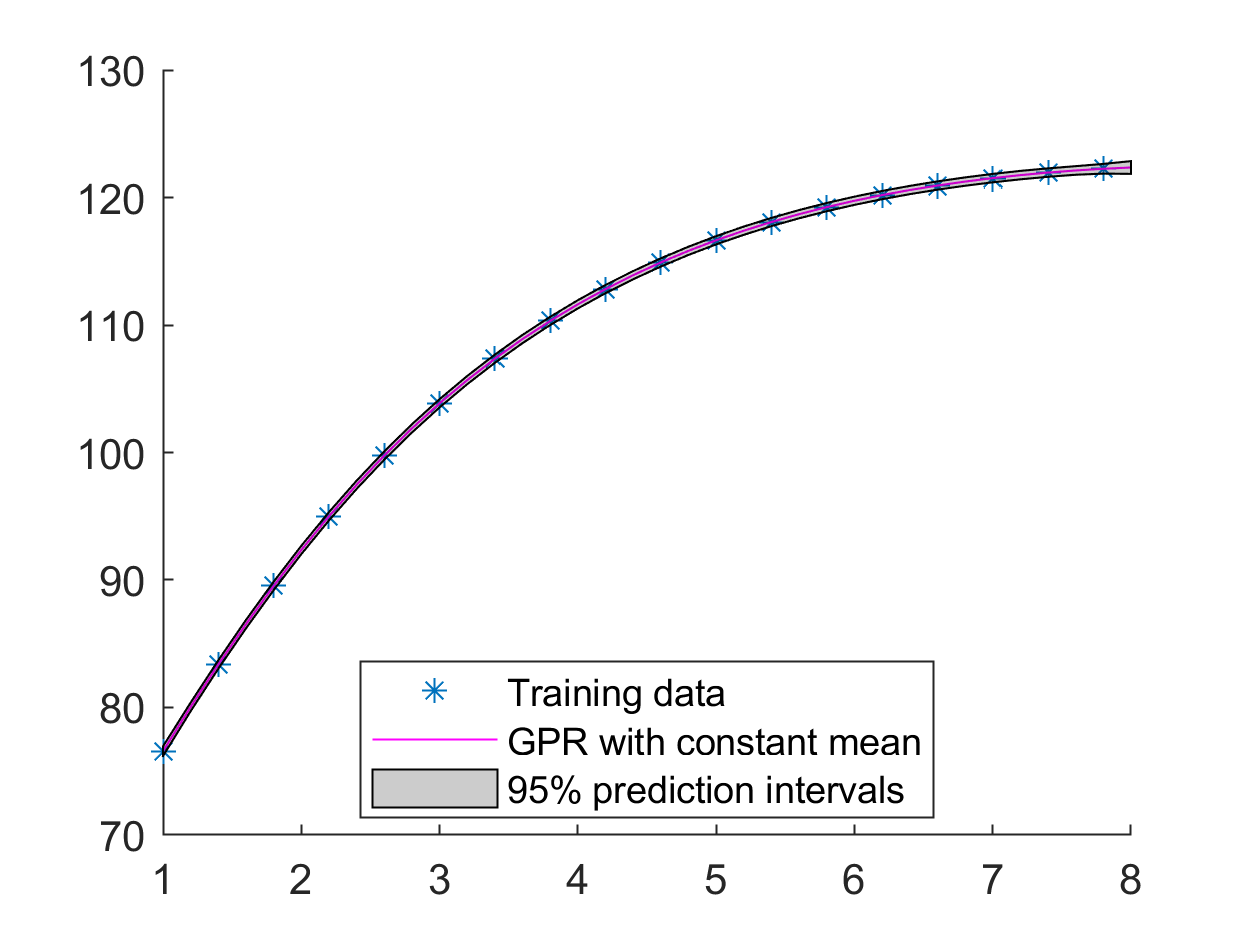}
    \caption{Mat\'{e}rn 5/2 Kernel}
  \end{subfigure}
  \begin{subfigure}{0.35\textwidth}
   \centering
\includegraphics[height=4cm,width=5cm]{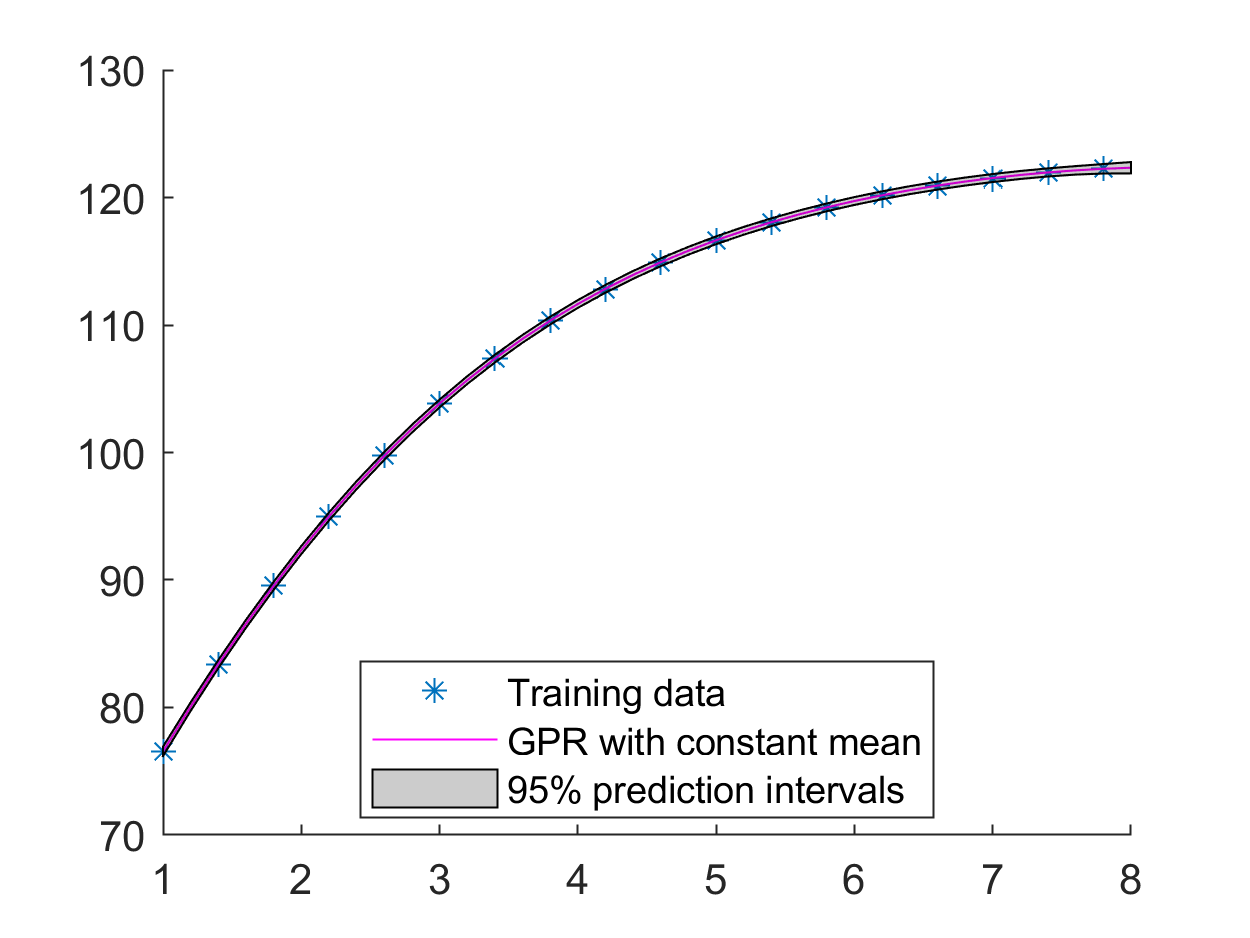}
    \caption{Squared Exponential Kernel}
  \end{subfigure}
  \caption{GPR corresponding to 1st coefficient of the Reduced 1st eigenvector of EVP~\eqref{eq:nafn1} using different kernels.}
  \label{nafn:ev1_coeff1}
\end{figure}

\begin{figure}
\centering
 \begin{subfigure}{0.35\textwidth}
   \centering
   \includegraphics[height=4cm,width=5cm]{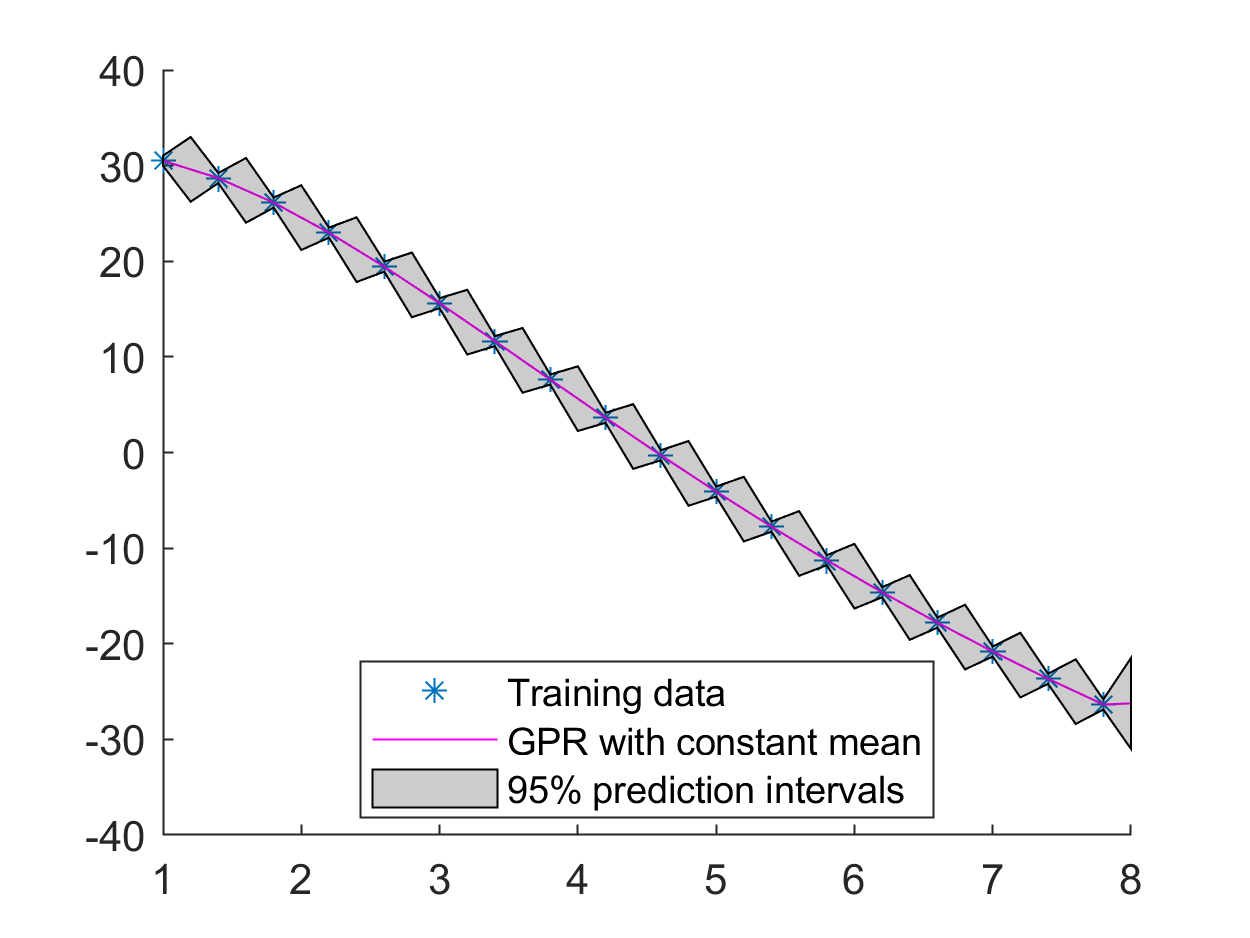}
    \caption{Exp Kernel}
  \end{subfigure}
 \begin{subfigure}{0.35\textwidth}
   \centering
\includegraphics[height=4cm,width=5cm]{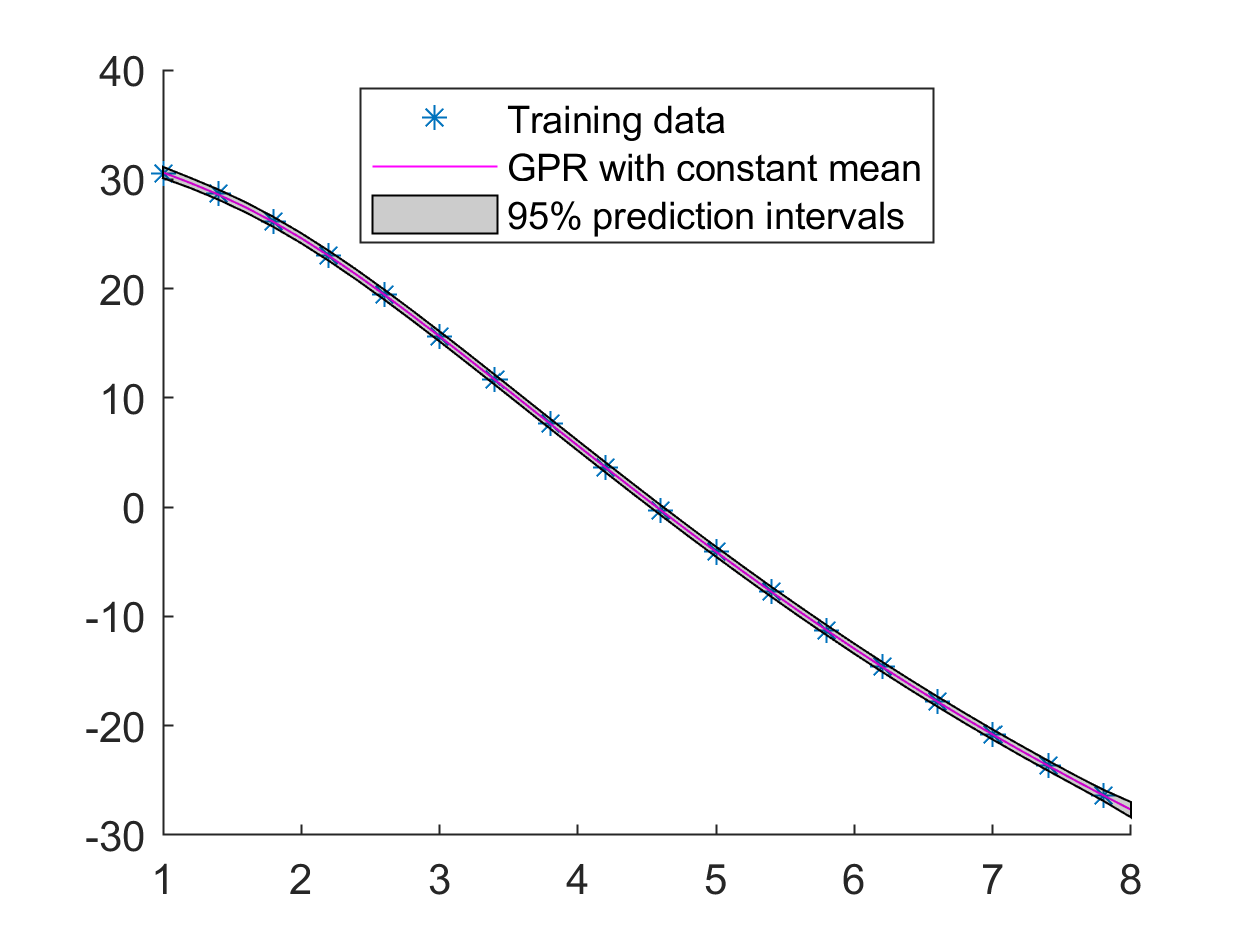}
    \caption{Mat\'{e}rn 3/2 Kernel}
  \end{subfigure}
  \begin{subfigure}{0.35\textwidth}
   \centering
\includegraphics[height=4cm,width=5cm]{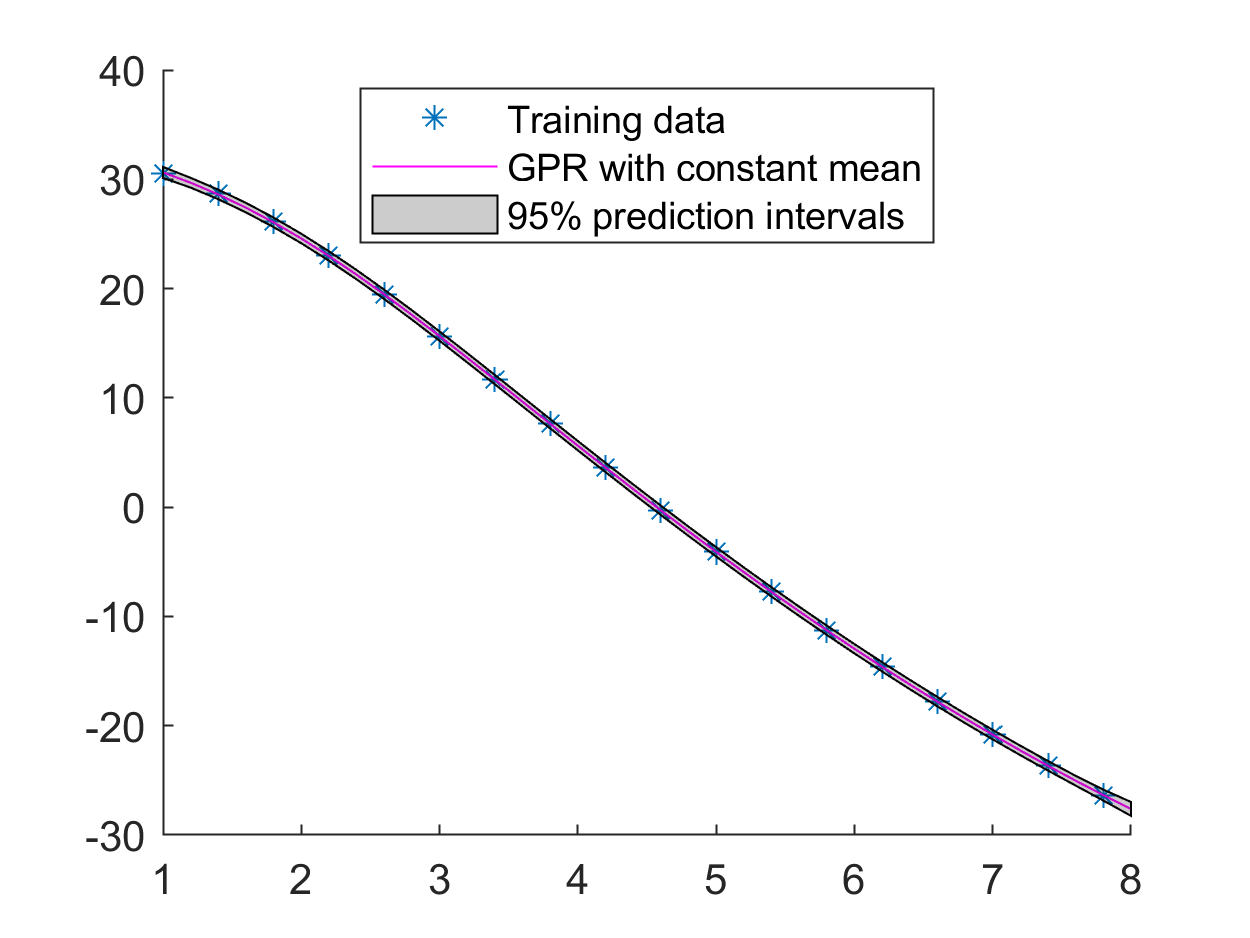}
    \caption{Mat\'{e}rn 5/2 Kernel}
  \end{subfigure}
  \begin{subfigure}{0.35\textwidth}
   \centering
\includegraphics[height=4cm,width=5cm]{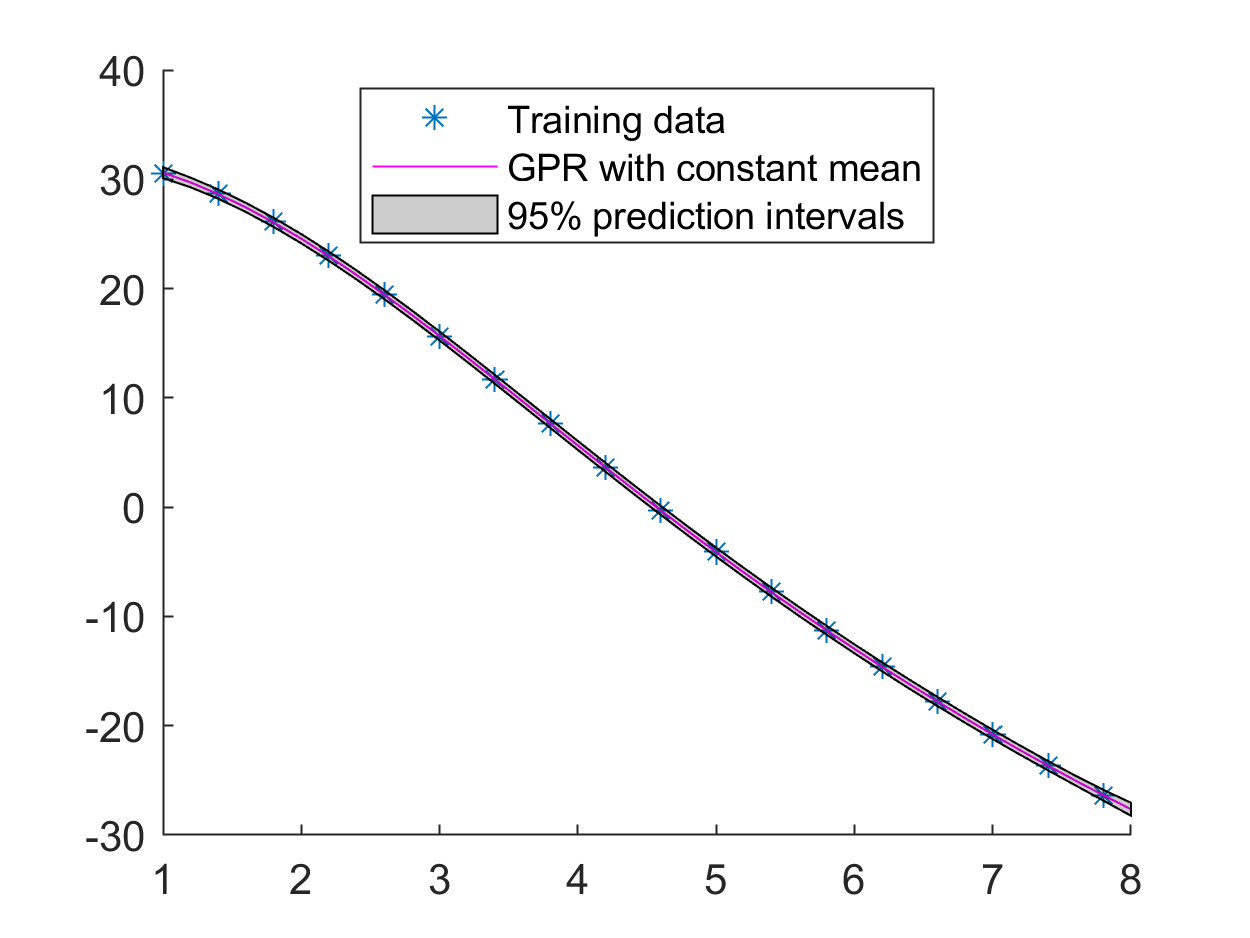}
    \caption{Squared Exponential Kernel}
  \end{subfigure}
  \caption{GPR corresponding to second coefficient of the Reduced 1st eigenvector of EVP~\eqref{eq:nafn1} using different kernels.}
  \label{nafn:ev1_coeff2}
\end{figure}

\begin{figure}
\centering
  \begin{subfigure}{0.35\textwidth}
   \includegraphics[height=4cm,width=5cm]{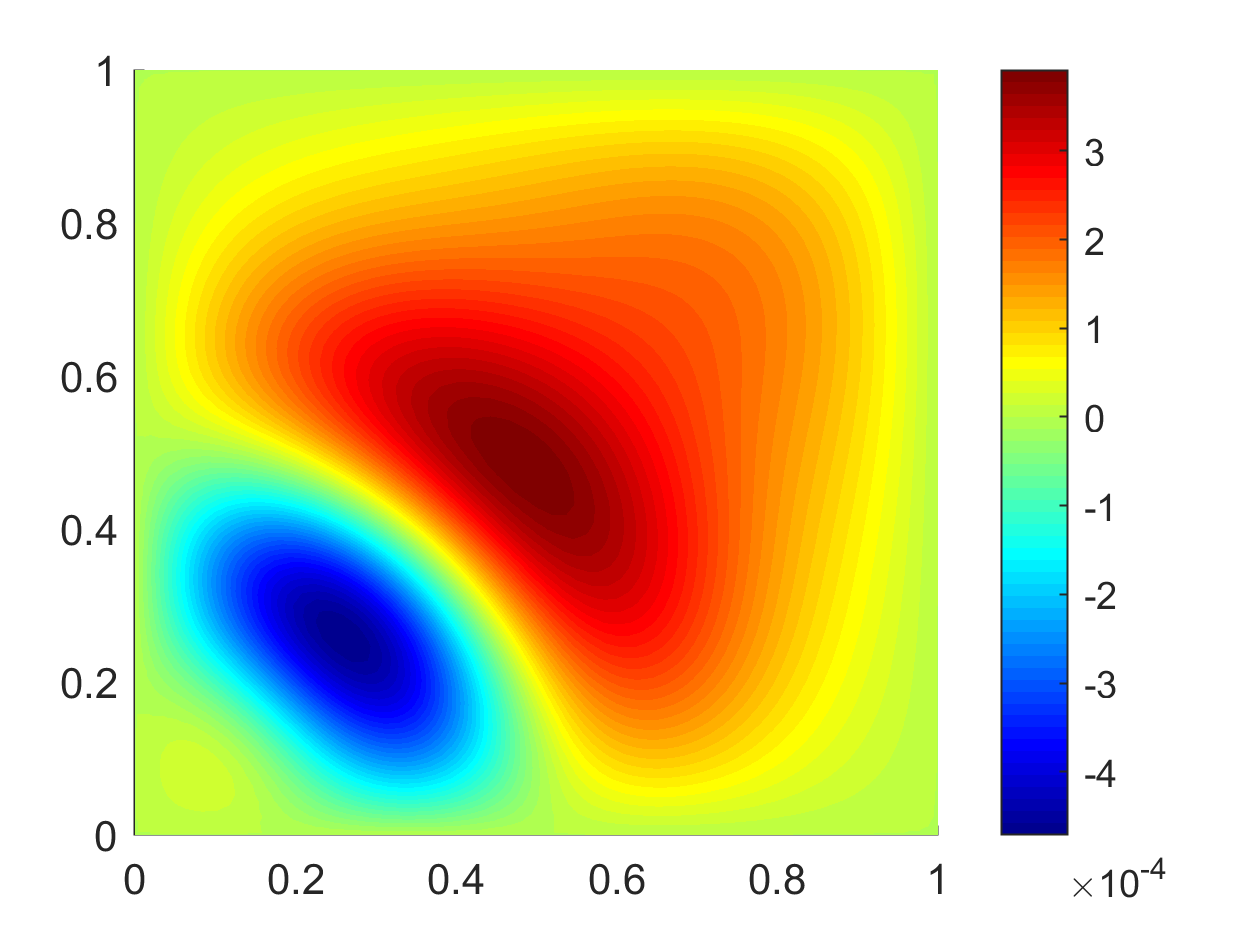}
    \caption{ Exp Kernel}
  \end{subfigure}
  \begin{subfigure}{0.35\textwidth}
   \includegraphics[height=4cm,width=5cm]{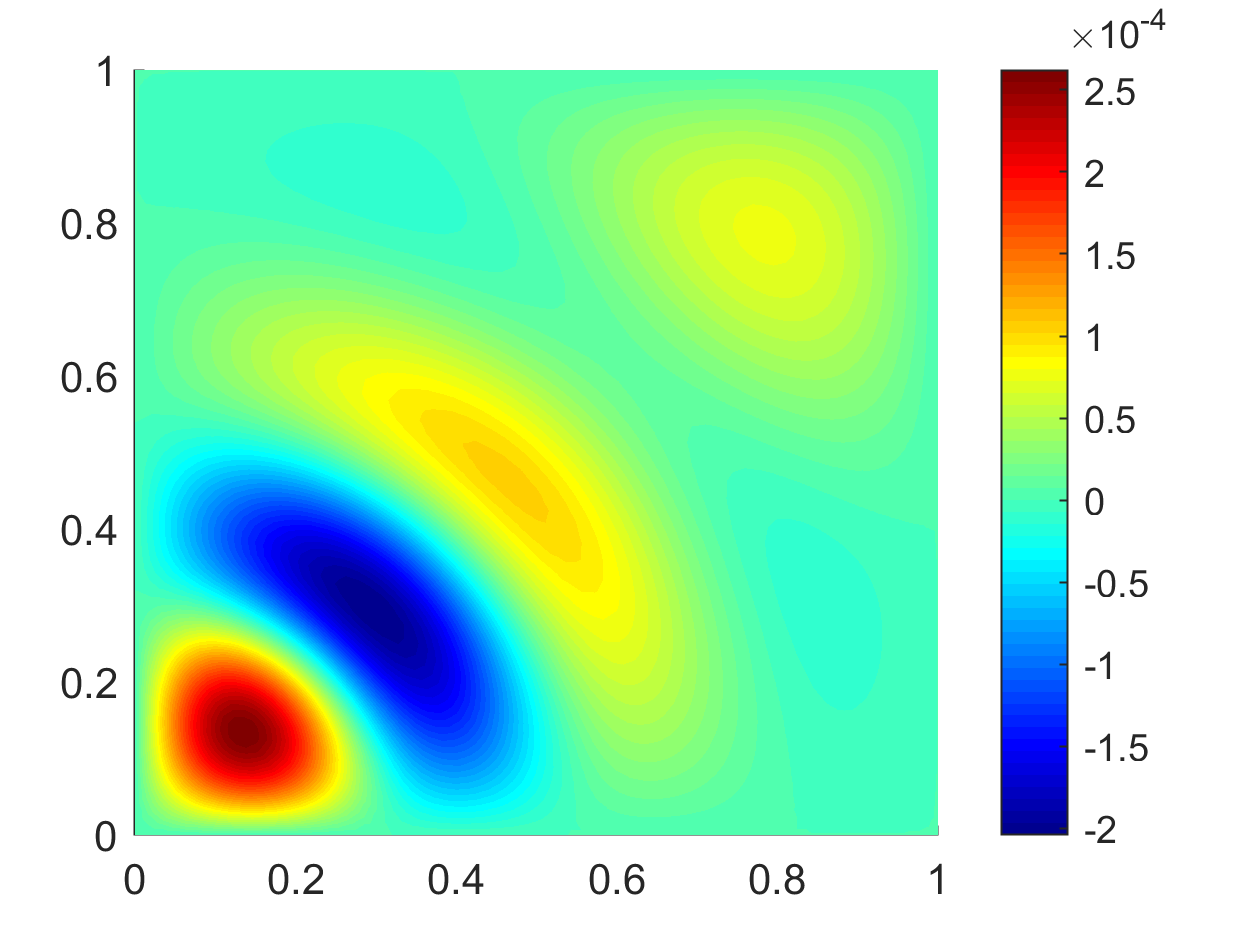}
    \caption{Mat\'{e}rn 3/2 Kernel}
  \end{subfigure}
  \begin{subfigure}{0.35\textwidth}
   \includegraphics[height=4cm,width=5cm]{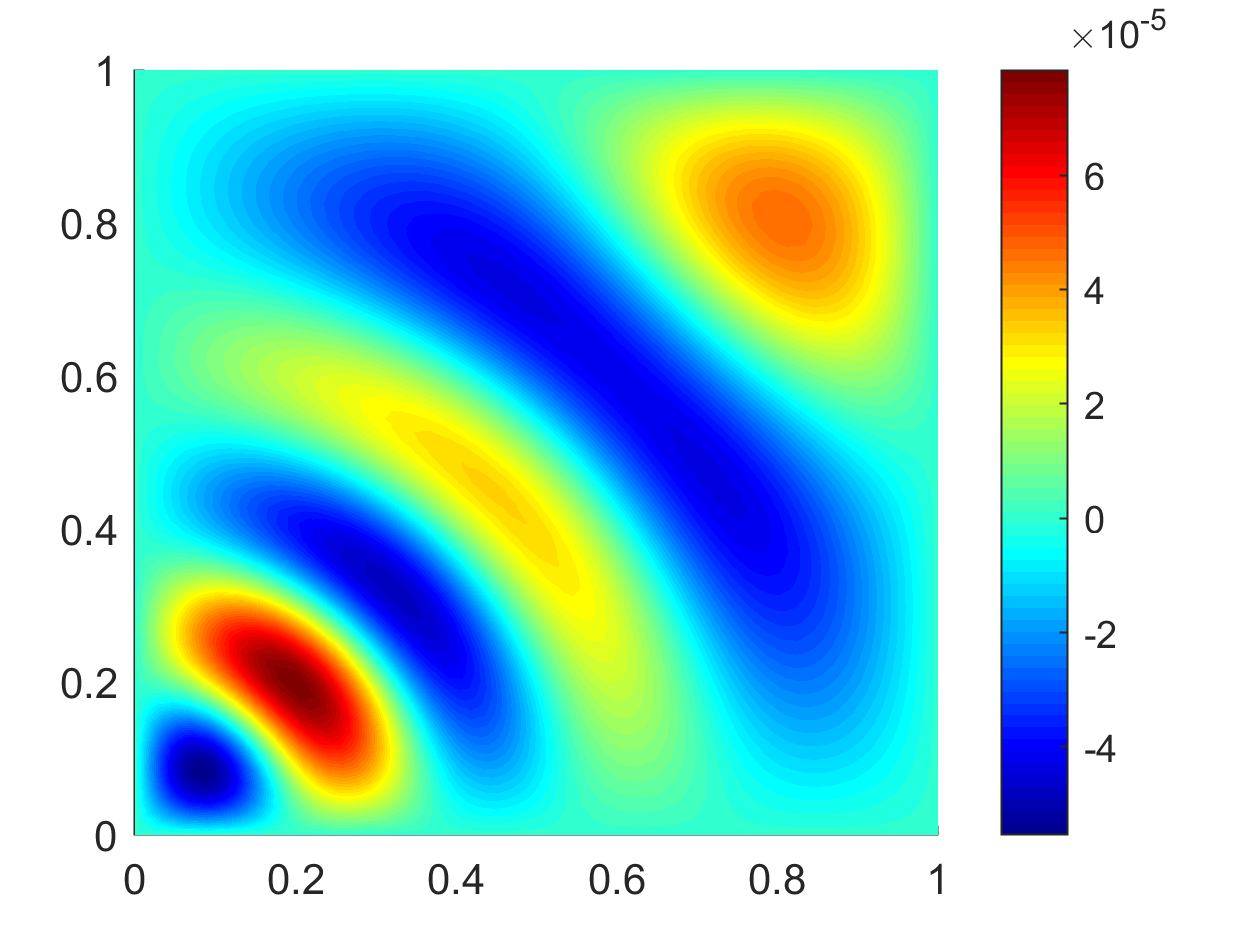}
    \caption{ Err Matrn 5/2 Kernel}
  \end{subfigure}
  \begin{subfigure}{0.35\textwidth}
   \includegraphics[height=4cm,width=5cm]{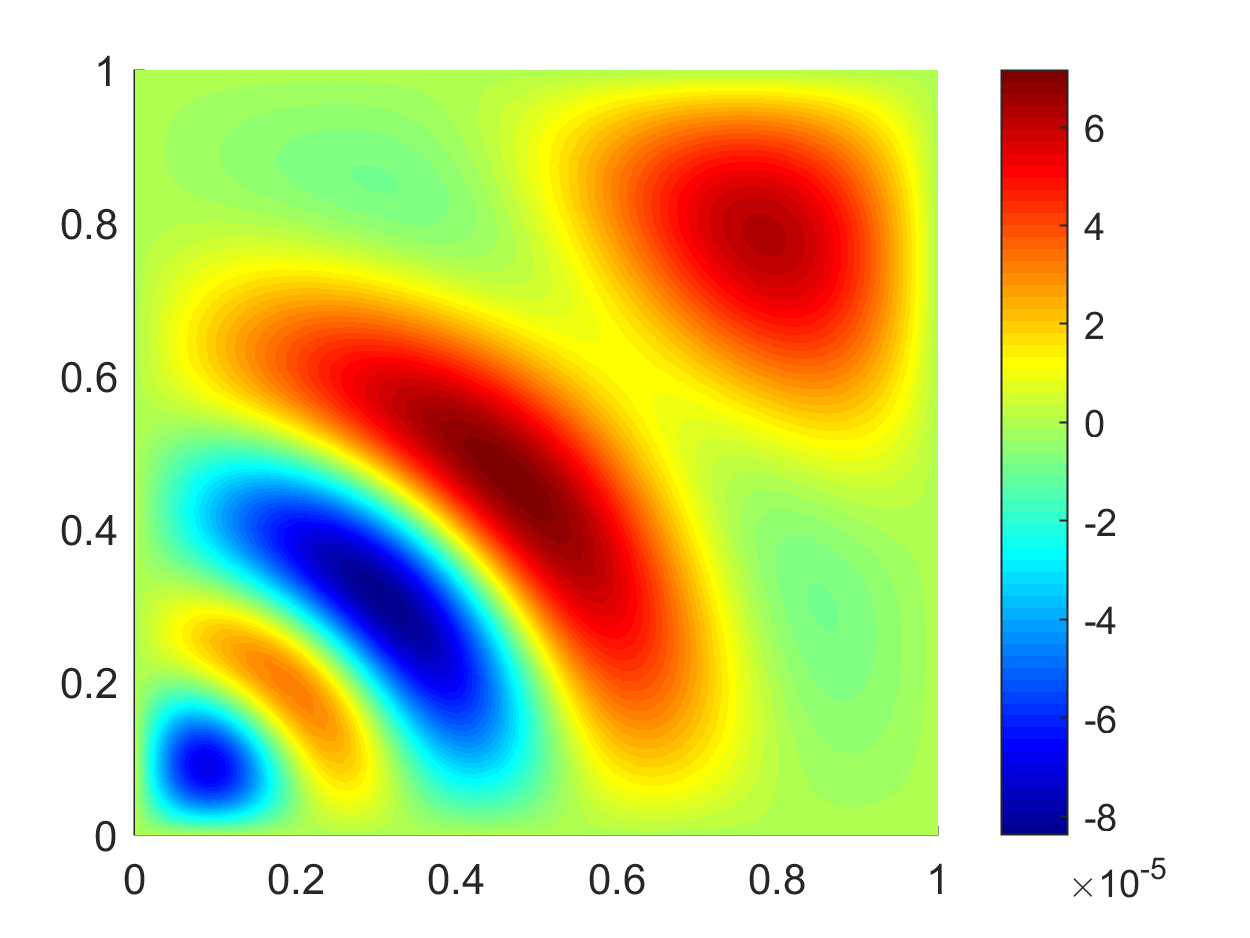}
    \caption{SE Kernel}
  \end{subfigure}
  \caption{Error between the FEM and GPR-based first eigenvectors of the EVP~\eqref{eq:nafn1} using different kernels at $\mu=7.6$.}
  \label{nafn1:evct1}
\end{figure}

{In Tab.~\ref{nafn1:ev}, we report the values of RMSRE for the first and second eigenvalue of the EVP~\eqref{eq:nafn1} using different sets of training points. To calculate the RMSRE, we have generated 40 random points from the parameter space $[1,8]$. From Tab.~\ref{nafn1:ev}, one can observe that the GPR with squared exponential performs better when we used a fewer number of training points. However, as we increased the number of points, the GPR with Mat\'{e}rn 5/2 kernel performs better. Note that in this case the curve corresponding to the eigenvalues is smooth. Consequently, it turns out that the GPRs with smooth kernels are better.}
{We now consider the eigenvectors associated with the first eigenvalue curve; as we already remarked, in this case there are no intersections so the parametric dependence is smooth. In this example our reduced space has dimension 5, so that we have five reduced coefficients corresponding to the first eigenvector. In Figs.~\ref{nafn:ev1_coeff1}--~\ref{nafn:ev1_coeff2}, we plot the GPR predictions for the first two reduced coefficients corresponding to the first eigenvector. The RMSRE for all the five coefficients corresponding to the first eigenvector is reported in Tab.~\ref{nafn:ev1_pb2_rrse} with different sets of training points. From the table, it can be observed that the GPR with the squared exponential kernel gives a smaller error when the number of training points is less but as we increase the number of training points Mat\'{e}rn 3/2  and Mat\'{e}rn 5/2 are better. This is because the graph of the coefficients is smooth once or twice but not infinitely differentiable.  Furthermore, we have shown as an example the error between the first eigenvector obtained by FEM and the corresponding eigenvectors using GPR with different kernels at $\mu=7.6$ in Fig.~\ref{nafn1:evct1}. }

\begin{table}
 \centering
  \begin{tabular}{|c|c|c|c|c|c|c|c|c|} 
  Eigenvalue& $\delta \mu$& Exponential  &Mat\'{e}rn 3/2 &Mat\'{e}rn 5/2  & Squared Exp.\\
    \hline
 $1^{st}$  
    & $0.4$ & $6.2\times 10^{-5}$ & $1.8\times 10^{-5}$ & $1.9\times 10^{-6}$ & $2.2\times 10^{-7}$\\
    & $0.2$ & $1.3\times 10^{-5}$ & $1.2\times 10^{-6}$ & $7.6\times 10^{-8}$ & $1.8\times 10^{-7}$\\
    &  $0.1$& $3.3\times 10^{-6}$ & $8.3\times 10^{-6}$ & $2.6\times 10^{-8}$ & $1.5\times 10^{-7}$\\
  %  & $0.05$& $8.3\times 10^{-7}$ & $6.4\times 10^{-10}$ & $1.9\times 10^{-8}$ & $3.5\times 10^{-8}$\\
    \hline
     $2^{nd}$  
    & $0.4$ & $4.9\times 10^{-5}$ & $1.2\times 10^{-5}$ & $4.9\times 10^{-7}$ & $2.1\times 10^{-8}$\\
    & $0.2$ & $1.1\times 10^{-5}$ & $8.9\times 10^{-7}$ & $1.7\times 10^{-8}$ & $1.8\times 10^{-8}$\\
    &  $0.1$& $2.6\times 10^{-6}$ & $1.0\times 10^{-7}$ & $3.0\times 10^{-9}$ & $1.7\times 10^{-8}$\\
  & $0.05$& $6.6\times 10^{-7}$ & $4.6\times 10^{-8}$ & $1.5\times 10^{-9}$ & $1,6\times 10^{-8}$\\
 \end{tabular}
\caption{{The values of RMSRE for the first and second eigenvalue of the EVP~\eqref{eq:nafn1} using GPR model with different covariance functions when training samples are $1:\delta \mu:8$.}}
\label{nafn1:ev}
 \end{table}

%%%%%%%%%%%%%%%%%%%%%%%%%%%%%% Remove all of these

\begin{table} %% RRSE
 \centering
  \begin{tabular}{|c|c|c|c|c|c|c|c|c|} 
 Coefficient & $\delta \mu$& Exponential  &Mat\'{e}rn 3/2 &Mat\'{e}rn 5/2  & Squared Exp.\\
\hline
 1st
 & $0.4$ &$7.2\times 10^{-5}$ & $1.7\times 10^{-5}$ & $6.1\times 10^{-7}$ & $2.0\times 10^{-8}$\\
 & $0.2$ &  $1.7\times 10^{-5}$ & $1.0\times 10^{-6}$ & $3.9\times 10^{-8}$ & $2.0\times 10^{-8}$\\
 & $0.1$&  $3.9\times 10^{-6}$ & $7.3\times 10^{-8}$ & $2.5\times 10^{-8}$ & $1.7\times 10^{-8}$\\
 & $0.05$&    $9.7\times 10^{-7}$ & $4.4\times 10^{-10}$ & $4.4\times 10^{-9}$ & $1.6\times 10^{-8}$\\
 \hline
 2nd 
 & $0.4$ &$3.3\times 10^{-4}$ & $5.2\times 10^{-5}$ & $8.9\times 10^{-6}$ & $5.7\times 10^{-7}$\\
 & $0.2$ &  $6.0\times 10^{-5}$ & $3.7\times 10^{-6}$ & $2.9\times 10^{-7}$ & $2.8\times 10^{-7}$\\
 & $0.1$&  $1.7\times 10^{-5}$ & $2.6\times 10^{-7}$ & $1.2\times 10^{-8}$ & $2.5\times 10^{-7}$\\
 & $0.05$&  $4.6\times 10^{-6}$ & $2.8\times 10^{-9}$ & $1.1\times 10^{-8}$ & $2.2\times 10^{-7}$\\
  \hline
 3rd
 & $0.4$ &$2.3\times 10^{-3}$ & $5.9\times 10^{-5}$ & $3.4\times 10^{-5}$ & $2.6\times 10^{-6}$\\
 & $0.2$ &  $1.5\times 10^{-3}$ & $2.1\times 10^{-6}$ & $1.3\times 10^{-6}$ & $2.2\times 10^{-6}$\\
 & $0.1$&  $4.4\times 10^{-4}$ & $6.1\times 10^{-8}$ & $2.0\times 10^{-7}$ & $1.4\times 10^{-5}$\\
 & $0.05$&  $8.2\times 10^{-5}$ & $9.1\times 10^{-10}$ & $9.2\times 10^{-8}$ & $4.6\times 10^{-6}$\\
 \hline
4th
 & $0.4$ & 0.0254 & 0.0125 & 0.0011 & 0.0060\\
 & $0.2$   &$3.9\times 10^{-3}$ & $7.9\times 10^{-5}$ & $4.9\times 10^{-7}$ & $2.1\times 10^{-5}$\\
 & $0.1$  &$1.3\times 10^{-3}$ & $3.8\times 10^{-6}$ & $3.2\times 10^{-7}$ & $2.9\times 10^{-5}$\\
 & $0.05$& $1.4\times 10^{-3}$ & $8.8\times 10^{-8}$ & $1.4\times 10^{-6}$ & $5.9\times 10^{-5}$\\
\hline
5th
 & $0.4$ & 0.0526 &  0.0669 & 0.1445 &0.0015\\
 & $0.2$ &   0.0832 & $3.8\times 10^{-3}$ & $6.5\times 10^{-4}$ & $2.3\times 10^{-4}$\\
 & $0.1$&  $2.4\times 10^{-3}$ & $3.2\times 10^{-3}$ & $4.5\times 10^{-3}$ & $1.1\times 10^{-4}$\\
 & $0.05$&   $6.5\times 10^{-4}$ & $3.5\times 10^{-7}$ & $1.7\times 10^{-5}$ & $1.6\times 10^{-4}$\\
 \end{tabular}
\caption{{The values of RMSRE for the reduced coefficients corresponding to the first eigenvector of the EVP~\eqref{eq:nafn1} using GPR model with different covariance functions when training samples are $1:\delta \mu:8$.}}
\label{nafn:ev1_pb2_rrse}
 \end{table}

%%%%%%%%%%%%%%%%%%%%%%%%%%%%%%%%%%%%%

\begin{figure}
\centering
 \begin{subfigure}{0.35\textwidth}
   \includegraphics[height=4cm,width=5cm]{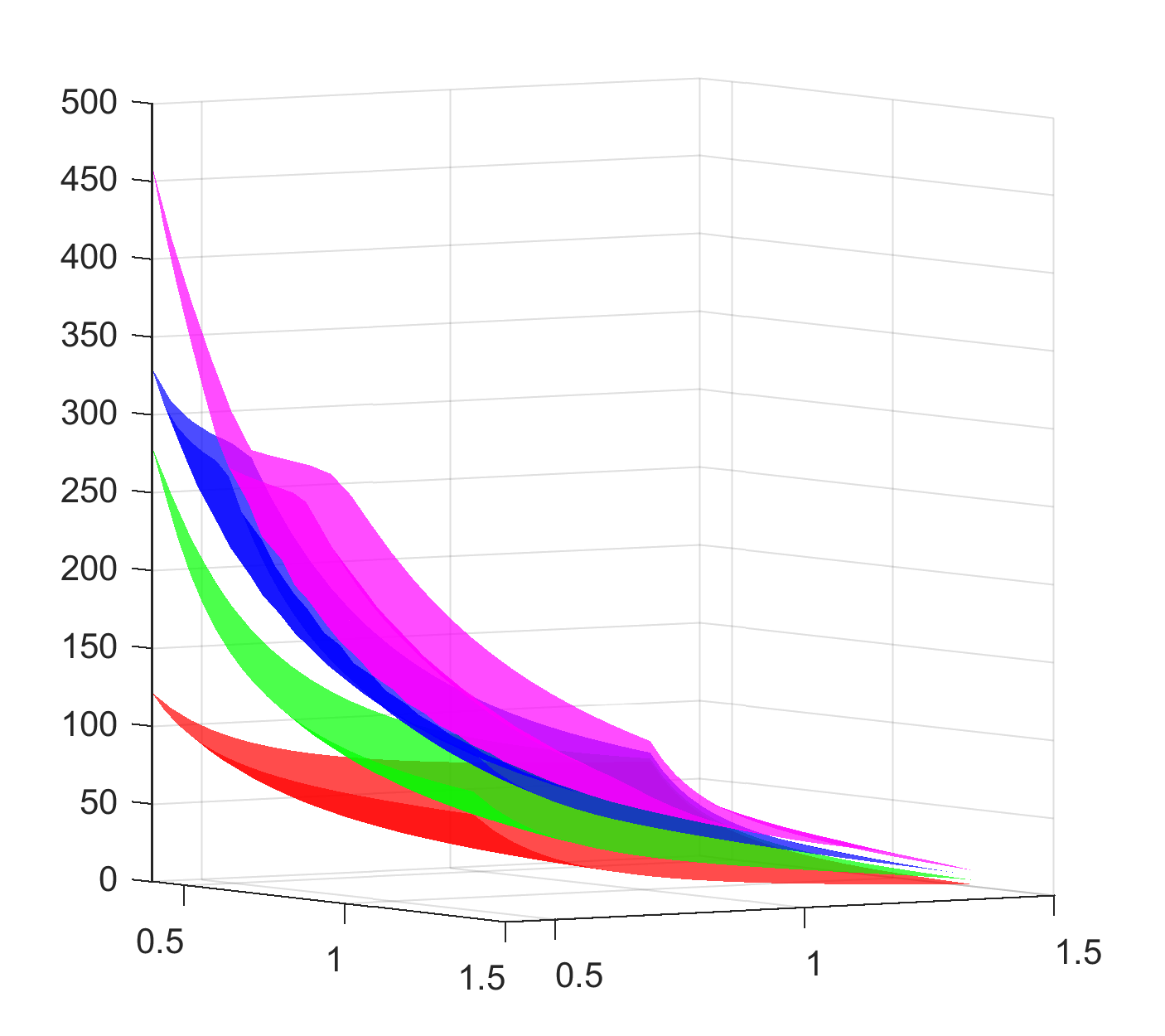}
    \caption{First four eigenvalues}
  \end{subfigure}
  \begin{subfigure}{0.35\textwidth}
   \includegraphics[height=4cm,width=5cm]{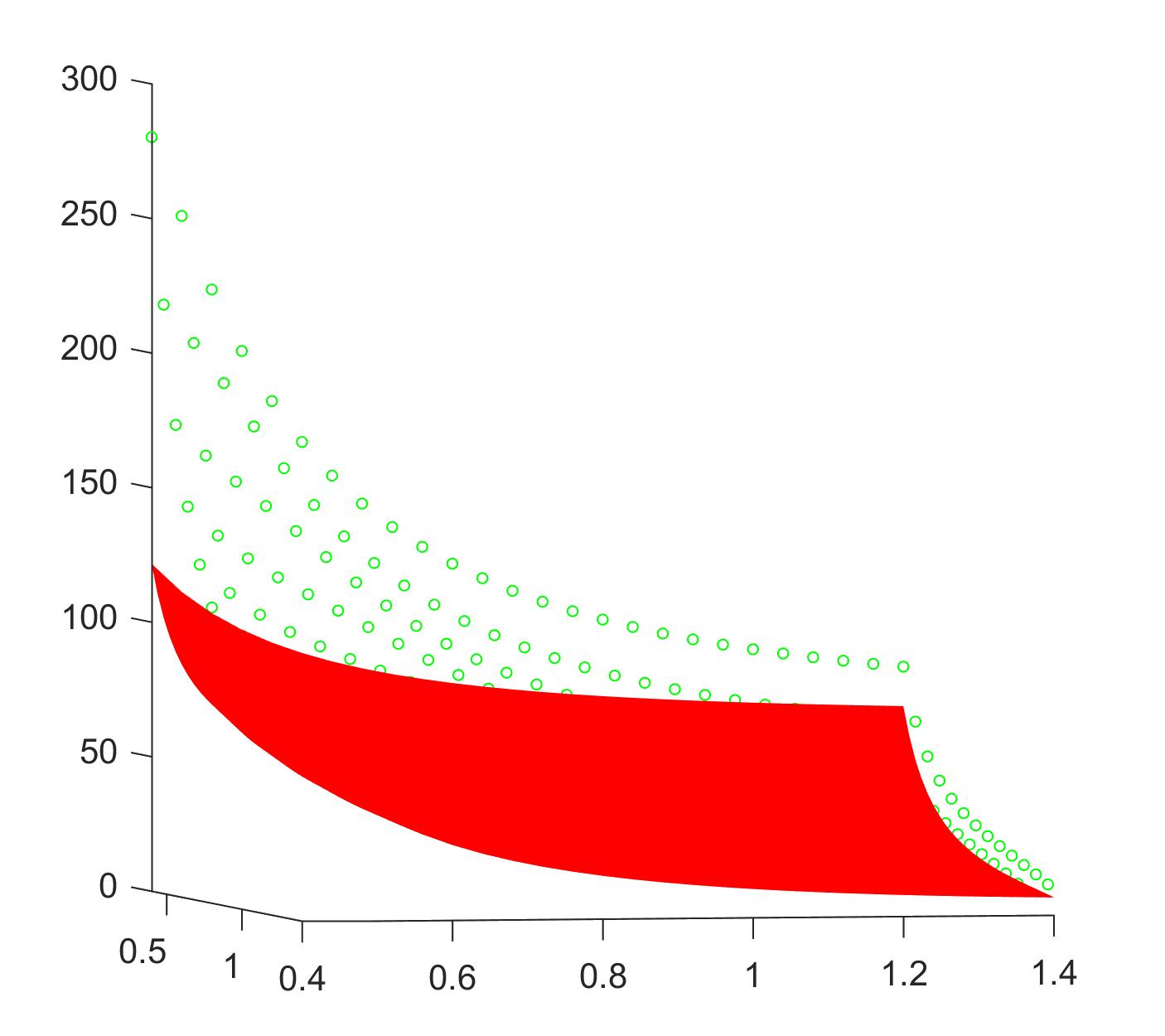}
    \caption{1st and 2nd eigenvalues}
  \end{subfigure}\\
  \begin{subfigure}{0.35\textwidth}
   \includegraphics[height=4cm,width=5cm]{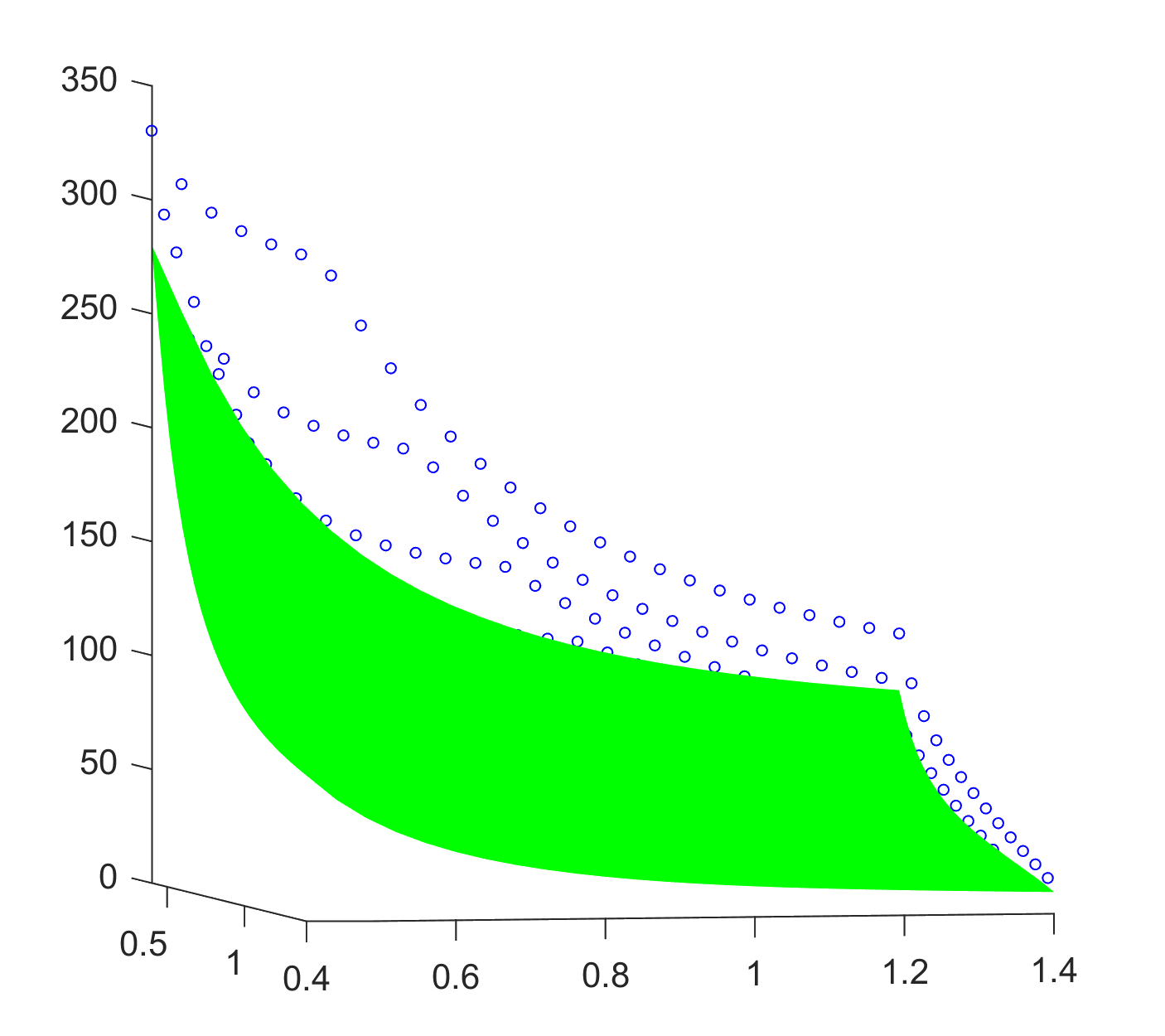}
    \caption{2nd and 3rd eigenvalues}
  \end{subfigure}
  \begin{subfigure}{0.35\textwidth}
   \includegraphics[height=4cm,width=5cm]{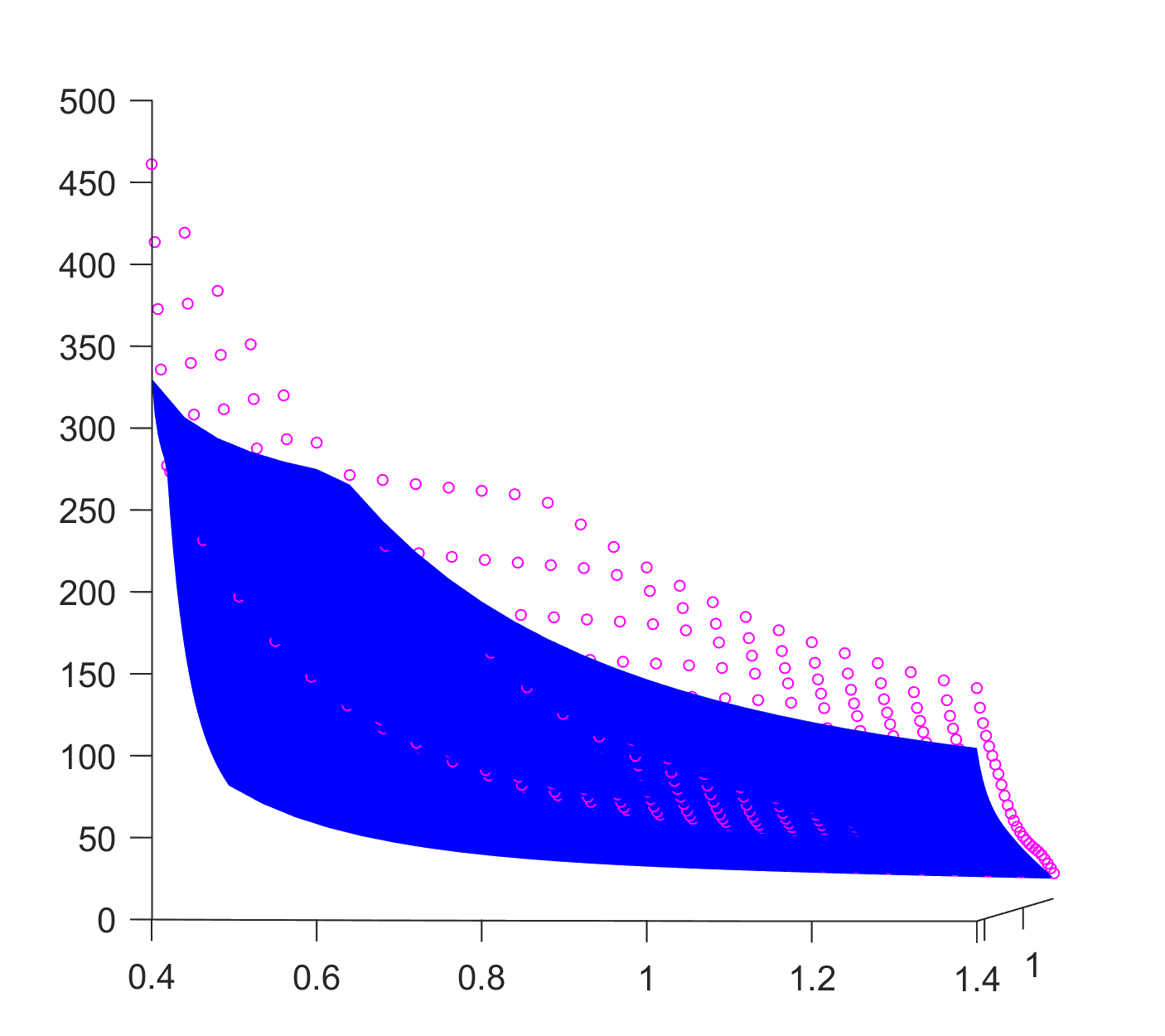}
    \caption{3rd and 4th eigenvalues}
  \end{subfigure}
  \caption{The eigenvalues of EVP~\eqref{eq:par2} as surfaces depending on the parameters. }
  \label{pltev:par2}
\end{figure}

\begin{figure}
\centering
   \begin{subfigure}{0.3\textwidth}
   \includegraphics[height=4cm,width=5cm]{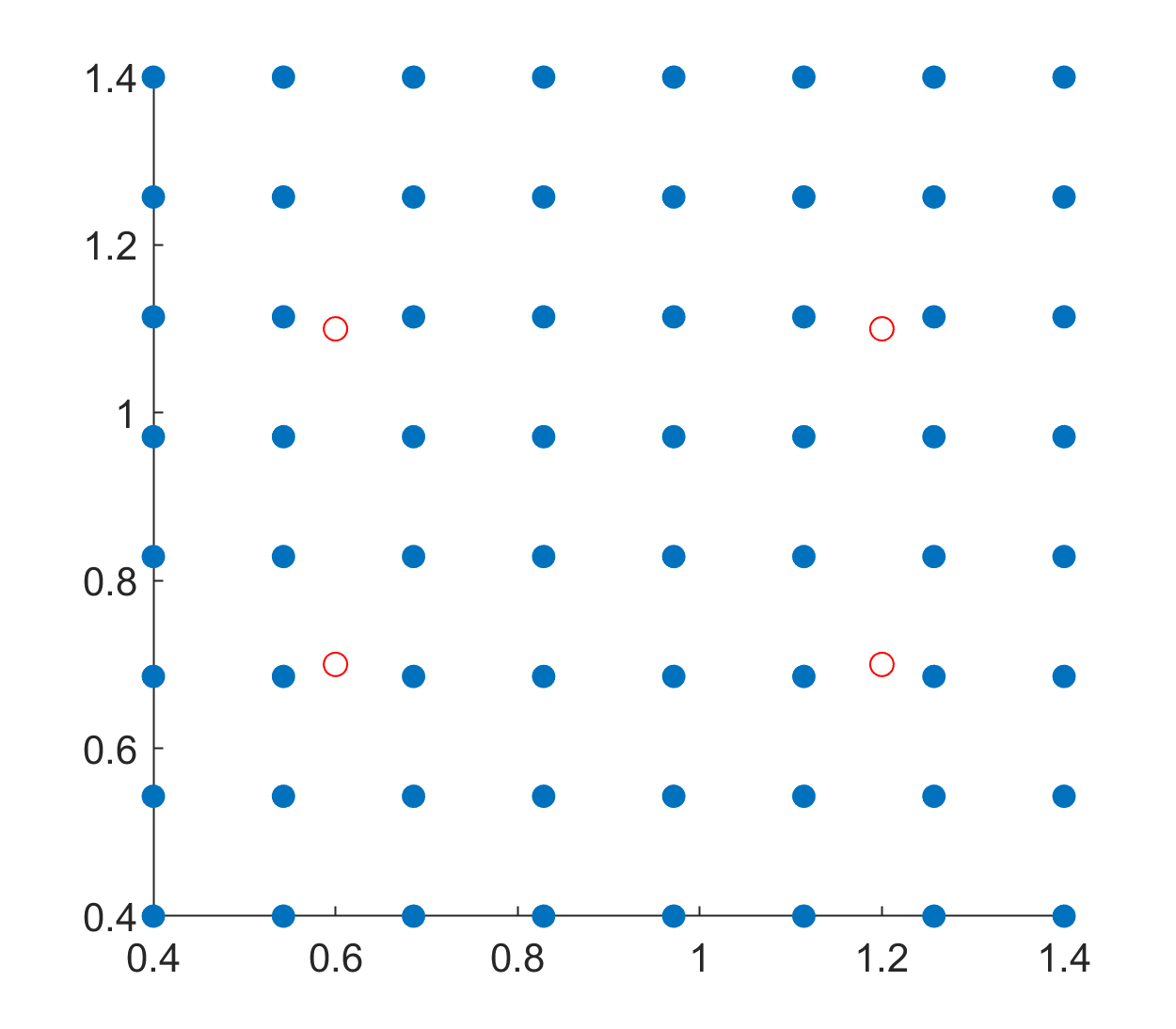}
    \caption{64 uniform Sample points}
    \label{uni3}
  \end{subfigure}
   \begin{subfigure}{0.3\textwidth}
   \includegraphics[height=4cm,width=5cm]{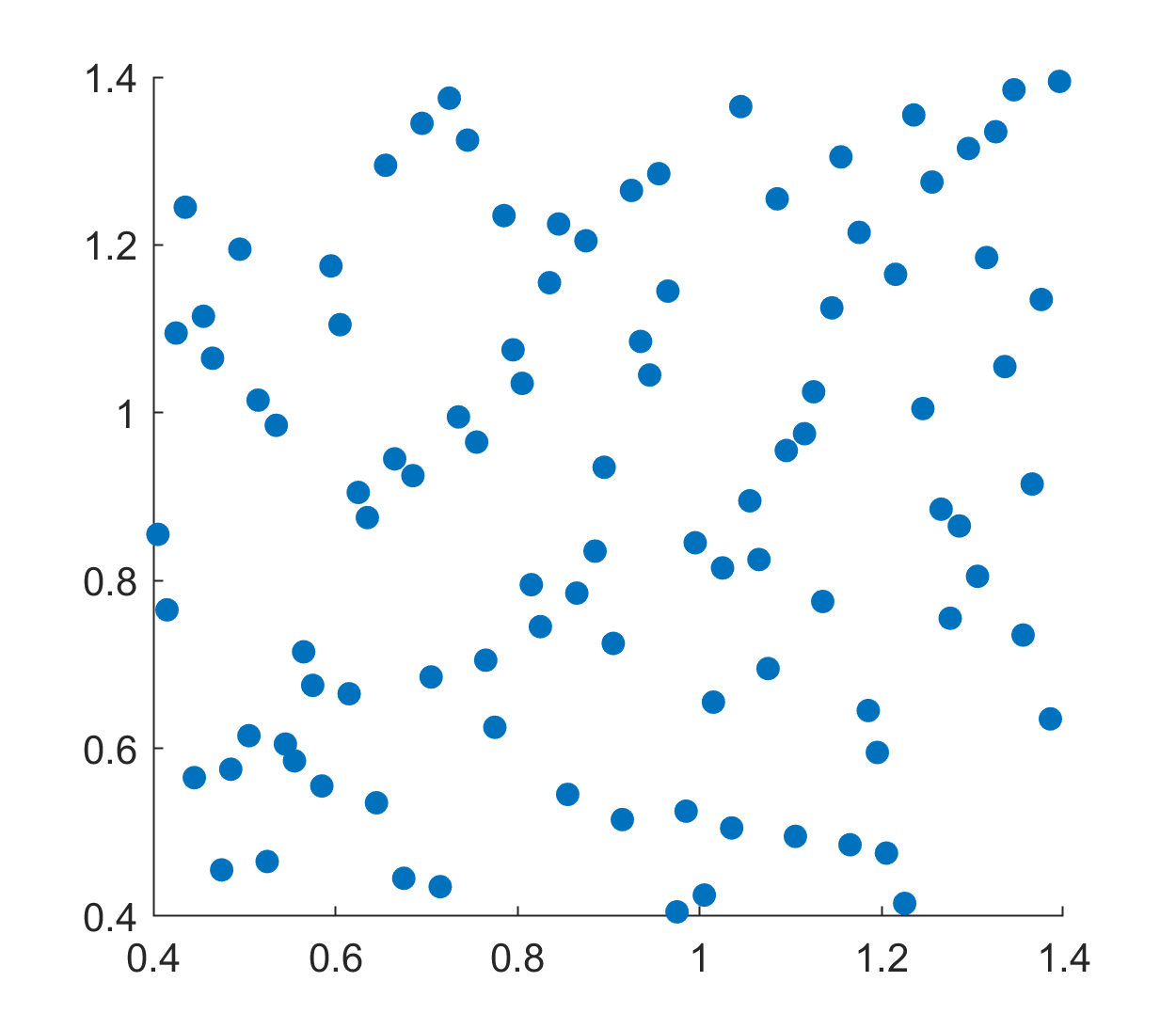}
    \caption{100 Latin hypercube Sample points}
    \label{lhs4}
  \end{subfigure}
  %  \begin{subfigure}{0.3\textwidth}
  %  \includegraphics[height=4cm,width=5cm]{sampling/plt_tst.png}
  %   \caption{30 test points}
  %   \label{30sample}
  % \end{subfigure}
   \begin{subfigure}{0.3\textwidth}
   \includegraphics[height=4cm,width=5cm]{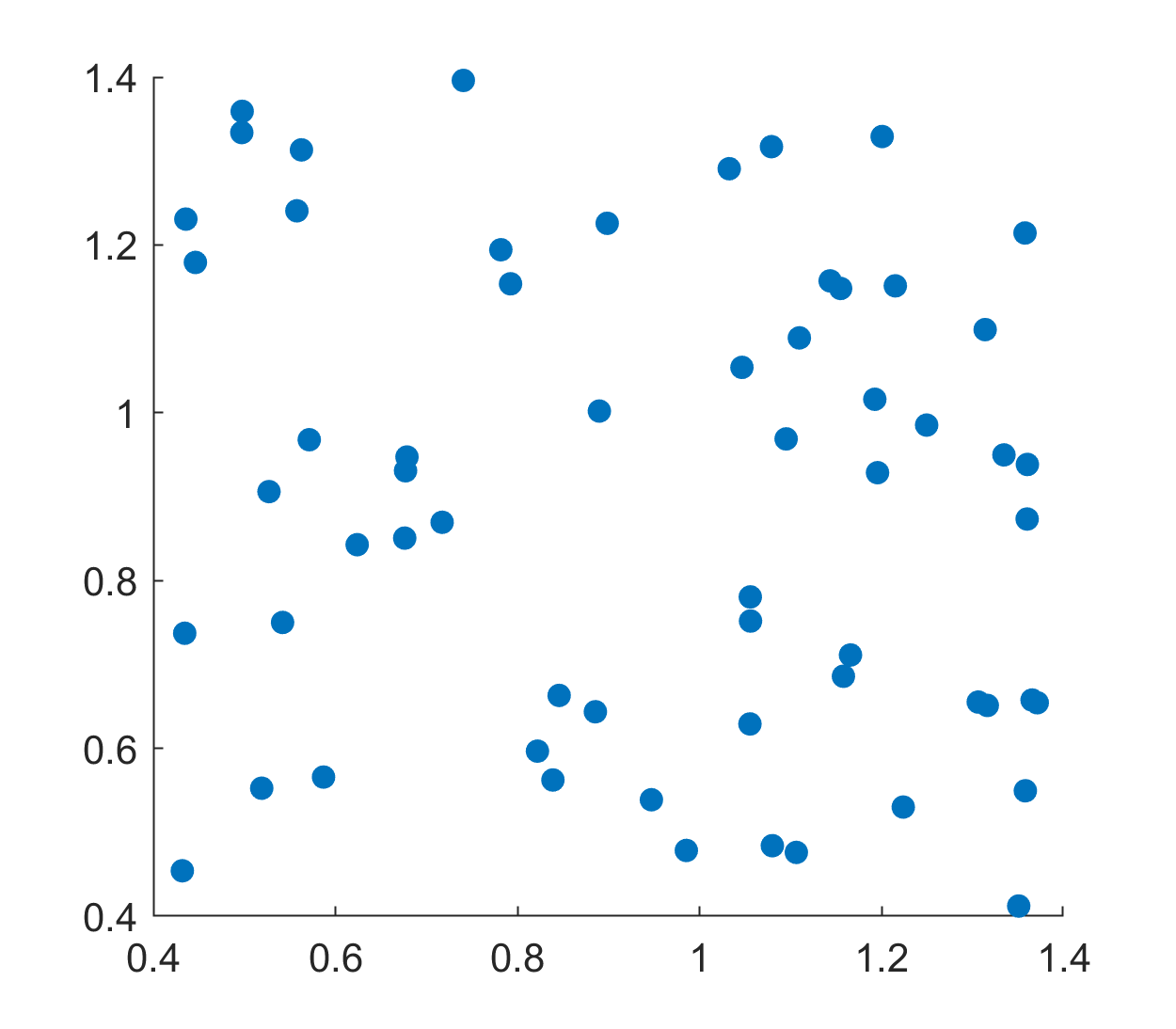}
    \caption{60 test points}
    \label{60sample}
  \end{subfigure}
  \caption{Sampling points and test points for EVP~\eqref{eq:par2}. }
  \label{sample:apr2}
\end{figure}

\begin{example}
Let us consider the eigenvalue problem
 \begin{equation}
\label{eq:par2}
\left\{
\begin{array}{ll}
     -\operatorname{div}(A(\bm{\mu})\nabla u(\bm{\mu}))=\lambda(\bm{\mu}) u(\bm{\mu})&  \textrm{ in }\Omega=(0,1)^2\\
     u(\bm{\mu})=0& \textrm{ on }\partial\Omega 
\end{array}
\right.
\end{equation}
where the diffusion $A(\bm{\mu})\in\mathbb R^{2\times 2}$ is given by the matrix
\begin{equation}\label{Amu}
    A(\bm{\mu})= 
    \begin{pmatrix}
        \frac{1}{\mu_1^2} &\frac{0.7}{\mu_2}\\
        \frac{0.7}{\mu_2} &\frac{1}{\mu_2^2}
    \end{pmatrix},
\end{equation}
with $\bm{\mu}=(\mu_1,\mu_2)\in \mathcal{M}\subset \mathbb{R}^2.$
The problem is symmetric and the parameter space $\mathcal{M}$ is chosen such that the matrix is positive definite. In particular, the matrix is positive definite for any nonzero value of $\mu_2$ and $\mu_1\in (-1.42,1.42)\setminus\{0\}$. In this example, we choose the parameter space to be $\mathcal{M}=[0.4,1.4]^2$. The surfaces for the first four eigenvalues are plotted separately and shown in Fig.~\ref{pltev:par2}. These eigensurfaces have the same regularity. It is thus sufficient to examine how GPR performs in one such curve.

% The eigenvlaue surfaces corresponding to the first and second are separate while the 3rd and fourth eigenvalue surfaces intersect with each other. 

To do so, we have sampled the parameter domain in two ways: 1) uniformly by selecting 64 equidistant sample points, and 2) using Latin-hypercube sampling to select $100$ sample points. These choices are shown in Fig.~\ref{sample:apr2}(a)-(b), respectively. To calculate the \re{RMSRE}, we use $60$ randomly selected test data inputs. The purpose is to check that the results do not depend on the choice of the test set. These are are respectively plotted in Fig.~\ref{sample:apr2}(c).
\end{example}

\begin{table}
 \centering
  \begin{tabular}{|c|c|c|c|c|c|c|} 
 Cov & Method & $\mu=(0.6,0.7)$ &$\mu=(0.6,1.1)$  & $\mu=(1.2,0.7)$  & $\mu=(1.2,1.1)$ &RMSRE  \\
 \hline
    & FEM &45.85197267& 34.64251756& 23.29248791& 12.18646344&   \\
 \hline
Exp      &GPR  &45.87852981& 34.67454252& 23.25188197& 12.15100067 & {$3.0\times 10^{-3}$}\\
& Rel. Err &$5.7\times 10^{-4}$&$9.2\times 10^{-4}$&$1.7\times 10^{-3}$&$2.9\times 10^{-3}$& \\%$4.7 \times 10^{-3}$\\
\hline
Mat\'{e}rn 3/2   & GPR &45.46178077 &34.20214955 &23.35084060 &12.20653421& {$1.5\times 10^{-3}$}\\
& Rel. Err &$8.5\times 10^{-3}$&$1.2\times 10^{-2}$&$2.5\times 10^{-3}$&$1.6\times 10^{-3}$& \\% $2.3 \times 10^{-3}$\\
\hline
Mat\'{e}rn 5/2   & GPR &45.59262509& 34.41018148& 23.26889817& 12.19913780& {$7.3\times 10^{-4}$}\\
& Rel. Err &$5.7\times 10^{-3}$&$6.7\times 10^{-3}$&$1.0\times 10^{-3}$&$1.0\times 10^{-3}$& \\%$1.3 \times 10^{-3}$\\
 \hline
SE       & GPR & 45.50385107& 34.36852082 &23.08793874 &12.05507251& {$1.9\times 10^{-4}$} \\
& Rel. Err &$7.6\times 10^{-3}$&$7.9\times 10^{-3}$&$8.8\times 10^{-3}$&$1.0\times 10^{-2}$& \\ %$1.1 \times 10^{-3}$\\
 \end{tabular}
\caption{Comparison of 1st eigenvalues of EVP~\eqref{eq:par2} using GPR model with different covariance functions using 64 uniform training samples.}% with 49 sample points
\label{par2:ev1}
 \end{table}

\begin{table} %% RRSE
 \centering
  \begin{tabular}{|c|c|c|c|c|c|c|c|c|} 
 Eigenvalue & No of Samples & Exponential  &Mat\'{e}rn 3/2 &Mat\'{e}rn 5/2  & Squared Exp.\\
\hline
 First
 & 64 &$3.0\times 10^{-3}$ & $1.5\times 10^{-3}$ & $7.3\times 10^{-4}$ & $1.9\times 10^{-4}$\\
 & 144 &$1.2\times 10^{-3}$ & $4.2\times 10^{-4}$ & $1.6\times 10^{-4}$ & $1.1\times 10^{-5}$\\
 & 256&  $6.2\times 10^{-4}$ & $1.7\times 10^{-4}$ & $5.6\times 10^{-5}$ & $5.4\times 10^{-6}$\\
 & 400 &    $3.7\times 10^{-4}$ & $7.9\times 10^{-5}$ & $2.2\times 10^{-5}$ & $1.4\times 10^{-6}$\\
 \hline
Second
 & 64 &$2.2\times 10^{-3}$ & $1.1\times 10^{-3}$ & $8.8\times 10^{-4}$ & $2.6\times 10^{-3}$\\
 & 144 &  $9.2\times 10^{-4}$ & $3.7\times 10^{-4}$ & $3.2\times 10^{-4}$ & $1.2\times 10^{-3}$\\
 & 256 &  $4.9\times 10^{-4}$ & $1.7\times 10^{-4}$ & $1.3\times 10^{-4}$ & $5.7\times 10^{-4}$\\
 & 400 &  $2.9\times 10^{-4}$ & $8.1\times 10^{-5}$ & $4.3\times 10^{-5}$ & $2.2\times 10^{-4}$\\
  \hline
 Third
 & 64 &$2.7\times 10^{-3}$ & $2.0\times 10^{-3}$ & $2.3\times 10^{-3}$ & $3.7\times 10^{-3}$\\
 & 144 &  $1.4\times 10^{-3}$ & $7.7\times 10^{-4}$ & $7.9\times 10^{-4}$ & $3.5\times 10^{-3}$\\
 & 256&  $8.1\times 10^{-4}$ & $3.8\times 10^{-4}$ & $4.1\times 10^{-4}$ & $2.5\times 10^{-3}$\\
 & 400&  $4.7\times 10^{-4}$ & $1.9\times 10^{-4}$ & $1.9\times 10^{-4}$ & $1.2\times 10^{-3}$\\
 \end{tabular}
\caption{{The values of RMSRE for the first three eigenvalues of EVP~\eqref{eq:par2} using GPR model with different covariance functions.}}
\label{par2:ev}
 \end{table}

In Tab.~\ref{par2:ev1}, we report the GPR predicted eigenvalues for the first curve of the problem \eqref{eq:par2} using the same selection of kernels. To compare the results, we have also reported the relative error in each case. In particular, we present the predictions at four different selected sample test points.%To measure the performance of the GPR models, we have selected $30$ random points from the parameter space and calculated the FEM eigenvalues at these points.
In the last column, the RMSRE values between the FEM eigenvalues and the GPR-based eigenvalues corresponding to the $60$ points are reported. The RMSRE values indicate that the GPR with squared exponential kernel produces better predictions for the considered eigenvalue surface.
% As for the second surface, the Mat\'{e}rn 3/2 kernel outperforms the rest. Finally, the exponential kernel is the best choice for the third eigenvalue surface. 
%Similar results are produced in~\ref{par2:ev1}.%-Tab.~\ref{par2:ev3}
w%here we instead use $60$ random test data points. The same conclusions can thus be seen to apply here regardless of the size of the test data set.

% Although the value of the error in changing but follow the same trend that is the Mat\'{e}rn kernels are performing better. From these results, we are unable to conclude which kernel is better for this problem. So, we increase the number of training points and use $100$ Latin Hypercube Sampling (LHS) samples. Note that LHS is a technique to select sparse points. }

{Next, we generate a GPR model on a training set containing more points to predict the same eigenvalue surface. In particular, we select $64$, $144$, $256$, and $400$ uniform sample points from the parameter space. In Tab.~\ref{par2:ev}, we report the RMSRE values for all three eigenvalues. From the table, it is observed that the GPR with squared exponential kernel gives better results for the first eigensurface also when we increase the number sample points. However, for the second and third eigensurfaces, the GPR with Mat\'{e}rn 5/2 kernel gives better results. This is because the surfaces are twice differentiable but not infinitely differentiable. }

% We similarly use our GPR eigenvalues of the problem \eqref{eq:par2} using GPR with four different kernels. 
{ The predictions corresponding to 100 Latin hypercube training points are reported in Tab.~\ref{par2:ev_l}.  The RMSRE values indicate that the GPR model with a Squared exponential (SE) kernel also gives better results than the other three for the first eigenvalues and Mat\'{e}rn 5/2  kernel gives better results for the second and third eigensurfaces. }% This supports the claim that a smooth kernel produces better results if the solution is smooth even in higher dimensions.}

\begin{table}
 \centering
  \begin{tabular}{|c|c|c|c|c|c|c|} 
 Eigenvalue & Exp & Mat\'{e}rn 3/2 &Mat\'{e}rn 5/2  & Squared exp.  \\
 \hline
First& {$3.7\times 10^{-3}$} &{$1.1\times 10^{-3}$}& {$3.9\times 10^{-4}$} &{$4.0\times 10^{-5}$} \\
Second &$3.5\times 10^{-3}$ & $1.3\times 10^{-3}$ & $7.6\times 10^{-4}$ & $8.6\times 10^{-4}$\\
Third & $2.7\times 10^{-3}$ & $1.2\times 10^{-3}$ & $1.1 \times 10^{-3}$ & $3.2\times 10^{-3}$\\
 \end{tabular}
\caption{{The values of} RMSRE for the eigenvalues of EVP~\eqref{eq:par2} using GPR model with different covariance functions using 100 LHS training samples.}% with 49 sample points
\label{par2:ev_l}
 \end{table}

\section{Conclusion}
In this paper, we have presented some preliminary discussion on how to apply Bayesian methods commonly used in machine learning, to the approximation of parametric PDE eigenvalues problems.
After a brief description of GPR, we critically compare GPR and splines for the approximation of data and functions. We plan to extend this discussion in subsequent investigations.

The main core of our numerical experiments deals with the comparison of different covariance functions in the application of GPR to parametric eigenvalue problems.

A general rule is that the GPR with an absolute exponential kernel outperforms if the output data is not smooth, while the GPR with the squared exponential kernel is better if the data is very smooth; otherwise, all four GPRs give comparable results. In the applications we have in mind, due to the presence of regions where the parametric dependence is smooth, together with other regions where the parametric dependence is not continuous, it is not immediate to define \emph{the best} kernel. Further studies and analysis will be performed to design an optimal strategy, possibly taking advantage of adaptive procedures and suitable combinations of kernels, and to identify appropriate criteria for the selection of the kernels.

\begin{acknowledgement}
 This research was supported by the Competitive Research Grants Program CRG2020 ``Synthetic data-driven model reduction methods for modal analysis'' awarded by the King Abdullah University of Science and Technology (KAUST).
Daniele Boffi is a member of the INdAM Research group GNCS and his research is partially supported by IMATI/CNR and by PRIN/MIUR.
\end{acknowledgement}

% \begin{funding}
%   Please insert information concerning research grant support here (institution and grant number). Please provide for each funder the funder’s DOI according to https://doi.crossref.org/funderNames?mode=list.
% \end{funding}
%\bibliographystyle{...}
%\bibliography{...}

\bibliographystyle{plain}
\bibliography{myref}

\end{document}